\documentclass[a4paper,11pt]{article}

\usepackage{amsfonts}
\usepackage{latexsym}
\usepackage{epsfig}
\usepackage{amssymb}
\usepackage{graphicx}
\usepackage{oldgerm}
\usepackage{theorem}

\newcommand{\mib}[1]{\mbox{\boldmath $#1$}}

\def\b{\mib{b}}

\theorembodyfont{\itshape}
\newtheorem{thm}{Theorem}[section]
\newtheorem{lem}[thm]{Lemma}

\newtheorem{prop}[thm]{Proposition}

\newcommand{\qed}{\hbox{\rule[-2pt]{3pt}{6pt}}}

\def\Det{{\rm Det}}

\def\x{\mib{x}}
\def\y{\mib{y}}

\def\w{\mib{w}}
\def\a{\mib{a}}
\def\b{\mib{b}}
\def\n{\mib{n}}

\def\f{\mib{f}} 
\def\X{\mib{X}}
\def\B{\mib{B}}

\def\vtheta{\mib{\theta}}

\def\N{{\bf N}}
\def\R{{\bf R}}
\def\W{{\bf W}}
\def\E{{\bf E}}
\def\O{{\bf 0}}

\def\tS{\widetilde{S}}
\def\mpN{\mathfrak{p}_{N}}

\def\mbK{\mathbb{K}}
\def\bK{\mib{K}}
\def\sK{{\cal K}}

\def\Ai{{\rm Ai}}

\def\mX{\mathfrak{X}}

\def\sG{{\cal G}}

\def\cR{\rho}
\def\cD{{\rm D}}
\def\cG{{\det {\bf M}}}

\newcommand{\SSC}[1]{\section{#1}\setcounter{equation}{0}}

\begin{document}

\title{\bf Noncolliding Brownian Motion \\
and Determinantal Processes}
\author{
Makoto Katori
\footnote{
Department of Physics,
Faculty of Science and Engineering,
Chuo University, 
Kasuga, Bunkyo-ku, Tokyo 112-8551, Japan;
e-mail: katori@phys.chuo-u.ac.jp
}
and 
Hideki Tanemura
\footnote{
Department of Mathematics and Informatics,
Faculty of Science, Chiba University, 
1-33 Yayoi-cho, Inage-ku, Chiba 263-8522, Japan;
e-mail: tanemura@math.s.chiba-u.ac.jp
}}
\date{12 November 2007}
\pagestyle{plain}
\maketitle
\begin{abstract}
A system of one-dimensional Brownian motions (BMs)
conditioned never to collide with each other
is realized as (i) Dyson's BM model,
which is a process of eigenvalues of hermitian
matrix-valued diffusion process in the 
Gaussian unitary ensemble (GUE),
and as (ii) the $h$-transform of absorbing BM
in a Weyl chamber, where the harmonic function $h$
is the product of differences 
of variables (the Vandermonde determinant).
The Karlin-McGregor formula gives determinantal
expression to the transition probability density
of absorbing BM. We show from the
Karlin-McGregor formula, if the initial state
is in the eigenvalue distribution of GUE,
the noncolliding BM is a determinantal process,
in the sense that any multitime correlation function is
given by a determinant specified by a matrix-kernel.
By taking appropriate scaling limits,
spatially homogeneous and inhomogeneous
infinite determinantal processes are derived.
We note that the determinantal processes related with
noncolliding particle systems have a feature in common
such that the matrix-kernels are expressed
using spectral projections of
appropriate effective Hamiltonians.
On the common structure of matrix-kernels,
continuity of processes in time
is proved and general property of the
determinantal processes is discussed.
\end{abstract}

\SSC{INTRODUCTION}

In the present paper, we discuss noncolliding Brownian 
motion (BM) in one-dimension with finite number $N$ of 
particles and its infinite particle limits $N \to \infty$.
The condition imposed to $N$ particles in the model,
{\it not to collide with each other}, causes ``entropy forces"
between all pairs of particles, which are repulsive 
long-ranged interactions proportional to the inverse
of distances between particles.
When we draw sample paths of $N$ particles of the system
on the spatio-temporal plane,
random sets of nonintersecting $N$ paths are obtained.
Viewing them as random patterns of polymers
or phase boundaries on a plane, 
the present system has been used as a model of
polymer networks \cite{deG68,EG95},
or a model showing wetting (or melting) transitions \cite{Fis84}
in statistical physics; see also \cite{K00,TY03}.
Recently many authors have reported that notion of
noncolliding BM and its discrete counterpart 
called {\it vicious walk} \cite{Fis84} is very useful
to analyze the polynuclear growth models 
\cite{PS02,Joh03,SI04,IS05},
time-dependent correlations of quantum spin chains \cite{Bog06},
traffic problems \cite{BBDS06},
and the Chern-Simons theory \cite{deH05,AHG06}.

\subsection{Observations of Three-Dimensional Bessel Process}

In order to demonstrate the important
connection between the random matrix theory \cite{Meh04}
and the noncolliding BM
here we show a couple of observations of
the three-dimensional Bessel process.
The noncolliding BM can be regarded as a multivariate
generalization of the three-dimensional Bessel process
given below.

Let $B_1(t), B_2(t), B_3(t)$ be
one-dimensional standard BMs
(see Section 2.1 for definition).
They are assumed to be independent and we consider
a $2 \times 2$ traceless hermitian matrix
\begin{equation}
M^{(1)}(t)=\left( \matrix{
B_1(t) & B_2(t)+{\rm i} B_3(t) \cr
B_2(t)-{\rm i} B_3(t) & -B_1(t) } \right),
\label{eqn:M1}
\end{equation}
where ${\rm i}=\sqrt{-1}$.
Since the four entries $(M^{(1)}_{jk}(t))_{1 \leq j, k \leq 2}$
are BMs, $M^{(1)}(t), t \in [0, \infty)$ is regarded
as a matrix-valued process, which describe a diffusion 
process in the space of $2 \times 2$ traceless
hermitian matrices, which is identified with
the three-dimensional real space
$\R^3$ ($\R$ denotes the set of all real numbers).
At each time $t \in [0, \infty)$, it will be diagonalized by
an appropriate unitary matrix and the eigenvalue
is given by $\pm X(t)$ with
\begin{equation}
X(t)=\sqrt{(B_{1}(t))^2+(B_{2}(t))^2+(B_{3}(t))^2}.
\label{eqn:Bessel1}
\end{equation}
If we consider a Brownian particle in $\R^3$,
$\B(t)=(B_1(t), B_2(t), B_3(t)), t \in [0, \infty)$,
the distance of the particle from the origin
({\it i.e.,} the radial coordinate of
$\B(t)$ ) is given by (\ref{eqn:Bessel1})
and thus it equals the eigenvalue process $X(t)$
associated with the matrix-valued process
$M^{(1)}(t)$.
This is called the three-dimensional Bessel process
in probability theory (see, {\it e.g.} 
\cite{KS91,BS02}), and a simple application
of the It\^o formula gives
its stochastic differential equation (SDE) as
\begin{equation}
dX(t)= d \widetilde{B}(t)+\frac{dt}{X(t)},
\quad t \in [0, \infty), 
\quad X(0)=x > 0,
\label{eqn:Bessel2}
\end{equation}
where $\widetilde{B}(t)$ is another one-dimensional
standard BM than the above $B_{j}(t), j=1,2,3$.
Corresponding to the SDE (\ref{eqn:Bessel2}),
the backward Kolmogorov (Fokker-Planck)
equation for the transition probability density
$p_{X}(t,y|x)$, starting from
$x >0$ at time $t=0$ and arriving at
$y > 0$ at time $t >0$, is given by
\begin{equation}
\frac{\partial}{\partial t}
p_X(t,y|x)=\frac{1}{2} 
\frac{\partial^2}{\partial x^2}
p_X(t, y|x)
+\frac{1}{x} \frac{\partial}{\partial x}
p_X(t,y|x),
\label{eqn:Bessel3}
\end{equation}
and its solution with $\lim_{t \to 0} 
p_{X}(t,y|x)=\delta(x-y)$ is obtained as
\begin{equation}
p_{X}(t,y|x)=\frac{h(y)}{h(x)} 
p_{\rm abs}(t,y|x)
\label{eqn:h-trans0}
\end{equation}
with
$$
p_{\rm abs}(t, y|x)=
p(t, y|x)-p(t, y|-x),
$$
where $h(x)=x$ and
$p(t,y|x)$ is the heat kernel given by (\ref{eqn:pA1})
in Section 2.1.
The reflection principle of BM can be used to prove that
$p_{\rm abs}(t, y|x)$ is the transition probability
density from $x>0$ to $y>0$ during time $t$ 
of the absorbing BM, in which an absorbing
wall is set at the origin and any particle
is absorbed if it arrives at the wall.
It is a matter of course that $h(x)=x$ is harmonic,
$\partial^2 h(x)/\partial x^2=0$,
but we dare to say that $p_{X}$
is the harmonic transform
($h$-transform) of $p_{\rm abs}$
looking at (\ref{eqn:h-trans0}).
We introduce another $2 \times 2$ matrix
\begin{equation}
M^{(2)}(t, \y|\x)
=\left( \matrix{
p(t,y_1|x_1) & p(t,y_1|x_2) \cr
p(t,y_2|x_1) & p(t,y_2|x_2) } \right)
\label{eqn:M2}
\end{equation}
for $\x=(x_1,x_2), \y=(y_1,y_2) \in \R^2$,
and consider its determinant
\begin{equation}
f_2(t, \y|\x) = \det \Big[ M^{(2)}(t, \y|\x) \Big].
\label{eqn:KM0}
\end{equation}
Then it is easy to see that
\begin{equation}
p_{\rm abs}(t, y|x)
=\sqrt{\frac{\pi t}{2}}
f_2 ( t/2, \{- y/2, y/2 \}
| \{ -x/2, x/2 \}),
\quad x, y >0.
\label{eqn:abs0}
\end{equation}
In summary the three-dimensional Bessel process
$X(t), t \in [0, \infty)$ 
has two different realizations;
(i) the eigenvalue-process of $2 \times 2$
hermitian-matrix valued process (\ref{eqn:M1}),
and (ii) the $h$-transform of the absorbing
BM with a wall at $x=0$.
We also observed that the transition probability
density of the absorbing BM has a determinantal 
expression of a $2 \times 2$ matrix (\ref{eqn:M2})-(\ref{eqn:abs0}).
If we consider the two-dimensional BM, 
it is represented by motion of
a point in the two-dimensional space
$(x_1, x_2) \in \R^2$.
We put an absorbing boundary on a line $x_2=x_1$ 
and trace the motion of the point
in the region 
$\W_2=\{\x=(x_1,x_2) \in \R^2 : x_1 < x_2\}$.
The transition probability density
from $\x=(x_1, x_2) \in \W_2$ to
$\y=(y_1, y_2) \in \W_2$ is generally given
by the determinant (\ref{eqn:KM0}).
As a special case of it (with a time-change $t \to t/2$),
(\ref{eqn:abs0}) is given.

\subsection{Dyson's BM Model, 
Karlin-McGregor Formula and Noncolliding BM}

The noncolliding BM, 
$\X(t)=(X_1(t), X_2(t), \cdots, X_N(t)), t \in [0, \infty)$
is a conditional diffusion process.
It has the following two kinds of realizations.
\begin{description}
\item{(i)} \quad
In order to generate the random matrix ensembles,
Dyson introduced $N \times N$ matrix-valued
diffusion processes \cite{Dys62}.
For the Gaussian unitary ensemble (GUE),
$N^2$ independent one-dimensional BMs are used to assign entries
of matrix to satisfy the condition that
the matrix is hermitian at any time.
Eigenvalues are real and define 
an $N$-particle system in one dimension
called Dyson's BM model
(with the parameter $\beta=2$ corresponding to GUE).
This process solves the SDE
(see \cite{KT04} for the proof
using generalized Bru's theorem)
\begin{equation}
dX_j(t)=dB_j(t) 
+\sum_{1 \leq k \leq N : k \not=j} 
\frac{dt}{X_{j}(t)-X_{k}(t)}, 
\quad 1 \leq j \leq N, \quad
t \in [0, \infty),
\label{eqn:Dyson1}
\end{equation}
where $\B(t)=(B_1(t), \cdots, B_{N}(t))$ is
an $N$-dimensional BM; 
$B_{j}(t), 1 \leq j \leq N$, are 
independent one-dimensional standard BMs.
The SDE (\ref{eqn:Dyson1}) is an $N$-variable
generalization of the SDE for the three-dimensional
Bessel process (\ref{eqn:Bessel2}).
It was proved that with probability one
Dyson's BM is non-colliding
\cite{RS93}.
\item{(ii)} \quad
We consider the following subset of $\R^{N}$,
\begin{equation}
\W_N=\Big\{\x =(x_1, \cdots, x_N) \in \R^{N} :
x_1 < x_2 < \cdots < x_N \Big\}.
\label{eqn:Weyl1}
\end{equation}
It is called the Weyl chamber of type $A_{N-1}$
\cite{Gra99,KT04}.
The absorbing BM is defined by putting absorbing walls
at all boundaries of the region $\W_N$, whose transition
probability density, 
when the process starts from $\x \in \W_N$ 
at time $t=0$ and arrives at $\y \in \W_N$
at time $t > 0$, is given by 
\begin{equation}
f_N(t, \y|\x)=
\det_{1 \leq j, k \leq N}
\Big[ p(t, y_j|x_{k}) \Big],
\quad \x, \y \in \W_N.
\label{eqn:f1}
\end{equation}
This determinantal expression is
known as the Karlin-McGregor formula
\cite{KM59}.
(Such a determinantal formula for nonintersecting paths
is known as the Lindstr\"on-Gessel-Viennot
formula \cite{Lin73,GV85} in the enumerative combinatorics,
see also 
\cite{deG68,Fis84,Ste90,EG95,KGV00,Joh02,KT03a,OR03,KT05,
Kra06,EK06}.)
The noncolliding BM is given by the 
$h$-transform of the absorbing BM \cite{Gra99}
\begin{equation}
p_N(t,\y|\x)=
\frac{h_N(\y)}{h_N(\x)} f_N(t, \y|\x),
\label{eqn:h-trans0b}
\end{equation}
where 
\begin{equation}
h_N(\x) = \prod_{1 \leq j< k \leq N}(x_k-x_j)
= \det_{1 \leq j, k \leq N}
[x_{j}^{k-1}]
\label{eqn:Vandermonde1}
\end{equation}
is the Vandermonde determinant \cite{Mac95,Ful97,Sta99}.
\end{description}

We note that there appear three different kinds of
matrices.
The matrices representing Dyson's matrix-valued
process (eq.(\ref{eqn:M1}) for the simplest case),
matrices in the Karlin-McGregor determinants
(eqs.(\ref{eqn:KM0}) and (\ref{eqn:f1})) 
and that in the Vandermonde determinant 
(\ref{eqn:Vandermonde1}).
Of course, the equivalence of Dyson's BM model
(the eigenvalue process of the first kind of matrices)
with the noncolliding BM implies
direct connection between the random matrix theory
and stochastic processes.
In the present paper, however, we will show that
the Karlin-McGregor formula is much more important.
In Section 3 we will show that the Vandermonde
determinant appears in the Schur
function expansion of the Karlin-McGregor
determinant. With the combination of
these two determinants, the orthogonal polynomial
method is applicable to study the processes.
This method has been developed to analyze
multi-matrix models in the random matrix theory \cite{Meh04}.

\subsection{Matrix-Kernels and Determinantal Processes}

In Section 2, we will explain that the Hermite 
orthonormal functions are useful to represent BMs.
Precise descriptions of facts briefly mentioned above
will be given in Section 3.
Staring from Karlin-McGregor's determinantal expression
of transition probability density,
we will prove in Section 4 that, 
if we specify the initial configuration as
the GUE-eigenvalue distribution, 
the generating function of multitime correlation functions
is given by a Fredholm determinant for
the noncolliding BM, and thus multitime correlation
functions are generally given by determinants
(Theorem \ref{thm:finite}).
The system whose spatial correlations are given by
determinants is usually called a determinantal point field
(or a determinantal point process)
in probability theory \cite{Sos00,ST03}. 
Theorem \ref{thm:finite} states
that the noncolliding BM is not only a determinantal
point field at any fixed time $0 \leq t < \infty$,
but it is also a determinantal point field
on the spatio-temporal plane.
We say that the noncolliding BM is a
(finite) {\it determinantal process} to express this 
situation.

Determinantal processes are generally determined by their
matrix-kernels (see, for example, \cite{TW04}).
The matrix-kernel of our finite noncolliding BM is
expressed by the Hermite orthonormal functions 
$\{\varphi_k\}_{k=0}^{\infty}$, 
which is called the extended Hermite kernel in \cite{TW04}.
In Section 5 we will show that the asymptotic properties
of $\{\varphi_k\}_{k=0}^{\infty}$ 
completely determine infinite particle limits
of the matrix-kernel
and then the determinantal process.
Appropriate scaling limits
are performed and two infinite particle systems are
derived from the noncolliding BM.
One of them is the spatially homogeneous 
infinite determinantal process with matrix-kernel
expressed by trigonometric functions
(Theorem \ref{thm:infinite_sin}) and
another is the spatially inhomogeneous 
one with matrix-kernel expressed by Airy functions
(Theorem \ref{thm:infinite_Ai}).
The former kernel is called the extended sine kernel
and the latter the extended Airy kernel in \cite{TW04}.
We will claim in Section 6 that these three
determinantal processes (one finite and two
infinite systems) and others reported in 
references \cite{PS02,Joh03,TW03,TW04,KNT04,KT07} 
have a common structure;
the matrix-kernels are expressed by
spectral projections associated with 
appropriate self-adjoint operators (effective Hamiltonians)
\cite{PS02}.
As explained by Spohn \cite{Spo87} and 
by Pr\"ahofer and Spohn \cite{PS02},
this common feature is shared also with the
1+1 dimensional Fermi field in quantum mechanics
(see also \cite{FS03,HO07}).
It may be due to the similarity 
between the Karlin-McGregor formula
for noncolliding systems
and the Slater determinant 
for free fermion systems with the 
Fermi exclusion principle \cite{Baym69}.
For finite and infinite
determinantal processes with matrix-kernels
associated with spectral projections, 
we will prove that the determinantal processes
are continuous in time
(Lemma \ref{lem:continuity}) and 
discuss the bilinear forms derived from 
correlation functions 
(Proposition \ref{thm:bilinear}).
Future problems are given in Section 8.

\SSC{BROWNIAN MOTION AND HERMITE POLYNOMIALS}
\subsection{Diffusion Equation and 
Hamiltonian of Harmonic Oscillator}

Let $(\Omega, {\cal F}, P)$ be the probability space.
One-dimensional standard BM
starting from a point $x_0 \in \R$ is defined as 
a real-valued stochastic process
$\{B(t,\omega) : t \in [0, \infty) \}$,
which satisfies the following conditions
(see, for example, \cite{KS91,GS92}).
Here $\omega$ is a label on sample path, 
$\omega \in \Omega$.
\begin{description}
\item{1.} \quad
$B(0, \omega)=x_0$ with probability 1.
\item{2.} \quad
For any fixed $\omega \in \Omega$, 
$B(t)$ is a continuous real-function of
$t$ with probability 1.
(With this property we say that the paths are
continuous in time.)
\item{3.} \quad
For any series of times
$t_0 \equiv 0 < t_1 < \cdots < t_M, M=1,2, \cdots$,
$\{B(t_{m+1})-B(t_{m}) \}_{m=0,1, \cdots, M-1}$
are independent and they are normally distributed
with mean 0 and variance $t_{m+1}-t_{m}$.
\end{description}
Then if we introduce an integral kernel
(Gaussian kernel)
\begin{equation}
p(t,x|x')=\frac{1}{\sqrt{2 \pi t}}
\exp \left\{ -\frac{(x-x')^2}{2t} \right\},
\quad t > 0, x, x' \in \R,
\label{eqn:pA1}
\end{equation}
the probability that the BM stays
in an interval $[a_m, b_m]$,
$-\infty < a_m < b_m < \infty$, at each time $t_m$,
$m=1,2, \cdots, M$, is given by
\begin{eqnarray}
&& P \Big(
B(t_m) \in [a_m, b_m], m=1,2, \cdots, M \Big)
\nonumber\\
&& \quad =
\int_{a_1}^{b_1} dx_1 \int_{a_2}^{b_2} dx_2
\cdots \int_{a_M}^{b_M} dx_M
\prod_{m=0}^{M-1} 
p(t_{m+1}-t_{m}, x_{m+1} | x_{m}).
\nonumber
\end{eqnarray}
That is, the transition
probability density of the BM
is given by (\ref{eqn:pA1}).
Since it satisfies the one-dimensional
diffusion equation (heat equation)
\begin{equation}
\frac{\partial}{\partial t} u(t,x)
= \frac{1}{2} \frac{\partial^2}{\partial x^2} u(t,x)
\label{eqn:diff1}
\end{equation}
with the initial condition
$u(0,x)=\delta(x-x')$, it is specially
called the heat kernel.

We consider the following transformation
of variables, 
$(t, x) \mapsto (\tau, \zeta)$ ;
\begin{equation}
\tau = \tau(t)=\log t, \quad
\zeta = \zeta(t,x)= \frac{x}{\sqrt{2t}}
\label{eqn:trans1}
\end{equation}
and set
\begin{equation}
u(t,x)=e^{-(\tau+\zeta^2)/2} U(\tau, \zeta).
\label{eqn:trans2}
\end{equation}
Then the diffusion equation (\ref{eqn:diff1})
is transformed to
\begin{equation}
\frac{\partial}{\partial \tau} U(\tau, \zeta)
= - \frac{1}{2} \left(
{\cal H}_{\rm H}-\frac{1}{2} \right) U(\tau, \zeta),
\label{eqn:diff2}
\end{equation}
where
\begin{equation}
{\cal H}_{\rm H}=-\frac{1}{2} \frac{\partial^2}{\partial \zeta^2}
+ \frac{1}{2} \zeta^2.
\label{eqn:H1}
\end{equation}
Note that ${\cal H}_{\rm H}$ is identified with the
Hamiltonian in the coordinate($\zeta)$-representation
of the one-dimensional harmonic oscillator 
in quantum mechanics,
if we set the mass of the oscillator $m=1$,
the circular frequency $\omega=1$, and $\hbar=1$.

Following Dirac's description
\cite{Dirac58,Baym69}, we consider the
real-valued Hilbert space
with basis $\{|\zeta \rangle : \zeta \in \R\}$,
which is orthonormal
$\langle \zeta | \zeta' \rangle 
=\delta(\zeta-\zeta')$
and complete
\begin{equation}
\int_{\R} d \zeta \,
|\zeta \rangle \langle \zeta | = 1.
\label{eqn:complete1}
\end{equation}
Let $\widehat{\cal H}_{\rm H}$ be the operator
such that
$
\langle \zeta' | \widehat{\cal H}_{\rm H} |\zeta \rangle
= \delta(\zeta'-\zeta) {\cal H}_{\rm H}
$
with (\ref{eqn:H1}).
We consider a state vector
$|\Psi(\tau) \rangle$, which follows the equation 
\begin{equation}
\frac{\partial}{\partial \tau}
|\Psi(\tau) \rangle
= - \frac{1}{2} 
\left( \widehat{\cal H}_{\rm H}-\frac{1}{2} \right)
|\Psi(\tau) \rangle.
\label{eqn:S1}
\end{equation}
If we multiply $\langle \zeta |$ to (\ref{eqn:S1})
from the left, we will have
the equation (\ref{eqn:diff2}) with
$U(\tau, \zeta)=\langle \zeta | \Psi(\tau) \rangle$.
Given $|\Psi(\tau') \rangle, \tau' \in (-\infty, \infty)$,
the solution of (\ref{eqn:S1})
for $\tau \geq \tau'$ is obtained as
$$
|\Psi(\tau) \rangle
= \exp \left\{
-\frac{1}{2}\left(\widehat{\cal H}_{\rm H}-\frac{1}{2} \right)
(\tau-\tau') \right\} |\Psi(\tau') \rangle.
$$
Then
\begin{eqnarray}
U(\tau, \zeta) 
&=& \langle \zeta |
\exp \left\{
-\frac{1}{2}\left(\widehat{\cal H}_{\rm H}-\frac{1}{2} \right)
(\tau-\tau') \right\} |\Psi(\tau') \rangle
\nonumber\\
&=& \int_{-\infty}^{\infty} d \zeta' \,
\langle \zeta | 
\exp \left\{
-\frac{1}{2}\left(\widehat{\cal H}_{\rm H}-\frac{1}{2} \right)
(\tau-\tau') \right\} | \zeta' \rangle
U(\tau', \zeta'),
\nonumber
\end{eqnarray}
where (\ref{eqn:complete1}) was used.
Insert this result into (\ref{eqn:trans2}), we have
\begin{eqnarray}
u(t, x) &=& \int_{-\infty}^{\infty} d \zeta' \,
e^{-(\tau+\zeta^2)/2+((\tau')+(\zeta')^2)/2}
\langle \zeta | 
\exp \left\{
-\frac{1}{2}\left(\widehat{\cal H}_{\rm H}-\frac{1}{2} \right)
(\tau-\tau') \right\} | \zeta' \rangle u(t', x')
\nonumber\\
&=& \int_{-\infty}^{\infty} dx' \,
p(t-t', x|x') u(t', x'),
\nonumber
\end{eqnarray}
where
\begin{eqnarray}
p(t-t', x|x')
&=& \frac{1}{\sqrt{2}}
e^{-\tau/2-\zeta^2/2+(\zeta')^2/2} 
\langle \zeta |
\exp \left\{
-\frac{1}{2}\left(\widehat{\cal H}_{\rm H}-\frac{1}{2} \right)
(\tau-\tau') \right\} | \zeta' \rangle
\nonumber\\
&=& \frac{1}{\sqrt{2t}} 
e^{-x^2/4t+(x')^2/4t'}
\left\langle \frac{x}{\sqrt{2t}} \right|
\exp \left\{
-\frac{1}{2}\left(\widehat{\cal H}_{\rm H}-\frac{1}{2} \right)
(\tau-\tau') \right\} 
\left| \frac{x'}{\sqrt{2t'}} \right\rangle. \qquad
\label{eqn:pA2}
\end{eqnarray}
This is the transition probability density
previously given as (\ref{eqn:pA1}).

\subsection{Hermite Polynomials and Equalities}

Let $\N=\{1,2,3, \cdots\}$ and 
$\N_0 = \N \cup \{0\}$.
The eigenvalues of
$\widehat{\cal H}_{\rm H}$ are 
$n+1/2$ with $n \in \N_{0}$; 
$
\widehat{\cal H}_{\rm H} |n \rangle
= \left(n+ 1/2 \right) | n \rangle.
$
Here $\{ |n \rangle : n \in \N_0 \}$
denotes the set of eigenvectors of $\widehat{\cal H}_{\rm H}$, 
which is orthonormal
$\langle n | n' \rangle = \delta_{n,n'}$
and complete
$ \sum_{n=0}^{\infty} |n \rangle \langle n |=1$.
Let $\{H_n(x) : n \in \N_0\}$ be the Hermite polynomials
\begin{eqnarray}
H_n(x) &=& e^{x^2} \left(-\frac{d}{dx}\right)^n e^{-x^2}
\nonumber\\
&=& n! \sum_{k=1}^{[n/2]}
(-1)^{n} \frac{(2x)^{n-2k}}
{k ! (n-2k)!},
\label{eqn:Hermite0}
\end{eqnarray}
where $[z]$ denotes the largest number
not greater than $z$.
They are orthogonal
$$
\int_{\R} dx \, e^{-x^2} H_n(x) H_m(x)
=h_n \delta_{n,m}, \quad n, m \in \N_0
$$
with 
$h_n=\sqrt{\pi} 2^{n} n!.$
We set
\begin{equation}
\varphi_n(\zeta)= 
\frac{1}{\sqrt{h_n}} e^{-\zeta^2/2} H_n (\zeta).
\label{eqn:Hermite1}
\end{equation}
Then $\{\varphi_{n}(\zeta) : n \in \N_0\}$ are
orthonormal and
$$
\langle \zeta | n \rangle
= \langle n | \zeta \rangle
= \varphi_n(\zeta)
$$
is established.

Now we define
\begin{equation}
\widehat{\cal H}_{\varphi} = 
\frac{1}{2} \left( 
\widehat{\cal H}_{\rm H}-\frac{1}{2} \right).
\label{eqn:Hphi1}
\end{equation}
Then
\begin{equation}
\widehat{\cal H}_{\varphi} |n \rangle
= \frac{n}{2} | n \rangle, \quad
n \in \N_0
\label{eqn:Hphi2}
\end{equation}
and (\ref{eqn:pA2}) gives
\begin{eqnarray}
p(t-t', x|x') 
&=& \frac{1}{\sqrt{2}}
e^{-\tau/2-\zeta^2/2+(\zeta')^2/2}
\sum_{n=0}^{\infty}
\langle \zeta | 
e^{-\widehat{\cal H}_{\varphi}\tau} 
|n \rangle \langle n |
e^{\widehat{\cal H}_{\varphi}\tau'}
|\zeta' \rangle
\nonumber\\
&=& \frac{1}{\sqrt{2}}
e^{-\tau/2-\zeta^2/2+(\zeta')^2/2}
\sum_{n=0}^{\infty}
\langle \zeta |e^{-n\tau/2}| n \rangle
\langle n | e^{n\tau'/2} |\zeta' \rangle \nonumber\\
&=& \frac{1}{\sqrt{2}}
e^{-\tau/2-\zeta^2/2+(\zeta')^2/2}
\sum_{n=0}^{\infty} e^{-n(\tau-\tau')/2}
\varphi_{n}(\zeta) \varphi_{n}(\zeta'),
\nonumber
\end{eqnarray}
that is
\begin{equation}
p(t-t',x|x')
=\frac{1}{\sqrt{2t}}
e^{-x^2/4t+(x')^2/4t'}
\sum_{n=0}^{\infty}
\left(\frac{t'}{t}\right)^{n/2}
\varphi_n \left(\frac{x}{\sqrt{2t}}\right)
\varphi_n \left(\frac{x'}{\sqrt{2t'}}\right).
\label{eqn:Mehler}
\end{equation}
This expression for the heat kernel
(\ref{eqn:pA1}) is called Mehler's formula
\cite{Bat53,Sze75}.

We introduce the vectors
\begin{eqnarray}
|t,x \rangle &=& e^{\zeta^2/2} 
e^{\widehat{\cal H}_{\varphi}\tau} |\zeta \rangle, 
\nonumber\\
\langle t, x | &=&
\frac{1}{\sqrt{2}} e^{-(\tau+\zeta^2)/2}
\langle \zeta | e^{-\widehat{\cal H}_{\varphi}\tau}.
\label{eqn:Vectors1}
\end{eqnarray}
It is easy to see that orthonornality and completeness
of $\{|\zeta \rangle : \zeta \in \R \}$ imply
$\langle t, x | t, x' \rangle =\delta(x-x')$
and
\begin{equation}
\int_{\R} dx \,
|t, x \rangle \langle t,x |=1.
\label{eqn:Vectors3}
\end{equation}
By (\ref{eqn:Vectors3}), we have the
equalities for
$0 \leq t_1 \leq t_2 \leq t_3$
\begin{eqnarray}
&& \int_{\R} d x_2 \langle t_3, x_3|t_2, x_2 \rangle
\langle t_2, x_2| t_1, x_1 \rangle
= \langle t_3, x_3|t_1, x_1 \rangle, 
\quad x_1, x_3 \in \R, \nonumber\\
&& \int_{\R} dx_1 \langle t_2, x_2|t_1, x_1 \rangle
\langle t_1, x_1|n \rangle
=\langle t_2, x_2| n \rangle, \quad n \in \N_0, 
\, x_2 \in \R, \nonumber\\
&& \int_{\R} dx_1 \int_{\R} dx_2 
\langle n|t_2, x_2 \rangle \langle t_2, x_2|t_1, x_1 \rangle
\langle t_1, x_1|n'\rangle
=\langle n| n' \rangle
=\delta_{n, n'}, \, n, n' \in \N_0. \nonumber
\end{eqnarray}
Though these seem to be trivial,
if we note that (\ref{eqn:pA2}) is rewritten as
\begin{equation}
p(t-t', x|x')=
\langle t,x | t', x' \rangle, \quad
0 \leq t' \leq t, \, x, x' \in \R,
\label{eqn:pA3}
\end{equation}
they become meaningful; for
$0 \leq t_1 \leq t_2 \leq t_3$ 
\begin{eqnarray}
\label{eqn:CK2}
&& \int_{-\infty}^{\infty} dx_2 \,
p(t_3-t_2, x_3|x_2) p(t_2-t_1, x_2|x_1)
= p(t_3-t_1, x_3|x_1), \quad
x_1, x_3 \in \R, \\
\label{eqn:invariance1b}
&& \int_{-\infty}^{\infty} dx_1 \, p(t_2-t_1, x_2|x_1)
\phi_n (t_1, x_1) =
\phi_n (t_2, x_2),
\quad n \in \N_0, x_2 \in \R, \\
\label{eqn:invariance2b}
&& \int_{-\infty}^{\infty} dx_1 \,
\int_{-\infty}^{\infty} dx_2 \,
\widehat{\phi}_n (t_2, x_2)
p(t_2-t_1, x_2|x_1) \phi_{n'} (t_1, x_1)
=\delta_{n,n'}, \, n, n' \in \N_{0},
\end{eqnarray} 
where
\begin{eqnarray}
\label{eqn:phin}
\phi_n(t,x) &\equiv& \langle t,x|n \rangle=
\frac{1}{\sqrt{2}} t^{-(n+1)/2}
e^{-x^2/4t} \varphi_{n}\left(
\frac{x}{\sqrt{2t}}\right), \\
\label{eqn:hatphin}
\widehat{\phi}_{n}(t,x) &\equiv& 
\langle n | t,x \rangle 
= t^{n/2}e^{x^2/4t} \varphi_n
\left( \frac{x}{\sqrt{2t}} \right),
\quad n \in \N_0.
\end{eqnarray}
The equation (\ref{eqn:CK2}) is called the
Chapman-Kolmogorov equation.
The equalities (\ref{eqn:invariance1b}) mean 
invariance of the functions
$\phi_{n}, n \in \N_0$, with respect to
the heat kernel.

\SSC{NONCOLLIDING SYSTMES}
\subsection{Application of Karlin-McGregor Formula}

Here we introduce the noncolliding BM
in a finite time-period $(0, T)$ with $T > 0$.
It is defined as an $N$-particle system of 
one-dimensional standard BMs 
conditioned not to collide with each other
in $(0,T)$.
As mentioned in Section 1.2, the
transition probability density of
the absorbing BM in the Weyl chamber
$\W_N$ of type $A_{N-1}$ is
given by (\ref{eqn:f1}) as an application
of Karlin-McGregor formula \cite{KM59}.
The probability that $\B(t)$ starting from
$\x' \in \W_N$ stays in $\W_N$ up to
at least time $t >0$ is given by
\begin{equation}
{\cal N}_N(t, \x')
=\int_{\W_N} f_N(t, \x|\x') d \x,
\quad \x' \in \W_N.
\label{eqn:NN1}
\end{equation}
The transition probability density 
of noncolliding BM is then given by
\begin{equation}
g_{N,T}(t', \x'; t, \x)=
\frac{{\cal N}_{N}(T-t, \x)}{{\cal N}_{N}(T-t', \x')}
f_N(t-t', \x|\x')
\label{eqn:gN1}
\end{equation}
for $0 < t' \leq t < T$,
$\x, \x' \in \W_N$
\cite{KT02,KT03a,KT03b}.
It should be noted that this process is in general
temporally inhomogeneous.
In the following, we will consider the $T \to \infty$ limit
to make the process be homogeneous in time.

\subsection{Schur Function Expansion}

By multilinearity of determinant
(\ref{eqn:f1}) with (\ref{eqn:pA1}),
$$
f_N(t, \x|\x')=
\left(\frac{1}{2\pi t}\right)^{N/2}
e^{-(|\x|^2+|\x'|^2)/2t}
\det_{1 \leq j, k \leq N} 
[e^{x_j x_k'/t}],
$$
where
$|\x|^2 \equiv \sum_{j=1}^{N}x_j^2$.
Consider $F_N(\x, \y)=\det_{1 \leq j, k \leq N}
[e^{x_j y_k}]$ for a pair of
multivariates $\x=\{x_j\}_{j=1}^{N}$ and
$\y=\{y_j\}_{j=1}^{N}$.
By definition of determinant,
$F(\x,\y)$ is skew-symmetric under any
exchange of indices of $\{x_j\}$ and also
it is for $\{y_j\}$;
Let ${\cal S}_N$ be the symmetric group of $N$ variables
(the set of all permutations of $N$ variables)
and for $\sigma \in {\cal S}_N$
write
$\sigma(\x)=(x_{\sigma(1)}, \cdots, x_{\sigma(N)})$.
Then for any $\sigma \in {\cal S}_{N}$
$F_{N}(\sigma(\x),\y)=F_{N}(\x, \sigma(\y))
={\rm sgn}(\sigma)F_N(\x, \y)$.
A fundamental skew-symmetric polynomials of multivariate
$\x$ is given by a product of differences
(the Vandermonde determinant) (\ref{eqn:Vandermonde1}).
The quotient of
$F_{N}(\x, \y)$ divided by
$h_N(\x) h_N(\y)$ is a symmetric function
both of $\x$ and $\y$.
The following lemma shows an
expansion of the symmetric part using
the Schur functions, which are labeled by partitions
$\mu=(\mu_1, \mu_2, \cdots)$,
sets of nonnegative integers
in decreasing order 
$\mu_1 \geq \mu_2 \geq \cdots$, and defined by
\begin{equation}
s_{\mu}(\x) = \frac{\displaystyle{
\det_{1 \leq j, k \leq N} [x_j^{\mu_{k}+N-k}] } }
{\displaystyle{\det_{1 \leq j, k \leq N} [x_j^{N-k}] } },
\label{eqn:Schur1}
\end{equation}
The non-zero $\mu_j$'s in a partition $\mu$ are called
parts of $\mu$ and the number of parts
is called length of $\mu$ and
denoted by $\ell(\mu)$
\cite{Mac95,Ful97,Sta99}.
The Schur function expansion is a special case of
character expansions (see
\cite{Bal00,Bal02,KO01,KT04,HO07}).
Let $\Gamma(a)=\int_{0}^{\infty} e^{-y} y^{a-1} dy$
for $a >0$ (the Gamma function) and note
$\Gamma(a+1)=a!$ if $a \in \N_{0}$.

\begin{lem}
\label{thm:Schur_exp1}
\begin{eqnarray}
\det_{1 \leq j, k \leq N}
[e^{x_j y_k}]
 &=& h_N(\x) h_N(\y) \sum_{\mu: \ell(\mu) \leq N}
\frac{s_{\mu}(\x) s_{\mu}(\y)}
{\prod_{k=1}^{N} \Gamma(\mu_k+N-k+1)}
\nonumber\\
&=& \frac{h_N(\x) h_N(\y)}
{\prod_{k=1}^{N} \Gamma(k)}
\times \left\{ 1+ {\cal O}(|\x|) \right\}
\quad \mbox{in} \quad |\x| \to 0.
\label{eqn:Schur_exp1}
\end{eqnarray}
\end{lem}
\noindent{\it Proof.} \quad
By multilinearity of determinant, we have
\begin{eqnarray}
\det_{1 \leq j, k \leq N} [e^{x_j y_k}] 
&=&
\det_{1 \leq j, k \leq N} \left[
\sum_{n=0}^{\infty} \frac{(x_{i}y_{k})^n}
{\Gamma(n+1)} \right] \nonumber\\
&=&
\sum_{\n=(n_{1}, n_{2}, \cdots, n_{N}) \in \N_{0}^{N}}
\prod_{m=1}^{N} \frac{1}{\Gamma(n_{m}+1)}
\det_{1 \leq j, k \leq N} \Bigg[ (x_{j} y_{k})^{n_{j}}
\Bigg].
\label{eqn:eqdet1}
\end{eqnarray}
We can see that for any symmetric function $f(\n)$ of 
$\n=(n_{1}, \cdots, n_{N})
\in \N_{0}^{N}$,
$$
\sum_{\n \in \N_{0}^{N}} f(\n)
\det_{1 \leq j, k \leq N} \Bigg[ (x_{j}y_{k})^{n_{j}} \Bigg]
= \sum_{\n \in \N_{0}^{N}} f(\n)
\frac{1}{N!} \sum_{\sigma \in {\cal S}_{N}} 
\det_{1 \leq j, k \leq N} 
\Bigg[ (x_{j} y_{k})^{n_{\sigma(j)}} \Bigg],
$$
and
$$
\sum_{\sigma \in {\cal S}_{N}} 
\det_{1 \leq j, k \leq N} 
\Bigg[ (x_{j} y_{k})^{n_{\sigma(j)}} \Bigg]
= \det_{1 \leq j, k \leq N} [x_{j}^{n_{k}}]
\det_{1 \leq \ell, m \leq N} [y_{\ell}^{n_{m}}].
$$
Since $\det_{1 \leq j, k \leq N} [x_{j}^{n_{k}}]=0$
if $n_{k_{1}}=n_{k_{2}}$ for any pair
$1 \leq k_{1} \not= k_{2} \leq N$,
(\ref{eqn:eqdet1}) equals
$$
\sum_{0 \leq n_{1} < n_{2} < \cdots < n_{N}}
\frac{\displaystyle{\det_{1 \leq j, k \leq N} (x_{j}^{n_{k}})
\det_{1 \leq \ell, m \leq N} (y_{\ell}^{n_{m}})}}
{\prod_{j=1}^{N} \Gamma(n_{j}+1)} .
$$
Here we change the variables in summation from
$\{n_{j}\}$ to $\{\mu_{j}\}$ by
$\mu_{j}=n_{j}-N+j, 1 \leq j \leq N$.
Using (\ref{eqn:Schur1}) we obtain the first equation
of (\ref{eqn:Schur_exp1}).
Since $s_{\mu}(\O)=0$ unless 
$\mu=\O \equiv (0,0, \cdots, 0) \in \N_0^{N}$,
and $s_{0}(\O)=1$, 
the estimation in $|\x| \to 0$ is given as shown by
the second equation of (\ref{eqn:Schur_exp1}).
\qed

By this lemma, we have the estimate
\begin{equation}
f_{N}(t, \x|\x')=\frac{1}{C_N}
t^{-N^2/2} h_N(\x) h_N(\x')
e^{-|\x|^2/2t}
\times \left\{ 1+ {\cal O}
\left( \frac{|\x'|}{\sqrt{t}} \right) \right\}
\quad \mbox{in} \quad
\frac{|\x'|}{\sqrt{t}} \to 0
\label{eqn:estimate1}
\end{equation}
with 
$C_N=(2\pi)^{N/2} \prod_{j=1}^{N} \Gamma(j)$.
The integral formula
$$
\int_{\R^N} e^{-a|\x|^2} |h_N(\x)|^{2 \gamma} d \x
=(2\pi)^{N/2} (2a)^{-N(\gamma(N-1)+1)/2}
\prod_{j=1}^{N} \frac{\Gamma(1+j \gamma)}{\Gamma(1+\gamma)}
$$
is found in \cite{Meh04} (eq.(17.6.7) p.321)
as a variation of the Selberg integral \cite{Sel44},
whose proof was given in \cite{Mac82}.
If we set $\gamma=1/2, a=1/2t$ and note that
the integral over $\R^{N}$ can be replaced by
the integral over $\W_N$ multiplied by $N!$, we have
\begin{equation}
\int_{\W_N} e^{-|\x|^2/2t} h_{N}(\x) d \x
= C_N' t^{N(N+1)/4}
\label{eqn:I1}
\end{equation}
with 
$C_N'=2^{N/2} \prod_{j=1}^{N} \Gamma(j/2)$.
Similarly by setting $\gamma=1$ and $a=1/2t$, we have
\begin{equation}
\int_{\W_N} e^{-|\x|^2/2t} (h_N(\x))^2 d \x
= C_N t^{N^2/2}.
\label{eqn:I2}
\end{equation}
Using (\ref{eqn:I1}) with (\ref{eqn:estimate1}), 
we obtain the asymptotics of ${\cal N}_{N}$,
\begin{equation}
{\cal N}_{N}(t, \x')=
\frac{C_N'}{C_N} t^{-N(N-1)/4} h_{N}(\x')
\times \left\{ 1+
{\cal O}\left( \frac{|\x'|}{\sqrt{t}} \right) \right\}
\quad \mbox{in} \quad
\frac{|\x'|}{\sqrt{t}} \to 0.
\label{eqn:estimate2}
\end{equation}
The integral (\ref{eqn:I2}) will be used shortly.

\subsection{Temporally Homogeneous Limit}

By the above estimate (\ref{eqn:estimate2}), we can
take the $T \to \infty$ limit in (\ref{eqn:gN1})
and obtain the transition probability density, which is
homogeneous in time,
{\it i.e.,} a function of time difference $t-t'$,
\begin{eqnarray}
p_{N}(t-t', \x|\x')
&=& \lim_{T \to \infty} g_{N,T}(t', \x'; t, \x)
\nonumber\\
&=& h_{N}(\x) f_{N}(t-t', \x|\x') 
\frac{1}{h_{N}(\x')}.
\label{eqn:pN1}
\end{eqnarray}
From now on, we consider the noncolliding BM,
which is defined by this transition probability density.
It is a temporally homogeneous process
and we denote it by
$\X(t)=(X_1(t),X_2(t), \cdots, X_N(t))$.

\vskip 0.5cm
\noindent{\bf Remark 1.} \quad
The product of differences (the Vandermonde determinant)
$h_N(\x)$
given by (\ref{eqn:Vandermonde1}) 
is a harmonic function in the sense
$$
\nabla^2 h_N(\x)= 
\sum_{j=1}^{N} \frac{\partial^2}{\partial x_j^2}
h_{N}(\x)=0,
$$
which has strictly positive values at
interior points of $\W_{N}$ and
zero at the boundary.
Eq. (\ref{eqn:pN1}) is considered as a transformation
from $f_N$ to $p_N$ associated with the
harmonic function, which is called 
the $h$-transform \cite{Doo84}.
That is, the temporally homogeneous noncolliding
BM is an $h$-transform
by $h_N$ of the absorbing BM
in the Weyl chamber $\W_N$ \cite{Gra99}.
It is easy to confirm that 
$p_N(t, \, \cdot \,|\x)$ 
satisfies the following backward
Kolmogorov equation
$$
\frac{\partial}{\partial t} u(t, \x)
= \frac{1}{2} \nabla^2 u(t, \x)
+ \nabla \log h_{N}(\x) \cdot
\nabla u(t, \x),
$$
which is a multivariate extension of (\ref{eqn:Bessel3}).
It implies that the process $\X(t)$
is a diffusion process, which solves the SDE
(\ref{eqn:Dyson1}) of 
Dyson's BM model
(with $\beta=2$), which describes the
time-evolution of eigenvalues of 
hermitian matrix-valued (GUE) diffusion process
\cite{Dys62,KT04}.
This equivalence between Dyson's BM model for GUE
(the eigenvalue part of the matrix-valued BM)
and the present noncolliding BM
(BM conditioned not to collide)
is a multi-dimensional extension of the
equivalence between the three-dimensional Bessel process
(the radial coordinate of the three-dimensional BM)
and the one-dimensional BM
conditioned to stay in the positive region
$B(t) >0$, as announced in Section 1. 
See also \cite{KT03a,KTNK03,KT04}.
\vskip 0.5cm

Let $\nu_{0}(\x)$ be the probability density at 
time $t_0 >0$.
Then the probability density of 
distribution $\X(t), t \geq t_0$
is given by
\begin{eqnarray}
\nu_{t}(\x) &=& \int_{\W_{N}} 
p_N(t-t_0, \x|\x') \nu_{0}(\x') d \x'
\nonumber\\
&=& h_{N}(\x) \int_{\W_N} f_{N}(t-t_0, \x|\x')
\frac{\nu_{0}(\x')}{h_{N}(\x')} d \x'.
\label{eqn:nut1}
\end{eqnarray}
By the estimate (\ref{eqn:estimate1}), we will see that
\begin{equation}
\nu_{t}(\x) \simeq
\frac{1}{C_{N}} t^{-N^2/2} e^{-|\x|^2/2t}
(h_{N}(\x))^2 \quad
\mbox{in} \quad t \to \infty,
\label{eqn:nuinf1}
\end{equation}
if the distribution $\nu_0$ has finite moment.
Note that the integral formula (\ref{eqn:I2})
guarantees that (\ref{eqn:nuinf1}) is normalized.
It should be noted that the distribution (\ref{eqn:nuinf1})
is equal to the eigenvalue distribution of
hermitian random-matrices in GUE
with variance $\sigma^2=t$ \cite{Meh04}.

\begin{prop}
\label{thm:GUE}
For any initial distribution having finite moment,
the asymptote of probability density of distribution
of the noncolliding BM,
$\X(t)$, in $t \to \infty$ is expressed by
the eigenvalue distribution of GUE with
variance $t$.
\end{prop}

\SSC{DETERMINANTAL PROCESS}
\subsection{GUE Initial Distribution}

Consider a sequence of times,
$0 < t_0 < t_1 < \cdots < t_M,
M=1,2, \cdots$ for observations
of distribution of $\X(t)$.
Given the initial distribution $\nu_0$
at time $t_0$, multitime probability density
is given using (\ref{eqn:pN1}) as
\begin{eqnarray}
p_{N}(t_0, \x^{(0)}; t_1, \x^{(1)}; 
\cdots ; t_{M}, \x^{(M)})
&=& \prod_{m=0}^{M-1} 
p_{N}(t_{m+1}-t_{m}, \x^{(m+1)} | \x^{(m)}) \nu_{0}(\x^{(0)})
\nonumber\\
&=& h_{N}(\x^{(M)}) 
\prod_{m=0}^{M-1} f_{N}(t_{m+1}-t_{m}, \x^{(m+1)}|
\x^{(m)}) 
\frac{\nu_0(\x^{(0)})}{h_{N}(\x^{(0)})}.
\qquad
\label{eqn:pNmulti1}
\end{eqnarray}

Now we assume that $\nu_{0}$ is the GUE-eigenvalue
distribution with variance $t_0$,
\begin{equation}
\nu_0(\x^{(0)})
=\frac{1}{C_{N}} t_0^{-N^2/2}
e^{-|\x^{(0)}|^2/2t_0}
(h_N(\x^{(0)}))^2.
\label{eqn:nu01}
\end{equation}
Then (\ref{eqn:pNmulti1}) becomes the
product of $f_N$'s multiplied by a factor
\begin{equation}
\frac{1}{C_{N}} t_{0}^{-N^2/2}
e^{-|\x^{(0)}|^2/2t_0}
h_{N}(\x^{(M)}) h_{N}(\x^{(0)}).
\label{eqn:factor1}
\end{equation}
By multilinearity of determinant and the fact
that the coefficient of the highest order
term of the Hermite polynomial
$H_n(x)$ is $(2x)^{n}$ (see (\ref{eqn:Hermite0})),
the following equality is established,
$$
\det_{1 \leq j, k \leq N} 
\Big[ H_{j-1}(x_{k}) \Big]
=h_{N}(2 \x).
$$
Then if we define the multivariate functions
\begin{eqnarray}
\mu_{0}(\x) &=& \det_{1 \leq j, k \leq N}
\Big[ \phi_{j-1}(t_0, x_k) \Big], \nonumber\\
\mu_{M}(\x) &=& \det_{1 \leq j, k \leq N}
\Big[ \widehat{\phi}_{j-1}(t_M, x_k) \Big],
\label{eqn:in_fi}
\end{eqnarray}
where $\phi_n$ and $\widehat{\phi}_n$ are
given by (\ref{eqn:phin}) and (\ref{eqn:hatphin}),
respectively,
the factor (\ref{eqn:factor1}) is readily shown 
to be equal to 
$\mu_0(\x^{(0)}) \mu_{M}(\x^{(M)})$.
That is, the multitime probability density
(\ref{eqn:pNmulti1}) 
is written as
$$
\mu_{M}(\x^{(M)})
\prod_{m=0}^{M-1} f_N(t_{m+1}-t_{m}, \x^{(m+1)}|\x^{(m)})
\mu_{0}(\x^{(0)}).
$$

This expression shows that the probability law is
invariant under any permutation
of indices of particles ;
$x_{j}^{(m)} \mapsto x_{\sigma(j)}^{(m)}, 0 \leq m \leq M,
\sigma \in {\cal S}_{N}$.
Then description will be easier if
we regard that particles are identical and
indistinguishable.
We denote by $\mX$ the space of countable subsets $\xi$ of $\R$
satisfying $\sharp (\xi \cap K) < \infty$
for any compact subset $K$.
The space $\mX$ is a Polish space with
the vague topology.
For $\x = (x_1, x_2, \dots, x_n)\in \bigcup_{\ell=1}^{\infty}
\W_\ell$,
we denote $\{x_j\}_{j=1}^n \in \mX$ simply by $\{\x\}$.
Then we consider the process on the set $\mX$, 
$\Xi^{\bf X}_N (t) = \{\X(t)\}, t \in [t_0, \infty)$, 
such that,
for any $M+1$ times 
$t_0 < t_1 < \cdots < t_{M-1} < t_M < \infty$ \, ($M=0,1,2, \cdots$)
the multitime probability density
is given of the form
\begin{eqnarray}
&& \mpN \Big(t_0, \{\x^{(0)}\}; t_1, \{\x^{(1)}\}; \cdots ;
t_{M-1}, \{\x^{(M-1)}\}; t_{M}, \{\x^{(M)}\} \Big)
\nonumber\\
&=&
\mu_{M}(\{\x^{(M)}\})
\prod_{m=0}^{M-1} 
f_N(t_{m+1}-t_{m}, \x^{(m+1)}|\x^{(m)})
\mu_{0}(\{\x^{(0)}\}),
\label{eqn:mpN1}
\end{eqnarray}

For $\x^{(m)} \in \R^N$, $0 \leq m \leq M$,
and $N'=1,2,\dots, N$, we put 
$\x^{(m)}_{N'} = \left(x_1^{(m)}, x_2^{(m)}, \dots, x_{N'}^{(m)}\right)$.
For a sequence $\{N_m \}_{m=0}^{M}$ of positive integers 
less than or equal to $N$,
we define the $(N_0, \dots, N_{M})$-multitime 
correlation function by
\begin{eqnarray}
&& \rho_N \Big(t_{0}, \{\x^{(0)}_{N_0}\}; t_{1}, \{ \x^{(1)}_{N_1}\} ; 
\dots; t_M, \{ \x^{(M)}_{N_{M}}\} \Big) 
\nonumber\\
&&=
\int_{\prod_{m=0}^{M} \R^{N-N_{m}}}
\mpN \Big( t_{0}, \{ \x^{(0)}\}; 
\dots; t_M, \{ \x^{(M)} \} \Big)
\prod_{m=0}^{M}
\frac{1}{(N-N_{m})!}
\prod_{j=N_{m}+1}^{N} dx_{j}^{(m)}.
\label{def:corr}
\end{eqnarray}
Expectations with respect to the configurations
$\{\X(t_0)\}, \{\X(t_1)\},\dots, \{\X(t_M)\}$ 
are denoted by $\E_{N}$ :
\begin{eqnarray}
&&\E_{N} 
\Big[ f(\{\X(t_0)\}, \{\X(t_1)\},\dots, \{\X(t_M)\}) \Big]
\nonumber\\
&=&
\left( \frac{1}{N!} \right)^{M+1} \int_{\R^{N(M+1)}}
f(\{\x^{(0)} \}, \dots, \{\x^{(M)}\}) 
\mpN \Big(t_{0}, \{\x^{(0)}\}; 
\dots; t_M, \{\x^{(M)}\} \Big)
\prod_{m=0}^{M} \prod_{j=1}^{N} d x_{j}^{(m)}.
\nonumber\\
\label{eqn:expect}
\end{eqnarray}

\noindent{\bf Remark 2.} \quad
Set $t'=0, t=t_0 >0, \x=\x^{(0)}$ in
(\ref{eqn:pN1}) and consider the limit
$\x' \to \O$.
By (\ref{eqn:estimate1}), we can see
\begin{equation}
\lim_{|\x'|\to 0}
p_{N}(t_0, \x^{(0)}|\x')=
\nu_0(\x^{(0)}).
\label{eqn:zerostart}
\end{equation}
This fact is essentially the same thing with
that stated as Proposition \ref{thm:GUE} 
for the scaling property of the process.
The GUE-eigenvalue distribution (\ref{eqn:nu01}) 
adopted here as the initial distribution 
at time $t_0 >0$ is immediately
realized if we start the present noncolliding
system at time 0 from $\{\O\}$.
The state $\{\O\}$, which is a boundary of $\W_N$,
is entrance
\cite{KT02,KT03a,KT03b}.
\subsection{Generating Function}

Let $C_{0}(\R)$ be the set of all continuous real functions
with compact supports. 
For $\f=(f_{0}, f_{1}, \cdots, f_{M}) \in C_{0}(\R)^{M}$,
and $\vtheta=(\theta_{0}, \theta_{1}, \cdots, \theta_{M})
\in \R^{M}$,
the generating function for multitime correlation functions
is defined
for the process $\{\X(t)\}, t \in [0,T]$ as
\begin{equation}
\Psi_{N} ( \f; \vtheta)
= \E_{N} \left[
\exp \left\{ \sum_{m=0}^{M}
\theta_{m} \sum_{j_{m}=1}^{N} 
f_{m}(X_{j_{m}}(t_{m})) \right\} \right].
\label{eqn:Psi0}
\end{equation}
Let 
$$
\chi_{m}(x)=e^{\theta_{m} f_{m}(x)}-1, 
\qquad 0 \leq m \leq M,
$$
and write (\ref{eqn:Psi0}) as 
$\Psi_{N}[\chi]$.
Then by the definition of
multitime correlation function (\ref{def:corr}), we have 
\begin{eqnarray}
\Psi_{N}[\chi]
&=& \sum_{N_{0}=0}^{N} \sum_{N_{1}=0}^{N} \cdots
\sum_{N_{M}=0}^{N}
\prod_{m=0}^{M}\frac{1}{N_m !}
\int_{\R^{N_{0}}} \prod_{j=1}^{N_0} d x_{j}^{(0)}
\int_{\R^{N_{1}}} \prod_{j=1}^{N_1} d x_{j}^{(1)}
 \cdots
\int_{\R^{N_{M}}} 
\prod_{j=1}^{N_{M}} d x_{j}^{(M)} \nonumber\\
&& \times \prod_{m=0}^{M} \prod_{j=1}^{N_{m}} 
\chi_{m} \Big(x_{j}^{(m)} \Big)
\rho_{N} \left(t_0, \{\x_{N_{0}}^{(0)} \}; 
t_{1}, \{\x_{N_{1}}^{(1)} \};
\dots ; t_{M}, \{\x_{N_{M}}^{(M)}\} \right).
\qquad
\label{eqn:kPhiA1}
\end{eqnarray}

By the definition (\ref{eqn:Psi0}) with (\ref{eqn:expect})
and (\ref{eqn:mpN1}), 
we have
\begin{eqnarray}
\Psi_{N}[\chi]
&=& \left(\frac{1}{N !}\right)^{M+1}
\int_{\R^{N(M+1)}} \prod_{m=0}^{M} \prod_{j=1}^{N} 
d x_{j}^{(m)} \, 
\det_{1 \leq j, k \leq N} \Big[
\widehat{\phi}_{j-1}(t_M, x_k^{(M)})(1+\chi_{M}(x_k^{(M)})) \Big]
\nonumber\\
&& \times
\prod_{m=0}^{M-1} \det_{1 \leq j, k \leq N}
\left[ p(t_{m+1}-t_{m}, x_{k}^{(m+1)}| x_{j}^{(m)})
(1+\chi_{m}(x_{k}^{(m)})) \right]
\det_{1 \leq j, k \leq N}
\Big[ \phi_{j-1}(t_0, x_k^{(0)}) \Big].
\nonumber
\end{eqnarray}
By definition of determinant, it is easy to prove
the following identity 
for square integrable continuous functions
$g_{j}, \bar{g}_{j}, 1 \leq j \leq N$,
\begin{equation}
\frac{1}{N!}
\int_{\R^{N}} 
\det_{1 \leq j, k \leq N} \Bigg[g_{j}(x_{k})\Bigg]
\det_{1 \leq j, k \leq N} \Bigg[ \bar{g}_{j}(x_{k}) \Bigg]
\prod_{j=1}^{N} d x_{j}
= \det_{1 \leq j, k \leq N} \left[
\int_{\R} g_{j}(x) \bar{g}_{k}(x) dx \right],
\label{eqn:Heine}
\end{equation}
which is called the Heine identity.
By repeated applications of this identity
we have
\begin{equation}
\Psi_{N}[\chi]= \det_{1 \leq j, k \leq N}
\Big[F_{jk}[\chi] \Big]
\label{eqn:F1}
\end{equation}
with
\begin{eqnarray}
F_{jk}[\chi] &=&
\int_{\R^{M+1}} \prod_{m=0}^{M} dx^{(m)} \,
\Big\{ \widehat{\phi}_{j-1}(t_M, x^{(M)})(1+\chi_{M}(x^{(M)})) \Big\} 
\nonumber\\
&& \times
\prod_{m=0}^{M-1} \Big\{
p(t_{m+1}-t_{m}, x^{(m+1)}| x^{(m)} ) 
(1+\chi_{m}(x^{(m)}))
\Big\}
\phi_{k-1}(t_0, x^{(0)}).
\nonumber
\end{eqnarray}
By the Chapman-Kolmogorov equation (\ref{eqn:CK2})
and the invariance (\ref{eqn:invariance2b}),
$F_{jk}[0]=\delta_{jk}$ and then
$\Psi_{N}[0]=1$, which implies that (\ref{eqn:mpN1})
is indeed normalized.
If we use the notations introduced in Section 2, 
it is written as
\begin{eqnarray}
F_{jk}[\chi]
&=& \int_{\R^{M}} \prod_{m=0}^{M} dx^{(m)} \,
\langle j-1 | t_{M}, x^{(M)} \rangle
\{1+\chi_M (x^{(M)}) \} \nonumber\\
&& \qquad \times \prod_{m=0}^{M-1}
\Big[ \langle t_{m+1}, x^{(m+1)}| t_m, x^{(m)} \rangle
\{1+\chi_m (x^{(m)})\} \Big]
\langle t_0, x^{(0)}|k-1 \rangle
\nonumber\\
&=& \langle j-1|k-1 \rangle \nonumber\\
&+& \sum_{\ell=0}^{M+1}
\sum_{M \geq m_1 > \cdots > m_{\ell} \geq 0}
\int_{\R^{\ell}} \prod_{n=1}^{\ell} dx^{(m_{n})}
\langle j-1 | t_{m_1}, x^{(m_1)} \rangle \chi_{m_1}(x^{(m_1)})
\langle t_{m_1}, x^{(m_1)}|t_{m_2}, x^{(m_2)} \rangle
\nonumber\\
&& \qquad \times \cdots \times
\langle t_{m_{\ell-1}}, x^{(m_{\ell-1})}|
t_{m_{\ell}}, x^{(m_{\ell})} \rangle
\chi_{m_{\ell}}(x^{(m_{\ell})})
\langle t_{m_{\ell}}, x^{(m_{\ell})} | k-1 \rangle.
\label{eqn:F3}
\end{eqnarray}
Now we introduce an indicator $\hat{1}_{+}$ such that
\begin{equation}
\langle t_m, x|\hat{1}_{+}|t_n, y \rangle =
\left\{
   \begin{array}{ll}
      \langle t_m, x| t_n, y \rangle & \quad \mbox{if} \quad m > n \\
      0 & \quad \mbox{otherwise}.
   \end{array}\right. 
\label{eqn:1+}
\end{equation}
Then (\ref{eqn:F3}) can be written of the form of
an infinite series
\begin{eqnarray}
F_{jk}[\chi] 
&=& \langle j-1|k-1 \rangle \nonumber\\
&+& \sum_{\ell=0}^{\infty}
\sum_{m_1=0}^{M}\cdots 
\sum_{m_{\ell}=0}^{M} 
\int_{\R^{\ell}} \prod_{n=1}^{\ell} dx^{(m_n)}
\langle j-1 | t_{m_1}, x^{(m_1)} \rangle \chi_{m_1}(x^{(m_1)})
\langle t_{m_1}, x^{(m_1)}|\hat{1}_+ |t_{m_2}, x^{(m_2)} \rangle
\nonumber\\
&& \qquad
\times \cdots \times
\langle t_{m_{\ell-1}}, x^{(m_{\ell-1})}| \hat{1}_{+} |
t_{m_{\ell}}, x^{(m_{\ell})} \rangle
\chi_{m_{\ell}}(x^{(m_{\ell})})
\langle t_{m_{\ell}}, x^{(m_{\ell})} | k-1 \rangle,
\nonumber
\end{eqnarray}
in which all terms with $\ell > M+1$
are zero by (\ref{eqn:1+}).
Let
\begin{equation}
\widehat{\chi}_{+}^{m,n}(x,y)
=\langle t_m, x| \hat{1}_{+}|t_n, y \rangle
\chi_{n}(y),
\quad m,n \in \{0,1, \cdots, M\},
\, x, y \in \R,
\label{eqn:hatchi1}
\end{equation}
and define an $M \times M$ matrix-kernel
$\widehat{\chi}_+$, whose $(m,n)$-element is given
by the integral kernel (\ref{eqn:hatchi1}).
The $(m,n)$-element of its $q$-th power
$(\widehat{\chi}_+)^{q}, q=2,3, \cdots,$
will be the following kernel given by
$(q-1)$-multiple integral,
$$
[(\widehat{\chi}_{+})^{q}]^{m, n}(x, y) 
= \sum_{m_1=0}^{M} \int_{\R} d x^{(m_1)}
\cdots \sum_{m_{q-1}=0}^{M}
\int_{\R} d x^{(m_{q-1})}
\widehat{\chi}_{+}^{m, m_1}(x, x^{(m_1)}) \times \cdots 
\times \widehat{\chi}_{+}^{m_{q-1}, n}
(x^{(m_{q-1})}, y). 
$$
Then we have
\begin{eqnarray}
F_{jk}[\chi] &=& \langle j-1| k-1 \rangle
\nonumber\\
&+& \sum_{m=0}^{M} \int_{\R} dx
\sum_{n=0}^{M} \int_{\R} dy
\langle j-1 | t_{m}, x \rangle \chi_{m}(x)
\Big\{ \delta_{m, n} \delta(x-y)
+ \sum_{q=1}^{\infty} 
[(\widehat{\chi}_{+})^{q}]^{m,n}(x, y) \Big\}
\langle t_n, y|k-1 \rangle.
\nonumber\\
\label{eqn:F5}
\end{eqnarray}
Let
$\hat{1}^{m,n}(x,y) \equiv \delta_{m,n} \delta(x-y)$.
Since it is easy to see that
\begin{equation}
\sum_{\ell=0}^{M} \int_{\R} d y
[\hat{1}-\widehat{\chi}_+]^{m, \ell}(x, y)
\Big\{ \delta_{\ell, n} \delta(y-z)
+ \sum_{q=1}^{\infty} 
[(\widehat{\chi}_{+})^{q}]^{\ell, n}(y, z) \Big\}
= \hat{1}^{m, n}(x, z),
\label{eqn:Inverse1}
\end{equation}
we can write
$$
\delta_{m, n} \delta(x-y)
+ \sum_{q=1}^{\infty} 
[(\widehat{\chi}_{+})^{q}]^{m,n}(x, y) 
=\left[ \frac{\hat{1}}
{\hat{1}-\widehat{\chi}_{+}} \right]^{m,n}(x, y).
$$
Then (\ref{eqn:F5}) becomes
\begin{equation}
F_{jk}[\chi] = 
\delta_{jk}+\sum_{m=0}^{M} \int_{\R} dx \,
B_{j}^{(m)}(x) C_{k}^{(m)}(x),
\label{eqn:F6}
\end{equation}
where
\begin{eqnarray}
B_{j}^{(m)}(x) &=&
\sum_{n=0}^{M} \int_{\R} dy \,
\langle j-1 | t_{n}, y \rangle \chi_{n}(y)
\left[ \frac{\hat{1}}{\hat{1}-\widehat{\chi}_{+}} 
\right]^{n,m}(y,x),
\nonumber\\
C_{k}^{(m)}(x) &=&
\langle t_m, x|k-1 \rangle.
\label{eqn:BC1}
\end{eqnarray}

\subsection{Fredholm Determinant and 
Determinantal Process}

Let $\widetilde{B}_{j}^{(m)}(x)$ and 
$\widetilde{C}_{j}^{(m)}(x)$,
$0 \leq m \leq M, 1 \leq j \leq N$, be
square integrable continuous functions.
Then the following formulae can be proved;
\begin{eqnarray}
&& \det_{1 \leq j, k \leq N}
\Big[ \delta_{j, k} 
+ \sum_{m=0}^{M} \int_{\R} dx \,
\widetilde{B}_j^{(m)}(x) \widetilde{C}_{k}^{(m)}(x) \Big]
\nonumber\\
&=& \det_{1 \leq j^{(m)}, j^{(n)} \leq N,
0 \leq m, n \leq M}
\left[ \delta_{m,n} \delta_{j^{(m)}, j^{(n)}}
+ \int_{\R} dx \, \widetilde{B}^{(m)}_{j^{(m)}}(x)
\widetilde{C}^{(m)}_{j^{(n)}}(x) \right]
\nonumber\\
&=& \sum_{N_0=0}^{N} \sum_{N_1=0}^{N} \cdots
\sum_{N_M=0}^{N} 
\prod_{m=0}^{M} \frac{1}{N_{m}!}
\int_{\R^{N_0}} \prod_{j=1}^{N_0} dx_{j}^{(0)}
\int_{\R^{N_1}} \prod_{j=1}^{N_1} dx_{j}^{(1)} 
\cdots
\int_{\R^{N_M}} \prod_{j=1}^{N_M} dx_{j}^{(M)}
\nonumber\\
&& \qquad \qquad 
\times \det_{1 \leq j \leq N_m,
1 \leq k \leq N_n, 0 \leq m, n \leq M}
\Big[ \sum_{p=1}^{N} \widetilde{C}_{p}^{(m)}(x_j^{(m)})
\widetilde{B}_{p}^{(n)}(x_{k}^{(n)}) \Big].
\label{eqn:Fred1}
\end{eqnarray}
The expansion formula of the last expression
defines the Fredholm determinant, which is abbreviated as
$$
\Det \left[
\hat{1}^{m,n}(x, y)+
\sum_{p=1}^{N} \widetilde{C}_{p}^{(m)}(x) 
\widetilde{B}_{p}^{(n)}(y) \right].
$$
The generating function (\ref{eqn:F1}) with 
(\ref{eqn:F6}) and (\ref{eqn:BC1}) is then expressed as
the Fredholm determinant, 
\begin{eqnarray}
&& \Psi_{N}[\chi] =
\Det
\left[ \hat{1}^{m, n}(x, z)
+\sum_{\ell=0}^{M} \int_{\R} dy
 \langle t_m, x| {\cal P}^N
|t_{\ell}, y \rangle \chi_{\ell}(y)
\left[ \frac{\hat{1}}{\hat{1}-\widehat{\chi}_{+}}
\right]^{\ell,n}(y, z) \right],
\nonumber\\
\label{eqn:Det1}
\end{eqnarray}
where
\begin{equation}
{\cal P}^N= \sum_{p=1}^{N} |p-1 \rangle
\langle p-1 |.
\label{eqn:hatS}
\end{equation}
Note that (\ref{eqn:Inverse1}) implies that
$$
\hat{1}^{m,n}(x, z)
= \sum_{\ell=0}^{M} \int_{\R} dy
\Big\{ \hat{1}^{m,\ell}(x, y)
-\langle t_m ,x |\hat{1}_+|t_{\ell}, y \rangle
\chi_{\ell}(y) \Big\}
\left[ \frac{\hat{1}}{\hat{1}-\widehat{\chi}_{+}}
\right]^{\ell, n}(y, z).
$$
Plugging this into (\ref{eqn:Det1}),
we have
\begin{eqnarray}
\Psi_{N}[\chi] 
&=& \Det \Bigg[
\sum_{\ell=0}^{M} \int_{\R} d y
\Big[ \hat{1}^{m, \ell}(x, y)
+ \langle t_m, x| ({\cal P}^N-\hat{1}_{+})
|t_{\ell}, y \rangle \chi_{\ell}(y) 
\Big] 
\left[ \frac{\hat{1}}{\hat{1}-\widehat{\chi}_{+}}
\right]^{\ell, n} (y, z) \Bigg]
\nonumber\\
&=& \Det \Big[ \hat{1}^{m,n}(x, y)
+\langle t_{m}, x|({\cal P}^N-\hat{1}_{+})
|t_{n}, y \rangle \chi_{n}(y) \Big]
\nonumber\\
&& \hskip 4cm \times
\Det \left[ \left[
\frac{\hat{1}}{\hat{1}-\widehat{\chi}_{+}}
\right]^{m,n}(x, y) \right]
\nonumber\\
&=& 
\Det \Big[ \hat{1}^{m,n}(x, y)
+\langle t_{m}, x|({\cal P}^N-\hat{1}_{+})
|t_{n}, y \rangle \chi_{n}(y) \Big].
\label{eqn:Det2}
\end{eqnarray}
Here
$$
\Det \left[ \left[
\frac{\hat{1}}{\hat{1}-\widehat{\chi}_{+}}
\right]^{m,n}(x, y) \right]
=1 \Big/
\Det\Big[ [\hat{1}-\widehat{\chi}_{+}]^{m,n}
(x, y) \Big],
$$
and we have used 
$\Det\Big[ [\hat{1}-\widehat{\chi}_{+}]^{m,n}
(x, y) \Big]=1$,
which is concluded from the fact that
$[\hat{1}-\widehat{\chi}_{+}]^{m,n}(x, y)=0$
for $m < n$ by the definition of $\widehat{\chi}_{+}$,
(\ref{eqn:hatchi1}) with (\ref{eqn:1+}), 
and it is $\delta(x-y)$ for $m=n$.

Let
\begin{eqnarray}
S^{m,n}_{N} (x,y) &=& \langle t_m, x| {\cal P}^{N} 
|t_n, y \rangle, \nonumber\\
\tS^{m,n}_{N}(x,y)
&=&
- \langle t_m, x| (\hat{1}_{+}-{\cal P}^N) | t_n, y \rangle.
\nonumber
\end{eqnarray}
Following the formulae given in Section 2,
$S_{N}^{m,n}$ is determined as 
\begin{eqnarray}
S^{m,n}_{N} (x,y) &=& \sum_{p=1}^{N}
\langle t_m, x|p-1 \rangle \langle p-1|t_n, y \rangle
\nonumber\\
&=& \frac{1}{\sqrt{2t_m}}
e^{-x^2/4t_m+y^2/4t_n}
\sum_{p=1}^{N}
\left(\frac{t_n}{t_m}\right)^{(p-1)/2}
\left\langle \frac{x}{\sqrt{2t_m}} \Big|
p-1 \right\rangle
\left\langle p-1 \Big|
\frac{y}{\sqrt{2t_n}} \right\rangle
\nonumber\\
&=& \frac{1}{\sqrt{2t_m}}
e^{-x^2/4t_m+y^2/4t_n} \sum_{k=0}^{N-1} 
\left(\frac{t_n}{t_m}\right)^{k/2}
\varphi_{k}\left(\frac{x}{\sqrt{2t_m}}\right)
\varphi_{k}\left(\frac{y}{\sqrt{2t_n}}\right).
\nonumber
\end{eqnarray}
Combination of this with (\ref{eqn:Mehler}) gives
$\widetilde{S}_N^{m,n}$ as
\begin{eqnarray}
&& \tS_N^{m,n}(x,y) = S^{m,n}_{N}(x,y)
-{\bf 1}_{\{m > n\}} p(t_{m}-t_{n}, x|y) 
\nonumber\\
&& \quad = \left\{
   \begin{array}{ll}
\displaystyle{\frac{1}{\sqrt{2t_m}}
e^{-x^2/4t_m+y^2/4t_n} \sum_{k=0}^{N-1} 
\left(\frac{t_n}{t_m}\right)^{k/2}
\varphi_{k}\left(\frac{x}{\sqrt{2t_m}}\right)
\varphi_{k}\left(\frac{y}{\sqrt{2t_n}}\right)}
& \quad \mbox{if} \quad m \leq n \\
\displaystyle{-\frac{1}{\sqrt{2t_m}}
e^{-x^2/4t_m+y^2/4t_n} \sum_{k=N}^{\infty} 
\left(\frac{t_n}{t_m}\right)^{k/2}
\varphi_{k}\left(\frac{x}{\sqrt{2t_m}}\right)
\varphi_{k}\left(\frac{y}{\sqrt{2t_n}}\right)}
& \quad \mbox{if} \quad m > n,
   \end{array} \right. 
\label{eqn:tSNmn1}
\end{eqnarray}
where ${\bf 1}_{\{\omega\}}$ is the
indicator of a condition $\omega$;
${\bf 1}_{\{\omega \}}=1$ if
$\omega$ is satisfied and 
${\bf 1}_{\{\omega \}}=0$ otherwise.

As shown above, not $\{\tS_N^{m,n}(x,y)\}$ themselves,
but determinants of matrices made of them are
observables.
By definition of determinant, 
factors $\{e^{-x^2/4t_m+y^{2}/4t_n}\}$ 
of $\{\tS_{N}^{m,n}(x,y)\}$ in 
(\ref{eqn:tSNmn1}) are completely cancelled out,
when we calculate determinants.
So here we define the following matrix-kernel 
by omitting these factors in $\{\tS_N^{m,n}(x,y)\}$,
\begin{equation}
\mbK_N(t_m,x;t_n,y) = \left\{
   \begin{array}{ll}
\displaystyle{
\frac{1}{\sqrt{2t_m}}
\sum_{k=0}^{N-1}
\left(\frac{t_n}{t_m}\right)^{k/2}
\varphi_{k}\left(\frac{x}{\sqrt{2t_m}}\right)
\varphi_{k}\left(\frac{y}{\sqrt{2t_n}}\right)}
& \quad \mbox{if} \quad m \leq n \\
\displaystyle{-\frac{1}{\sqrt{2t_m}}
\sum_{k=N}^{\infty} 
\left(\frac{t_n}{t_m}\right)^{k/2}
\varphi_{k}\left(\frac{x}{\sqrt{2t_m}}\right)
\varphi_{k}\left(\frac{y}{\sqrt{2t_n}}\right)}
& \quad \mbox{if} \quad m > n,
   \end{array} \right. 
\label{eqn:KN1}
\end{equation}
and rewrite (\ref{eqn:Det2}) as
\begin{equation}
\Psi_{N}[\chi]=
\Det \Big[ \delta_{m,n} \delta (x-y)
+\mbK_N(t_m, x; t_n, y)
\chi_n (y) \Big].
\label{eqn:Det3}
\end{equation}
This Fredholm determinant is 
by definition expanded as
\begin{eqnarray}
&& \Det \Big[
\delta_{m,n} \delta(x-y)
+\mbK_{N}(t_m, x; t_n, y) \chi_{n}(y) \Big]
\nonumber\\
&=& \sum_{N_{0}=0}^{N} \sum_{N_{1}=0}^{N} \cdots
\sum_{N_{M}=0}^{N}
\prod_{m=0}^{M}\frac{1}{N_m !}
\int_{\R^{N_{0}}} \prod_{j=1}^{N_0} d x_{j}^{(0)}
\int_{\R^{N_{1}}} \prod_{j=1}^{N_1} d x_{j}^{(1)}
 \cdots
\int_{\R^{N_{M}}} 
\prod_{j=1}^{N_{M}} d x_{j}^{(M)} \nonumber\\
&& \quad \quad
\times \prod_{m=0}^{M} \prod_{j=1}^{N_{m}} 
\chi_{m} \Big(x_{j}^{(m)} \Big) 
\det_{1 \leq j \leq N_{m}, 1 \leq k \leq N_{n},
0 \leq m, n \leq M} \Bigg[
\mbK_{N}(t_m, x_{j}^{(m)}; t_n, x_{k}^{(n)} )
\Bigg]. 
\label{eqn:FredExp1}
\end{eqnarray}
Comparison of (\ref{eqn:FredExp1}) and 
(\ref{eqn:kPhiA1}) determines all of the
multitime correlation functions.
Now we summarize the above results as a theorem.

\begin{thm}
\label{thm:finite}
The temporally homogeneous noncolliding BM, 
$\Xi^{\bf X}_{N}(t)=\{\X(t) \}$, 
starting from the GUE-eigenvalue distribution (\ref{eqn:nu01})
at time $t_0 > 0$,
is a finite determinantal process in the following sense.
\begin{description}
\item{(i)} \quad
The multitime generating function is given by
the Fredholm determinant (\ref{eqn:Det3}),
where the matrix-kernel 
$\mbK_{N}$ is given by (\ref{eqn:KN1}).
\item{(ii)} \quad
Any multitime correlation function is given by a determinant;
for any $M \geq 0$, any sequence $\{N_m \}_{m=0}^{M}$ of 
positive integers less than or equal to $N$,
any time sequence $t_0 < t_1 < \cdots < t_M < \infty$,
the $(N_0, \dots, N_{M})$-multitime 
correlation function is given by
\begin{equation}
\rho_{N} 
\left(t_0, \{\x^{(0)}_{N_1}\}; t_1, \{\x^{(1)}_{N_1}\}; 
\dots; t_{M}, \{\x^{(M)}_{N_{M}}\} \right) 
=\det_{1 \leq j \leq N_{m}, 1 \leq k \leq N_{n},
0 \leq m, n \leq M} \Bigg[
\mbK_{N}(t_m, x_{j}^{(m)}; t_n, x_{k}^{(n)} )
\Bigg]. 
\label{eqn:THfinite}
\end{equation}
\end{description}
\end{thm}
\vskip 0.5cm

The following relations hold for the Hermite
polynomials,
\begin{eqnarray}
\label{eqn:Hermite2}
&& H_{k+1}(x)=2x H_{k}(x)-2kH_{k-1}(x), \\
\label{eqn:Hermite3}
&& \frac{d}{dx} H_k(x) = 2k H_{k-1}(x),
\qquad k=1,2,3, \cdots.
\end{eqnarray}
From (\ref{eqn:Hermite2}), the Christoffel-Darboux
formula is derived for the Hermite
orthonormal functions $\{\varphi_k(x)\}_{k \in \N_0}$,
\begin{equation}
\sum_{k=0}^{N-1} \varphi_{k}(x) \varphi_{k}(y)
=\sqrt{\frac{N}{2}}
\frac{\varphi_N(x) \varphi_{N-1}(y)
-\varphi_{N-1}(x) \varphi_{N}(y)}{x-y}
\label{eqn:CD1}
\end{equation}
for $x \not= y$.
Eq. (\ref{eqn:Hermite3}) can be used to evaluate
the limit $y \to x$ in (\ref{eqn:CD1})
and we find
$$
\sum_{k=0}^{N-1} \Big\{ \varphi_k(x) \Big\}^2
=N \Big\{ \varphi_{N}(x) \Big\}^2
-\sqrt{N(N+1)} \varphi_{N+1}(x) \varphi_{N-1}(x).
$$
Then the matrix-kernel have the following
simpler expressions if $m=n$,
\begin{equation}
\mbK_{N}(t_m, x; t_m, y) = \left\{
   \begin{array}{l}
\displaystyle{
\sqrt{\frac{N}{2}}
\frac{\varphi_{N}(x/\sqrt{2t_m}) \varphi_{N-1}(y/\sqrt{2t_m})
-\varphi_{N-1}(x/\sqrt{2t_m}) \varphi_N(y/\sqrt{2t_m})}
{x-y}
}  \\
 \hskip 6cm \mbox{if} \quad x \not= y  \\
\displaystyle{
\frac{1}{\sqrt{2 t_{m}}}
\left[ N \left\{
\varphi_N \left(\frac{x}{\sqrt{2t_m}}\right) \right\}^2
-\sqrt{N(N+1)}
\varphi_{N-1}\left(\frac{x}{\sqrt{2t_m}}\right)
\varphi_{N+1}\left(\frac{x}{\sqrt{2t_m}}\right)
\right]
} \\
\hskip 6cm \mbox{if} \quad x=y. 
   \end{array} \right. 
\label{eqn:KN2}
\end{equation}

\SSC{INFINITE PARTICLE SYSTEMS}
\subsection{Wigner's Semicircle Law}

The density of $\Xi^{\bf X}_{N}(t)$ is given by
\begin{eqnarray}
\rho_{N}(t,x) &=& \frac{1}{\sqrt{2t}}
\sum_{k=0}^{N-1} \left\{ \varphi_{k}
\left(\frac{x}{\sqrt{2t}}\right) \right\}^2
\nonumber\\
&=& 
\frac{1}{\sqrt{2 t}}
\left[ N \left\{
\varphi_N \left(\frac{x}{\sqrt{2t}}\right) \right\}^2
-\sqrt{N(N+1)}
\varphi_{N-1}\left(\frac{x}{\sqrt{2t}}\right)
\varphi_{N+1}\left(\frac{x}{\sqrt{2t}}\right)
\right],
\nonumber
\end{eqnarray}
as a special case 
($M=0, N_{0}=1$ with setting $t_0=t$)
of Theorem \ref{thm:finite} with (\ref{eqn:KN2}).
It is easy to confirm that
$\int_{-\infty}^{\infty} \rho_{N}(t,x) dx=N$ by the
orthonormality of $\varphi_{k}(x)$.
The following estimations for asymptote in $N \to \infty$
are established \cite{Bat53,Sze75}.
Let $\varepsilon$ and $\omega$ be the fixed
positive numbers. We have
\begin{eqnarray}
&(i)& \quad \varphi_{N}\Big( \sqrt{2N+1} \cos \phi \Big)
= \frac{1}{\sqrt{\pi \sin \phi}}
\left( \frac{2}{N} \right)^{1/4} \nonumber\\
&& \qquad \qquad \times
\left\{ \sin \left[
\left( \frac{N}{2}+\frac{1}{4} \right)
(\sin 2 \phi - 2 \phi)+ \frac{3}{4} \pi \right]
+{\cal O}\left( \frac{1}{N} \right) \right\},
\quad
\varepsilon \leq \phi \leq \pi -\varepsilon 
\nonumber\\
&(ii)& \quad 
\varphi_{N}\Big( \sqrt{2N+1} \cosh \phi \Big)
= \frac{1}{\sqrt{2 \pi \sinh \phi}}
\left( \frac{1}{2N} \right)^{1/4} \nonumber\\
&& \qquad \qquad \times
\exp \left[
\left( \frac{N}{2}+\frac{1}{4} \right)
(2 \phi - \sinh 2 \phi )+ \frac{3}{4} \pi \right]
\left\{ 1+ {\cal O} \left( \frac{1}{N} \right) \right\},
\quad \varepsilon \leq \phi \leq \omega.
\nonumber
\end{eqnarray}
Using them, we will have the asymptote
of the density profile in $N \to \infty$,
\begin{equation}
\rho_N(t,x) \simeq \left\{
   \begin{array}{ll}
\displaystyle{
\frac{1}{\pi\sqrt{2t}}
\sqrt{2N- \frac{x^2}{2t}}} \qquad &
 \mbox{if} \quad -2 \sqrt{Nt} \leq x \leq 2 \sqrt{Nt}  \\
 & \\
0 &
 \mbox{otherwise.} 
   \end{array} \right. 
\label{eqn:dens2}
\end{equation}
The distribution of $N$ particles has a finite
support, whose interval $\propto \sqrt{N}$,
and thus $\rho \sim \sqrt{N} \to \infty$
as $N \to \infty$ for fixed $ 0 < t < \infty$.
If we set $x= 2 \sqrt{Nt} \xi$, we see
\begin{equation}
\lim_{N \to \infty}
\frac{1}{N} \rho_{N}(t, x) dx =
\left\{
   \begin{array}{ll}
\displaystyle{
\frac{2}{\pi} \sqrt{1-\xi^2}} \, d \xi  &
 \mbox{if} \quad -1 \leq \xi \leq 1  \\
 & \\
0 &
 \mbox{otherwise},
   \end{array} \right. 
\label{eqn:Wigner1}
\end{equation}
which is known as Wigner's semicircle law \cite{Meh04}.
See also \cite{RS93}.
In the following we consider scaling limit,
in which long-term limit $t \to \infty$ is taken 
at the same time with $N \to \infty$.

\subsection{Bulk Scaling Limit and Homogeneous Infinite System}

First we consider the central region $x \simeq 0$
in the semicircle-shaped profile of particle density
in the scaling limit
\begin{equation}
t \simeq N \to \infty.
\label{eqn:BulkLimit}
\end{equation}
In this limit the system becomes homogeneous also
in space with a constant density
$\rho=1/\pi$.
We call this the bulk scaling limit.

\begin{thm}
\label{thm:infinite_sin}
For any $M \geq 0$, any sequence $\{N_m \}_{m=0}^{M}$
of positive integers, and any strictly increasing
sequence $\{ s_m \}_{m=1}^{M}$ of
positive numbers
\begin{eqnarray}
&& \lim_{N \to \infty}
\rho_{N}(N, \{\x^{(0)}_{N_0}\} ;
N+2s_1, \{\x^{(1)}_{N_1}\}; \cdots ;
N+2s_M, \{\x^{(M)}_{N_M}\})
\nonumber\\
&& \qquad \qquad =
\det_{1 \leq j \leq N_{m}, 1 \leq k \leq N_{n},
0 \leq m, n \leq M}
\Big[ 
\bK(s_m, x_{j}^{(m)}; s_n, x_{k}^{(n)} )
 \Big]
\nonumber\\
&& \qquad \qquad \equiv \rho_{\rm sin}
\Big(0, \xi^{(0)}_{N_0}; s_1, \xi^{(0)}_{N_1}; 
\cdots ; s_M, \xi^{(M)}_{N_M} \Big),
\label{eqn:bulk1}
\end{eqnarray}
where 
\begin{equation}
\bK(t,y;s,x) = \left\{
   \begin{array}{ll}
\displaystyle{
\frac{1}{\pi} \int_{0}^{1} d u \,
e^{(s-t)u^2} \cos (u(y-x))
} 
& \mbox{if} \quad t < s  \\
& \\
\displaystyle{
\frac{\sin(y-x)}{\pi(y-x)}
}
& \mbox{if} \quad t=s \\
& \\
\displaystyle{
-\frac{1}{\pi} \int_{1}^{\infty} d u \,
e^{-(t-s)u^2} \cos (u(y-x))
}
& \mbox{if} \quad t > s.
   \end{array} \right. 
\label{eqn:Kbulk1}
\end{equation}
\end{thm}
\vskip 0.5cm
\noindent{\it Proof.} \quad
For any $u \in \R$, the formula
\begin{eqnarray}
&& \lim_{\ell \to \infty} (-1)^{\ell} \ell^{1/4}
\varphi_{2 \ell} 
\left(\frac{u}{2 \sqrt{\ell}} \right)
= \frac{1}{\sqrt{\pi}} \cos u, 
\nonumber\\
&& \lim_{\ell \to \infty} (-1)^{\ell} \ell^{1/4}
\varphi_{2 \ell+1} 
\left(\frac{u}{2 \sqrt{\ell}} \right)
= \frac{1}{\sqrt{\pi}} \sin u
\label{eqn:formula1}
\end{eqnarray}
are known \cite{Bat53,Sze75}.
We note that
\begin{eqnarray}
\left(\frac{t_n}{t_m}\right)^{\alpha}
&=& \left(\frac{N+2s_n}{N+2s_m}\right)^{\alpha}
\nonumber\\
&=& \left\{ \left(1+\frac{2s_n}{N}\right)^{N}
\left(1+\frac{2s_{m}}{N}\right)^{-N} 
\right\}^{\alpha/N}
\nonumber\\
&\simeq& e^{2(s_n-s_m)\alpha/N}
\nonumber
\end{eqnarray}
for $N \gg 1$ with fixed number $\alpha$.
Then (\ref{eqn:KN1}) with $m \leq n$
is evaluated in $N \to \infty$ as
\begin{eqnarray}
&& \mbK_N(t_m,y; t_n, x) \nonumber\\
&& \qquad \simeq \frac{1}{\pi N}
\sum_{\ell=0}^{N/2-1}
e^{2(s_n-s_m) \ell/N}
\sqrt{\frac{N}{2 \ell}}
\left\{ \cos \left(\sqrt{\frac{2\ell}{N}} \, y \right)
\cos \left(\sqrt{\frac{2\ell}{N}} \, x \right) 
+\sin \left(\sqrt{\frac{2\ell}{N}} \, y \right)
\sin \left(\sqrt{\frac{2\ell}{N}} \, x \right) \right\}
\nonumber\\
&& \qquad \simeq \frac{1}{2\pi}
\int_{0}^{1} d \lambda \,
e^{(s_n-s_m) \lambda} \frac{1}{\sqrt{\lambda}}
\Big\{ \cos(y \sqrt{\lambda}) \cos(x \sqrt{\lambda})
+ \sin(y \sqrt{\lambda}) \sin(x \sqrt{\lambda}) \Big\}
\nonumber\\
&& \qquad = \frac{1}{\pi} \int_{0}^{1}
d u \,
e^{(s_n-s_m) u^2}
\cos (u(y-x)).
\nonumber
\end{eqnarray}
In particular, when $m=n$, {\it i.e.,}
$s_n-s_m=0$,
the integration is readily performed to have
$\int_{0}^{1} d u \, \cos(u(y-x))
=\sin(y-x)/(y-x)$.
Similar evaluation in $N \to \infty$ can be done
also for (\ref{eqn:KN1}) with $m > n$.
\qed
\vskip 0.5cm

In the bulk scaling limit (\ref{eqn:BulkLimit}),
the temporally and spatially homogeneous infinite
particle system is obtained, whose multitime correlation functions
are given by (\ref{eqn:bulk1}).
The matrix-kernel (\ref{eqn:Kbulk1}) is called
the extended sine kernel in \cite{TW04}.
In the present paper, we will call it simply
``sine kernel".
The system with the sine kernel was
studied by Spohn \cite{Spo87}, Osada \cite{Osa96,Osa04},
and Nagao and Forrester \cite{NF98b},
as an infinite particle limit of 
Dyson's BM model
with $\beta=2$ \cite{Dys62}.
See also \cite{TW04,AM05}.

\subsection{Soft-edge Scaling Limit and 
Spatially Inhomogeneous Infinite System}

Next we consider the scaling limit
\begin{equation}
t \simeq N^{1/3}
\quad \mbox{and} \quad
x \simeq 2 N^{2/3}.
\label{eqn:soft1}
\end{equation}
Since (\ref{eqn:soft1}) gives
$x^2/2t \simeq 2N$, 
the vicinity of the right edge 
of semicircle-shaped profile (\ref{eqn:dens2})
will be closed up, and we will obtain
a spatially inhomogeneous infinite particle system
in this scaling limit.
Following the random matrix theory \cite{Meh04},
we call (\ref{eqn:soft1}) the soft-edge scaling limit.

In order to describe the limit, we introduce 
the Airy function
\begin{equation}
{\rm Ai}(x)=\frac{1}{2 \pi} \int_{-\infty}^{\infty} dk \
{\rm e}^{i (x k + k^3/3)}.
\label{eqn:Airy}
\end{equation}
It is a solution of equation
\begin{equation}
\frac{d^2}{dx^2} {\rm Ai}(x)=x {\rm Ai}(x),
\label{eqn:Aeq1}
\end{equation}
which behaves as
\begin{eqnarray}
{\rm Ai}(x) &\simeq& \frac{1}{2 \sqrt{\pi} x^{1/4}}
\exp\left(-\frac{2}{3} x^{3/2} \right), \nonumber\\
{\rm Ai}(-x) &\simeq& \frac{1}{\sqrt{\pi} x^{1/4}}
\cos \left(\frac{2}{3} x^{3/2}-\frac{\pi}{4} \right)
\quad \mbox{in} \quad x \to \infty. \nonumber
\end{eqnarray}
In the proof of the following theorem,
we will use the formula
\begin{equation}
\lim_{\ell \to \infty}
2^{-1/4} \ell^{1/12}
\varphi_{\ell}
\left( \sqrt{2\ell}+\frac{u}{\sqrt{2}}\ell^{-1/6} \right)
=\Ai(u)
\quad \mbox{for} \quad u \in \R.
\label{eqn:formula2}
\end{equation}
Let
\begin{equation}
a_N(s)=2N^{2/3}+2N^{1/3} s -s^2,
\label{eqn:aNs}
\end{equation}
and
$
\x_{N'}(s)=
(a_N(s)+x_1, a_N(s)+x_2, \cdots, a_N(s)+x_{N'}).
$

\begin{thm}
\label{thm:infinite_Ai}
For any $M \geq 0$, any sequence $\{N_m \}_{m=0}^{M}$
of positive integers, and any strictly increasing
sequence $\{ s_m \}_{m=1}^{M}$ of
positive numbers
\begin{eqnarray}
&& \lim_{N \to \infty}
\rho_{N} (N^{1/3}, \{\x^{(0)}_{N_0}(0)\} ;
N^{1/3}+2s_1, \{\x^{(1)}_{N_1}(s_1)\}; \cdots ;
N^{1/3}+2s_M, \{\x^{(M)}_{N_M}(s_M)\})
\nonumber\\
&& \qquad \qquad = 
\det_{1 \leq j \leq N_{m}, 1 \leq k \leq N_{n},
0 \leq m, n \leq M}
\Big[ 
\sK(s_m, x_{j}^{(m)}; s_n, x_{k}^{(n)} )
\Big]
\nonumber\\
&& \qquad \qquad \equiv \rho_{\rm Ai}
\Big(0, \xi^{(0)}_{N_0}; s_1, \xi^{(0)}_{N_1}; 
\cdots ; s_M, \xi^{(M)}_{N_M} \Big),
\label{eqn:edge1}
\end{eqnarray}
where 
\begin{equation}
\sK(t,y; s,x) = \left\{
   \begin{array}{ll}
\displaystyle{
\int_{-\infty}^{0} d \lambda \,
e^{(s-t)\lambda} \Ai(y-\lambda) \Ai(x-\lambda)
} 
& \mbox{if} \quad t \leq s  \\
& \\
\displaystyle{
-\int_{0}^{\infty} d \lambda \,
e^{-(t-s)\lambda} \Ai(y-\lambda) \Ai(x-\lambda)
}
& \mbox{if} \quad t > s.
   \end{array} \right. 
\label{eqn:KAiry1}
\end{equation}
\end{thm}
\vskip 0.5cm
\noindent{\it Proof.} \quad
Putting $k=N-p-1$ in the summation of
(\ref{eqn:KN1}) for $m \leq n$, we have
$$
\mbK_{N}(t_m,y; t_n,x)=
\left( \frac{t_n}{t_m} \right)^{(N-1)/2}
\frac{1}{\sqrt{2t_m}}
\sum_{p=0}^{N-1}
\left(\frac{t_n}{t_m}\right)^{-p/2}
\varphi_{N-p-1} \left(\frac{y}{\sqrt{2t_m}}\right)
\varphi_{N-p-1} \left(\frac{x}{\sqrt{2t_n}}\right).
$$
With the same reason mentioned above eq. (\ref{eqn:KN1}),
we can omit the factor $(t_n/t_m)^{(N-1)/2}$.
Since, when we set $t_m=N^{1/3}+2s_m$,
$$
\frac{a_N(s_m)+x}{\sqrt{2t_m}}
=\sqrt{2N}+\frac{1}{\sqrt{2}}N^{-1/6} x
+{\cal O}(N^{-1/2}),
$$
we can use the formula (\ref{eqn:formula2});
\begin{eqnarray}
\varphi_{N-p-1} \left( \frac{x}{\sqrt{2t_m}}\right)
&\simeq& \varphi_{N-p-1}
\left( \sqrt{2N}+\frac{1}{\sqrt{2}} N^{-1/6}x \right)
\nonumber\\
&\simeq& \varphi_{N-p-1}
\left( \sqrt{2(N-p-1)} + \frac{1}{\sqrt{2}}
(N-p-1)^{-1/6}
\left\{ x+ \frac{p}{N^{1/3}} \right\} \right).
\nonumber\\
&\simeq& 2^{1/4} N^{-1/12} \Ai \left(x+\frac{p}{N^{1/3}}\right).
\nonumber
\end{eqnarray}
For 
$$
\left( \frac{t_n}{t_m}\right)^{-p/2}
= \left[
\left( \frac{1+2s_n/N^{1/3}}{1+2s_m/N^{1/3}}
\right)^{N^{1/3}/2} \right]^{-p/N^{1/3}}
\simeq e^{-(s_n-s_m)p/N^{1/3}}
\quad \mbox{in} \quad N \to \infty,
$$
we have
\begin{eqnarray}
&&\mbK_{N}(N^{1/3}+2s_m, a_N(s_m)+y;
N^{1/3}+2s_n, a_N(s_n)+x)
\nonumber\\
&& \qquad \sim
\frac{1}{N^{1/3}} \sum_{p=0}^{N-1}
e^{-(s_n-s_m)p/N^{1/3}}
\Ai \left(y+\frac{p}{N^{1/3}}\right)
\Ai \left(x+\frac{p}{N^{1/3}}\right)
\nonumber\\
&& \qquad \simeq \int_{0}^{\infty} d u \,
e^{-(s_n-s_m)u}
\Ai(y+u) \Ai(x+u)
\quad \mbox{in} \quad N \to \infty.
\nonumber
\end{eqnarray}
Note that the factor $(t_n/t_m)^{(N-1)/2}$,
which is irrelevant in calculating determinants,
was omitted in the second line in the 
above equations.
Put $u=-\lambda$ to obtain the
expression (\ref{eqn:KAiry1}).
Similar evaluation in $N \to \infty$
of (\ref{eqn:KN1}) can be done
also for $m>n$.
\qed

The infinite system obtained by the soft-edge scaling
limit (\ref{eqn:soft1}) is temporally homogeneous, but
spatially inhomogeneous as shown by 
the ``Airy kernel" (\ref{eqn:KAiry1}).
Pr\"ahofer and Spohn \cite{PS02} and Johansson \cite{Joh03}
studied the right-most path in the present system
and called it the Airy process $A(t)$.
For a given $t >0$, $A(t)$ has 
distribution of the celebrated Tracy-Widom
distribution, which is governed by the 
Painlev\'e II equation \cite{TW94}.
Recently Tracy and Widom derived a system of
partial differential equations (PDE), 
which govern the Airy kernel (\ref{eqn:KAiry1})
\cite{TW03}.
They also discussed other determinantal processes
by PDE \cite{TW04}.
See also \cite{AM99,AM05}.

\SSC{DETERMINANTAL PROCESSES ASSOCIATED \\
WITH SPECTRAL PROJECTIONS}

\subsection{Spectral Projections}

First we note that, following the notations in Section 2,
if we set
$\tau=\log t, \tau'=\log t',
\zeta=y/\sqrt{2t}, \zeta'=x/\sqrt{2t'}$, and
$\mbK_N(t,y;s,x)dy=\widetilde{\mbK}_N
(\tau, \zeta; \tau', \zeta') d \zeta$
with $d \zeta=dy/\sqrt{2t}$, the matrix-kernel
of the determinantal process of noncolliding
BM with finite number of
particles $N < \infty$ given by (\ref{eqn:KN1}) 
is rewritten as
\begin{equation}
\widetilde{\mbK}_N(\tau, \zeta; \tau', \zeta') = \left\{
   \begin{array}{ll}
\displaystyle{
\langle \zeta | e^{(\tau'-\tau) 
\widehat{\cal H}_{\varphi}}
{\cal P}_{\varphi} | \zeta' \rangle
}
& \quad \mbox{if} \quad \tau \leq \tau' \\
\displaystyle{
-\langle \zeta | e^{-(\tau-\tau') 
\widehat{\cal H}_{\varphi}}
(1-{\cal P}_{\varphi}) | \zeta' \rangle
}
& \quad \mbox{if} \quad \tau > \tau',
   \end{array} \right. 
\label{eqn:tKN1}
\end{equation}
where $\widehat{\cal H}_{\varphi}$ is given by
(\ref{eqn:Hphi1}),
and ${\cal P}_{\varphi}$
is a projection operator defined by
\begin{equation}
{\cal P}_{\varphi}
= \sum_{0 \leq k \leq N-1}
|k \rangle \langle k |.
\label{eqn:projection1}
\end{equation}
It is called the extended Hermite kernel in \cite{TW04}.

The trigonometric functions of the form
\begin{eqnarray}
S_{-}(\sqrt{\lambda}x) &=&
\frac{1}{\sqrt{2\pi} \lambda^{1/4}} 
\sin(\sqrt{\lambda}x),
\nonumber\\
S_{+}(\sqrt{\lambda}x) &=&
\frac{1}{\sqrt{2\pi} \lambda^{1/4}} 
\cos(\sqrt{\lambda}x),
\nonumber
\end{eqnarray}
can be regarded as the generalized eigenfunctions of the
Hamiltonian
\begin{equation}
{\cal H}_{\rm sin} = 
- \frac{\partial^2}{\partial x^2}
\label{eqn:Hsin1}
\end{equation}
with spectrum $\lambda >0$.
Here $S_-$ is an odd function
(parity ${\sf p}=-$) and
$S_+$ is an even function
(parity ${\sf p}=+$), respectively.
Consider an operator $\widehat{\cal H}_{\rm sin}$
such that $\langle x| \widehat{\cal H}_{\rm sin} | y \rangle
=\delta(x-y) {\cal H}_{\rm sin}$ with
(\ref{eqn:Hsin1}) 
and introduce a set of its eigenvectors
$\{|\lambda, {\sf p} ; {\rm sin} \rangle :
\lambda > 0, {\sf p}=\pm \}$;
$$
\widehat{\cal H}_{\rm sin} |\lambda, {\sf p} ; {\rm sin} \rangle
= \lambda | \lambda, {\sf p} ; {\rm sin} \rangle.
$$
Since 
$\langle x | \lambda, {\sf p} ; {\rm sin} \rangle
=\langle \lambda; {\sf p}; {\rm sin}|x \rangle
= S_{{\sf p}}(\sqrt{\lambda}x)$,
we can confirm the completeness of the set
$$
\sum_{{\sf p}=\pm} \int_{0}^{\infty} d \lambda \,
\langle x | \lambda, {\sf p} ; {\rm sin} \rangle
\langle \lambda, {\sf p} ; {\rm sin} | y \rangle
=\frac{1}{2\pi} \int_{-\infty}^{\infty} d u \,
\cos(u(x-y))
= \delta(x-y).
$$

If we change the variable in the Airy differential equation
(\ref{eqn:Aeq1}) by $x \to u-\lambda$,
we have
$$
-\frac{d^2}{dx^2} {\rm Ai}(x-\lambda)+x {\rm Ai}(x-\lambda)
=\lambda {\rm Ai}(x-\lambda).
$$
That is, the Hamiltonian
\begin{equation}
{\cal H}_{\rm Ai}=-\frac{\partial^2}{\partial x^2} + x
\label{eqn:HAi1}
\end{equation}
has $\R$ as spectrum and 
the Airy functions of the form ${\rm Ai}(x-\lambda)$
are its generalized eigenfunctions.
We can consider the corresponding operator
$\widehat{\cal H}_{\rm Ai}$ and its
eigenvectors $\{|\lambda ; {\rm Ai} \rangle :
\lambda \in \R \}$,
$$
\widehat{\cal H}_{\rm Ai} |\lambda; {\rm Ai} \rangle
= \lambda |\lambda ; {\rm Ai} \rangle,
$$
where
$\langle x | \lambda; {\rm Ai} \rangle 
=\langle \lambda; {\rm Ai}|x \rangle = {\rm Ai}(x-\lambda)$.
We find the completeness
$$
\int_{-\infty}^{\infty} d \lambda \,
\langle x | \lambda ; {\rm Ai} \rangle
\langle \lambda ; {\rm Ai} | y \rangle
=\int_{-\infty}^{\infty} d \lambda \ {\rm Ai}(x-\lambda)
{\rm Ai}(y-\lambda)
= \delta(x-y).
$$
\vskip 0.5cm
\noindent{\bf Remark 3.} \quad
The $2(\nu+1)$-dimensional squared Bessel process
(BESQ$_{\nu}$) $Y^{(\nu)}(t)$ is defined as a
solution of the SDE
\begin{equation}
d Y^{(\nu)}(t)= 2 \sqrt{Y^{(\nu)}(t)} d B(t)
+ 2 (\nu+1) dt, \quad
\nu > -1,
\label{eqn:BESQ1}
\end{equation}
where $B(t)$ is a one-dimensional standard BM
\cite{RY98,BS02}.
Its forward and backward Kolmogorov 
(Fokker-Planck) equations are given as
\begin{equation}
\frac{\partial}{\partial t} u(t, x)
= 2 x \frac{\partial^2}{\partial x^2}
u(t,x) \mp 2(\nu \mp 1) 
\frac{\partial}{\partial x} u(t,x),
\label{eqn:BESQ2}
\end{equation}
where $-$ for the forward and $+$ for the
backward equations, respectively.
We will see that the Laguerre polynomials 
$\{L_{k}^{\nu}(x)\}$ 
play for BESQ$_{\nu}$ the similar role to 
the Hermite polynomials $\{H_{k}(x)\}$ for
the BM shown in Section 2.
By considering the noncolliding system of BESQ$_{\nu}$s
instead of BMs \cite{KO01}, 
we can derive a finite determinantal 
process whose matrix-kernel is described using
the orthonormal Laguerre functions
(extended Laguerre kernel \cite{TW04}).
If we take the so-called hard-edge scaling limit,
a spatially inhomogeneous infinite particle system
is obtained, which is the determinantal process
associated with the matrix-kernel
\begin{equation}
\bK^{(\nu)}(t,y;s,x) = \left\{
   \begin{array}{ll}
\displaystyle{
\int_{0}^{1} d \lambda \,
e^{(s-t)\lambda} J_{\nu}(2 \sqrt{\lambda y})
J_{\nu}(2 \sqrt{\lambda x})
} 
& \mbox{if} \quad t < s  \\
& \\
\displaystyle{
\frac{J_{\nu}(2 \sqrt{y}) \sqrt{x} J_{\nu}'(2 \sqrt{x})
-\sqrt{y} J_{\nu}'(2\sqrt{y}) J_{\nu}(2\sqrt{x})}{y-x}
}
& \mbox{if} \quad t=s \\
& \\
\displaystyle{
- \int_{1}^{\infty} d \lambda \,
e^{-(t-s) \lambda} J_{\nu}(2 \sqrt{\lambda y})
J_{\nu}(2 \sqrt{\lambda x})
}
& \mbox{if} \quad t > s,
   \end{array} \right. 
\label{eqn:KBessel1}
\end{equation}
where $J_{\nu}(x)$ is the Bessel function,
$J_{\nu}(x)=\sum_{n=0}^{\infty} (-1)^{n}(x/2)^{2n+\nu}
/\{\Gamma(n+1) \Gamma(n+1+ \nu)\}$,
and $J_{\nu}'(x)=d J_{\nu}(x)/dx$.
This kernel was called the extended Bessel kernel
in \cite{TW04}. 
See also \cite{Osa04,TW04,KT07,TW07}.
We can see that $J_{\nu}(2 \sqrt{\lambda x})$ is the
generalized eigenfunction of the Hamiltonian
\begin{equation}
{\cal H}_{J}=
- \frac{\partial}{\partial x} x
\frac{\partial}{\partial x}
+ \frac{\nu^2}{4x},
\label{eqn:HJ1}
\end{equation}
with spectrum $\lambda \geq 0$.
We introduce the corresponding operator
$\widehat{\cal H}_{J}$ and its eigenvectors
$\{|\lambda; J \rangle : \lambda \geq 0\}$,
$\widehat{\cal H}_{J}|\lambda; J \rangle
= \lambda |\lambda ; J \rangle$,
where $\langle x | \widehat{\cal H}_{J}
| y \rangle = \delta(x-y){\cal H}_{J}$,
$\langle x| \lambda; J \rangle
=\langle \lambda; J | x \rangle
= J_{\nu}(2 \sqrt{\lambda x})$
with the completeness
$$
\int_{0}^{\infty} d \lambda \,
\langle x | \lambda; J \rangle
\langle \lambda; J | y \rangle
= \int_{0}^{\infty} d \lambda \,
J_{\nu}(2 \sqrt{\lambda x})
J_{\nu}(2 \sqrt{\lambda y})
=\delta(x-y).
$$
\vskip 0.5cm
Here we note the fact that
the sine kernel $\bK$ given by (\ref{eqn:Kbulk1}), 
the Airy kernel $\sK$ given by (\ref{eqn:KAiry1}),
and the Bessel kernel given by (\ref{eqn:KBessel1}) are
expressed in the same way as (\ref{eqn:tKN1}) for 
the Hermite kernel $\widetilde{\mbK}_N$,
$$
K(t,y; s,x) = \left\{
   \begin{array}{ll}
\displaystyle{
\langle y | e^{(s-t) \widehat{\cal H}}
{\cal P} | x \rangle
} 
& \mbox{if} \quad t \leq s  \\
& \\
\displaystyle{
- \langle y|e^{-(t-s) \widehat{\cal H}}
(1-{\cal P}) | x \rangle
}
& \mbox{if} \quad t > s,
   \end{array} \right. 
$$
if we assign the Hamiltonian 
$\delta(x-y) {\cal H}
=\langle x | \widehat{\cal H} | y \rangle$
and projection operator as follows
instead of $\widehat{\cal H}_{\varphi}$ and 
(\ref{eqn:projection1}); 
\begin{description}
\item{(i)} \quad
for the sine kernel, set
${\cal H}={\cal H}_{\rm sin}$ with (\ref{eqn:Hsin1})
and
$\displaystyle{
{\cal P}_{\rm sin} = \sum_{{\sf p}=\pm} \int_0^{1} d \lambda \,
|\lambda, {\sf p} ; {\rm sin} \rangle
\langle \lambda, {\sf p} ; {\rm sin} |,
}$
\item{(ii)} \quad
for the Airy kernel, set
${\cal H}={\cal H}_{\rm Ai}$ with (\ref{eqn:HAi1})
and
$\displaystyle{
{\cal P}_{\rm Ai} = \int_{-\infty}^{0} d \lambda \,
|\lambda ; {\rm Ai} \rangle
\langle \lambda ; {\rm Ai} |.
}$
\item{(iii)} \quad
for the Bessel kernel, set
${\cal H}={\cal H}_{J}$ with (\ref{eqn:HJ1})
and
$\displaystyle{
{\cal P}_{J} = \int_{0}^{1} d \lambda \,
|\lambda ; J \rangle
\langle \lambda ; J |.
}$
\end{description}

\subsection{Effective Hamiltonians and Matrix-Kernels} 

In the previous subsection, we claimed that the structure of
the matrix-kernel of determinantal correlation functions
is common both in finite particle systems and infinite
particle systems.
It should be noted, however, that even from the same
finite system ({\it e.g.} noncolliding BM governed by
the Hamiltonian ${\cal H}_{\varphi}$),
depending on scaling limits, different kinds of
infinite determinantal systems are derived,
in which the Hamiltonian is replaced by
appropriate {\it effective Hamiltonians}
({\it e.g.} 
${\cal H}_{\rm sin}$ and ${\cal H}_{\rm Ai}$)
and spectral projection operator is modified.

In the present subsection, we give a possible
general consideration on the common structure of determinantal
processes.
First we note the following fact for a 
general form of effective Hamiltonian
$
{\cal H}=-a(x) \partial^2/\partial x^2
-b(x) \partial/\partial x
-c(x)$, 
where $a(x), b(x), c(x)$ are sufficiently smooth functions
with $a(x) \not=0$ in an interval $\Lambda \subset \R$.
If we change the variable $x \mapsto z$ following
$$
z(x)=\int_0^{x} a(y)^{-1/2} dy,
$$
${\cal H} \to \widetilde{\cal H}=
-\partial^2/\partial z^2
-\widetilde{b}(z) \partial/ \partial z -c(z)$
with
$\widetilde{b}(z(x))=
a(x) z''(x)+b(x)z'(x)=
a(x)^{-1/2}\{-a'(x)/2+b(x)\}$.
Then we define
$$
r(z)=\exp \left\{
-\frac{1}{2} \int_0^{z} \widetilde{b}
(u) du \right\}
$$
and if we perform a similarity transformation
$\widetilde{\cal H} \to
\overline{\cal H}=r^{-1} \widetilde{\cal H} r$,
the term of first derivative can be eliminated
and we have the form
$\overline{\cal H}=-\partial^2/\partial z^2
+q(z)$.
The transformation ${\cal H} \to \overline{\cal H}$
is called the Liouville transformation.
Then, without loss of generality, 
we can assume effective Hamiltonians of the form
(the Sturm-Liouville operator \cite{Tit62})
\begin{equation}
{\cal H}=-\frac{\partial^2}{\partial x^2}
+q(x)
\label{eqn:Hgen1}
\end{equation}
defined on $\Lambda \subset \R$.
\vskip 0.5cm
\noindent {\bf Example 1.} \quad
The effective Hamiltonians
${\cal H}_{\varphi}, {\cal H}_{\rm sin},
{\cal H}_{\rm Ai}$, and ${\cal H}_{J}$ are
transformed to the form (\ref{eqn:Hgen1})
with the following $q(x)$, respectively
\begin{equation}
q(x) : \quad \frac{1}{16}(x^2-4), \quad
0, \quad x, \quad
\left( \nu^2-\frac{1}{4} \right) \frac{1}{x^2},
\label{eqn:qx}
\end{equation}
where $\Lambda=\R$ for the first three cases and
$\Lambda=\R_{+} \equiv \{x \in \R : x >0\}$
for the last case.
\vskip 0.5cm
Let $\widehat{\cal H}$ be the operator corresponding to
(\ref{eqn:Hgen1});
$\langle x | \widehat{\cal H} | y \rangle
=\delta(x-y){\cal H}$, and define
\begin{equation}
\delta_{t}(x,y) = \delta_{t}(y,x)=
\langle y | e^{-t \widehat{\cal H}} |x \rangle.
\label{eqn:deltat1}
\end{equation}
By definition, it solves the equation
\begin{equation}
\frac{\partial}{\partial t}
\delta_{t}(x,y) = - {\cal H} \delta_{t}(x,y)
\label{eqn:GtA1}
\end{equation}
with
\begin{equation}
\lim_{t \to 0} \delta_{t}(x,y)=\delta(x-y).
\label{eqn:Gd-I}
\end{equation}
Assume that $\widehat{\cal H}$ has a
distribution $\omega(d\lambda)$ 
of spectrum 
$\sigma(\widehat{\cal H})
=\{\lambda: \widehat{\cal H}|\lambda \rangle
=\lambda |\lambda \rangle\}$
with the complete set of eigenvectors
$\{ |\lambda \rangle :
\lambda \in \sigma(\widehat{\cal H}) \}$.
Then the spectral representation of
$\delta_t(x,y)$ is given by
\begin{eqnarray}
\delta_{t}(x,y) &=& \langle y | e^{-t
\widehat{\cal H}}
\int_{\sigma(\widehat{\cal H})} \omega(d \lambda)
|\lambda \rangle \langle \lambda | x \rangle
\nonumber\\
&=& \int_{\sigma(\widehat{\cal H})} \omega(d \lambda)
e^{-\lambda t} \Phi_{\lambda}(x)
\Phi_{\lambda}(y)
\nonumber
\end{eqnarray}
with $\Phi_{\lambda}(x)=\langle x | \lambda \rangle
=\langle \lambda | x \rangle$.
We assume that 
\begin{equation}
^{\exists}C > 0 \quad
\mbox{s.t.} \quad 
\int_{\Lambda} dy \,
(x-y)^{4} \delta_{t}(x,y)
\leq C t^{2}, \quad
t \in (0, 1], \quad
^{\forall}x \in \Lambda,
\label{eqn:Bound1}
\end{equation}
where $C$ does not depend of $t$.

\vskip 0.5cm
\noindent {\bf Example 2.} \quad
For the four examples (\ref{eqn:qx}), we have 
the following explicit expressions of
$\delta_t(x,y)$;
\begin{eqnarray}
q(x)=\frac{1}{16}(x^2-4): \quad && \nonumber\\
\delta_{t}(x,y) &=& 
\frac{1}{2} \sum_{n=0}^{\infty}
\varphi_n \left(\frac{x}{2}\right)
\varphi_{n} \left( \frac{y}{2} \right)
e^{-nt/2}
\nonumber\\
&=& 
\frac{1}{2 \sqrt{\pi (1-e^{-t})}}
\exp \left\{
-\frac{1}{8}(x-y)^2 \coth \left(\frac{t}{2}\right)
-\frac{1}{4} xy \tanh \left(\frac{t}{4}\right) \right\},
\nonumber\\
&& \hskip 6cm x, y \in \Lambda=\R,
\nonumber\\
 q(x)=0: \hskip 2cm && \nonumber\\
\delta_{t}(x,y) &=& 
\sum_{{\sf p} = \pm} \int_{0}^{\infty} d \lambda \,
\langle y | e^{-t \widehat{\cal H}_{\rm sin}} |
\lambda, {\sf p} ; {\rm sin} \rangle
\langle \lambda, {\sf p} ; {\rm sin} | x \rangle
\nonumber\\
&=& \frac{1}{2\pi} \int_{0}^{\infty}
\frac{d \lambda}{\sqrt{\lambda}} e^{-\lambda t}
\Big\{ \cos(\sqrt{\lambda} x) \cos(\sqrt{\lambda} y) 
+ \sin(\sqrt{\lambda} x) \sin(\sqrt{\lambda} y) \Big\}
\nonumber\\
&=& \frac{1}{\pi} \int_{0}^{\infty}
du \, e^{-u^2 t} \cos (u(x-y))
\nonumber\\
&=& \frac{1}{\sqrt{4 \pi t}}
\exp \left\{- \frac{(x-y)^2}{4t} \right\},
\quad x, y \in \Lambda = \R,
\nonumber\\
 q(x)=x : \hskip 2cm && \nonumber\\
\delta_{t}(x,y) &=&
\int_{-\infty}^{\infty} 
\Ai(x-\lambda) \Ai(y-\lambda) e^{-\lambda t} d \lambda
\nonumber\\ 
&=& \frac{1}{\sqrt{ 4 \pi t}}
\exp \left\{ -\frac{(x-y)^2}{4t} 
-\frac{t(x+y)}{2}+\frac{t^3}{12} \right\},
\quad x, y \in \Lambda=\R,
\nonumber\\
q(x)=\left( \nu^2-\frac{1}{4} \right) \frac{1}{x^2} : 
&& \nonumber\\
\delta_{t}(x,y) &=& 
\frac{1}{2} \int_{0}^{\infty} \sqrt{x} 
J_{\nu} (\sqrt{\lambda} x)
\sqrt{y} J_{\nu}(\sqrt{\lambda} y) e^{-\lambda t} d \lambda
\nonumber\\
&=& 
\frac{\sqrt{xy}}{2t} e^{-(x^2+y^2)/4t}
I_{\nu}\left(\frac{xy}{2t}\right),
\quad \nu > -1, \quad x , y \in \Lambda=\R_+,
\nonumber
\end{eqnarray}
where $I_{\nu}(z)$ is the modified Bessel function
given by
$
I_{\nu}(z)=
\sum_{n=0}^{\infty} (z/2)^{2n+\nu}/
\{n! \Gamma(\nu+n+1)\}.
$
We can confirm that (\ref{eqn:Bound1}) is satisfied
in these four cases.
\vskip 0.5cm

Now we consider a subset of 
spectrum $\sigma(\widehat{\cal H})$,
$\sigma_{-}(\widehat{\cal H})=\{\lambda \in \sigma(\widehat{\cal H}) :
\lambda \leq \lambda_{*} \}$
with a specified level $\lambda_{*} \in \sigma(\widehat{\cal H})$,
and define the projection operator onto $\sigma_{-}(\widehat{\cal H})$
\begin{equation}
{\cal P} = \int_{\sigma_{-}(\widehat{\cal H})} 
\omega(d \lambda) |\lambda \rangle
\langle \lambda |.
\label{eqn:projgen}
\end{equation}
so that
\begin{eqnarray}
\sG_{t}(x,y) &=& \sG_{t}(y,x) 
= \langle y | e^{t \widehat{\cal H}} 
{\cal P} | x \rangle,
\nonumber\\
\overline{\sG}_{t}(x,y) &=& \overline{\sG}_{t}(y,x) 
= - \langle y | e^{-t \widehat{\cal H}} 
(1-{\cal P}) | x \rangle,
\label{eqn:sG1}
\end{eqnarray}
and 
\begin{eqnarray}
\label{eqn:G-K}
K(x,y) &=& \lim_{t \to 0} \sG_t(x,y)
= \lim_{t \to 0} \sG_{-t}(x,y) 
= \langle y | {\cal P} |x \rangle, \\
\label{eqn:rhox1}
\rho(x) &=& K(x,x) 
= \langle x | {\cal P}|x \rangle.
\end{eqnarray}
By definition 
\begin{equation}
\overline{\sG}_t(x,y)=\sG_{-t}(x,y)-
\delta_{t}(x,y),
\label{eqn:expandG}
\end{equation}
\begin{eqnarray}
&& \frac{\partial}{\partial t} \sG_{t}(x,y)
= {\cal H} \sG_{t}(x,y), \nonumber\\
&& \frac{\partial}{\partial t} \overline{\sG}_{t}(x,y)
= -{\cal H} \overline{\sG}_{t}(x,y),
\label{eqn:Gt2}
\end{eqnarray}
and
\begin{equation}
\lim_{t\to 0} \frac{\partial^2}{\partial y^2} \sG_t(x,y)
= \lim_{t\to 0} \frac{\partial^2}{\partial y^2} \sG_{-t}(x,y)
= \frac{\partial^2}{\partial y^2} K(x,y).
\label{eqn:dG-dG}
\end{equation}
Moreover, by the completeness of $\{|x \rangle : x \in \Lambda \}$, 
$\int_{\Lambda} dx |x \rangle \langle x|=1$,
we have the relations
\begin{eqnarray}
&&\int_{\Lambda} dy \
\delta_t(x,y) \sG_t(y,x) = \rho (x),
\label{eqn:DG=rho}
\\
&&\int_{\Lambda} dy \
\delta_t(x,y) \sG_t(y,z) = K (x,z).
\label{eqn:DG=K}
\end{eqnarray}

We define matrix-kernel $K$ by 
\begin{equation}
K(t,y; s,x) = \left\{
   \begin{array}{ll}
\displaystyle{
\sG_{s-t}(y,x)
} 
& \mbox{if} \quad t \leq s  \\
& \\
\displaystyle{
\overline{\sG}_{t-s}(y,x)
}
& \mbox{if} \quad t > s,
   \end{array} \right. 
\label{eqn:Kgen1}
\end{equation}
and consider the determinantal process
$\Xi(t)=\{\X(t)\}$, whose 
multitime correlation function is given by
\begin{equation}
\rho(0, \xi_{N_0}^{(0)}; t_1, \xi_{N_1}^{(1)};
\cdots; t_M, \xi_{N_M}^{(M)})
=\det_{1 \leq j \leq N_{m},
1 \leq k \leq N_{n},
0 \leq m,n \leq M}
\Bigg[ K(t_m, x_j^{(m)}; t_n, x_k^{(n)}) \Bigg]
\label{eqn:rhogen1}
\end{equation}
for any $M \geq 0$, any sequence $\{N_m \}_{m=0}^{M}$
of positive integers, and any series of 
observation times
$0 < t_1 < \cdots < t_M$.

The invariant measure of the process is
the determinantal point field $\mu$ associated with $K(x,y)$
given by (\ref{eqn:G-K}).

The two-time correlation function 
$\cR(0,\{\x_m \} ; t, \{\y_n \})$ 
of the system is given by 
$$
\cR(0,\{\x_{m}\}; t,\{\y_n\})
= \det \overline{\bf M}(t, \y_{n}| \x_{m}),
\quad t >0, \, m,n \geq 0,
$$
where
\begin{equation}
\overline{\bf M}(t, \y_{n}| \x_{m})
= \left( \matrix{
\rho(x_{1}) & \cdots & K(x_{1}, x_{m}) &
\sG_{t}(x_{1}, y_{1}) & \cdots
& \sG_{t}(x_{1}, y_{n}) \cr
K(x_{2}, x_{1},) & \cdots & K(x_{2}, x_{m}) &
\sG_{t}(x_{2}, y_{1}) & 
\cdots & \sG_{t}(x_{2}, y_{n}) \cr
   &   \cdots &  & & \cdots &  \cr
K(x_{m}, x_{1}) 
& \cdots & \rho(x_{m}) & 
\sG_{t}(x_{m}, y_{1}) & \cdots
& \sG_{t}(x_{m}, y_{n}) \cr
\overline{\sG}_{t}(y_{1}, x_{1}) & \cdots & 
\overline{\sG}_{t}(y_{1}, x_{m}) & 
\rho(y_{1}) & 
\cdots & K(y_{1}, y_{n}) \cr
\overline{\sG}_{t}(y_{2}, x_{1}) & \cdots &
\overline{\sG}_{t}(y_{2}, x_{m}) & K(y_{2}, y_{1}) 
& \cdots & K(y_{2}, y_{n}) \cr
   &   \cdots &  & & \cdots & \cr
\overline{\sG}_{t}(y_{n}, x_{1}) & \cdots 
& \overline{\sG}_{t}(y_{n}, x_{m}) 
& K(y_{n}, y_{1}) & 
\cdots & \rho (y_{n}) 
} \right). \nonumber\\
\end{equation}
For a matrix $A=(a_{ij})_{i \in I, j \in J}$
with index sets $I, J$,
we denote its submatrix as
$A_{I' J'} \equiv (a_{ij})_{i \in I', j \in J'}$
for $I' \subset I, J' \subset J$,
and the complementary submatrix as
$A^{I' J'} \equiv (a_{ij})_{i \in I \setminus I', 
j \in J \setminus J'}$.
By using the relation (\ref{eqn:expandG}), 
$\cR(0,\{\x_{m}\};t,\{\y_{n}\})$ is expanded as
\begin{eqnarray}
\cR(0,\{\x_{m}\};t, \{\y_{n}\})
&=& \det {\bf M}(t, \y_{n}| \x_{m})
\nonumber\\
&+& \sum_{\ell=1}^{m\wedge n} 
\sum_{1\leq a_1<\cdots<a_\ell \leq m}
\sum_{1\leq b_1<\cdots<b_\ell \leq n}
(-1)^{\ell+\sum_{i=1}^\ell (a_i+m+b_i)}
\nonumber\\
&& \qquad \qquad \times
\det {\bf D}(t, \y_{n}|\x_{m})_{\{\a_\ell\} \{\b_\ell\}}
\det {\bf M}(t, \y_{n}|\x_{m})^{\{m+\b_\ell\} \{\a_\ell\}},
\nonumber\\
\label{eqn:expansion1}
\end{eqnarray}
where
\begin{eqnarray}
{\bf M}(t, \y_{n}| \x_{m})
&=& \left( \matrix{
\rho(x_{1}) & \cdots & K(x_{1}, x_{m}) &
\sG_{t}(x_{1}, y_{1}) & \cdots
& \sG_{t}(x_{1}, y_{n}) \cr
K(x_{2}, x_{1},) & \cdots & K(x_{2}, x_{m}) &
\sG_{t}(x_{2}, y_{1}) & 
\cdots & \sG_{t}(x_{2}, y_{n}) \cr
   &   \cdots &  & & \cdots &  \cr
K(x_{m}, x_{1}) 
& \cdots & \rho(x_{m}) & 
\sG_{t}(x_{m}, y_{1}) & \cdots
& \sG_{t}(x_{m}, y_{n}) \cr
\sG_{-t}(y_{1}, x_{1}) & \cdots & 
\sG_{-t}(y_{1}, x_{m}) & 
\rho(y_{1}) & 
\cdots & K(y_{1}, y_{n}) \cr
\sG_{-t}(y_{2}, x_{1}) & \cdots &
\sG_{-t}(y_{2}, x_{m}) & K(y_{2}, y_{1}) 
& \cdots & K(y_{2}, y_{n}) \cr
   &   \cdots &  & & \cdots & \cr
\sG_{-t}(y_{n}, x_{1}) & \cdots 
& \sG_{-t}(y_{n}, x_{m}) 
& K(y_{n}, y_{1}) & 
\cdots & \rho (y_{n}) 
} \right), \nonumber\\
{\bf D}(t, \y_{n}|\x_{m})
&=& \Big( \delta_{t}(x_i, y_j) \Big)_{1 \leq i \leq m, 
1 \leq j \leq n}.
\nonumber
\end{eqnarray}

It was shown by Shirai and Takahashi \cite{ST03} that 
the Palm measure $\mu^{z}$ coincides with the determinantal 
point field
associated with the kernel $K^{z}$ defined by  
$$
K^z(x,y)=\frac{1}{K(z,z)}
\det\left( \matrix{
K(x,y) & K(x,z) 
\cr
K(z,y)& K(z,z)
} \right).
$$
Note that
\begin{equation}
\rho^z(\{\x_m\})
=\frac{1}{\rho(z)}\rho(\{\x_m\} \cup \{z\}).
\label{eqn:palm}
\end{equation}

For $f,g \in C_0(\R)$, let
\begin{equation}
\langle f, \xi \rangle = \sum_{x \in \xi}f(x)
\quad \mbox{and} \quad
\langle g, \eta \rangle = \sum_{y \in \eta} g(y),
\label{eqn:bracketA}
\end{equation}
and set $\alpha, \beta \in {\bf C}$.
Define
$F(\xi)= e^{\alpha \langle f, \xi \rangle}$,
$G(\xi)=e^{\beta \langle g, \xi \rangle}$.
Then the two-time generating function is
defined by 
\begin{equation}
\Phi_{t}(f,g; \alpha, \beta) 
\equiv {\bf E} \Big[F( \Xi(0)) G(\Xi(t)) \Big] 
= {\bf E} \Big[
e^{\alpha \langle f, \Xi(0) \rangle +
\beta \langle g, \Xi(t) \rangle} \Big],
\label{eqn:ttchara0}
\end{equation}
for $t>0$.
Let 
$$
\chi_{1}(x)=e^{\alpha f(x)}-1, \quad
\chi_{2}(y)=e^{\beta g(y)}-1,
$$
then we can see
\begin{equation}
\Phi_{t}(f,g; \alpha, \beta) 
= \sum_{m=0}^{\infty} \sum_{n=0}^{\infty} 
\frac{1}{m! n!} 
\int_{{\Lambda}^{m}} d \x_m \int_{{\Lambda}^{n}} d \y_n \
\prod_{j=1}^{m} \chi_{1}(x_{j})
\prod_{k=1}^{n} \chi_{2}(y_{k})
\cR(0, \{\x_{m}\};t, \{\y_{n}\}). 
\label{eqn:ttchara1}
\end{equation}

\SSC{CHARACTERIZATION OF DETERMINANTAL PROCESSES}
In the previous section, we introduced a class of
determinantal processes associated with spectral projections
defined by effective Hamiltonians.
Here we give properties of determinantal processes
of this class, which can be derived from the
common structure of correlation functions.

\subsection{Continuity}
Let $C_0^{\infty}(\R)$ be the set of all 
infinitely differentiable real functions with
compact supports and set
$\langle f, \xi \rangle \equiv 
\sum_{x \in \xi} f(x)$ for
$f \in C_0^{\infty}(\R), \xi \in \mX$.
By a criterion of Kolmogorov 
(see, for example,\cite{Kal97,Bil99}), 
the following lemma implies that
$\sum_{x \in \Xi(t)} f(x)$ is continuous
in time with probability one for any 
$f \in C_{0}^{\infty}(\R)$.
Since $\mX$ is separable with the vague topology,
we can choose a countable set
$\{f_j\}_{j=1}^{\infty} \subset C_0^{\infty}(\R)$ such that
$\xi_n \to \xi$ in $n \to \infty$ on $\mX$
if and only if 
$\langle \xi_n, f_j \rangle 
\to \langle \xi, f_j \rangle$ in 
$n \to \infty, \, \forall j \geq 1$.
Then it implies that the determinantal process
$\Xi(t)$ is continuous in the vague topology 
with probability one.
That is, if the condition (\ref{eqn:Bound1}) 
is satisfied,
it can be expressed by
$$
\Xi(t)=\sum_{j=1}^{\infty} 
\delta_{X_j(t)}
$$
with some real-valued continuous
processes $X_{j}(t), j \in \N$.

\begin{lem}
\label{lem:continuity}
Let $\Xi(t)$ be the determinantal process,
whose multitime correlation functions are
given by (\ref{eqn:rhogen1}) 
with (\ref{eqn:sG1}) and (\ref{eqn:Kgen1})
associated with an effective Hamiltonian 
(\ref{eqn:Hgen1}) defined on $\Lambda \subset \R$.
Assume that (\ref{eqn:Bound1}) is satisfied. 
Then for any $f \in C_0^{\infty}(\R)$
\begin{equation}
\E \Bigg[
\Big|\langle f, \Xi(t)\rangle
-\langle f, \Xi(0)\rangle \Big|^4 \Bigg]
\leq C t^2,
\quad t \in (0,1],
\label{eqn:KolCri}
\end{equation}
where $C$ does not depend on $t$.
\end{lem}
\vskip 0.5cm
\noindent{\bf Remark 4.} \quad
Theorems \ref{thm:infinite_sin} and \ref{thm:infinite_Ai}
are the limit theorems of the processes
\begin{eqnarray}
&&\Xi_{N}^{\rm bulk}(t)
\equiv \Big\{ X_1(N+2t), \cdots,
X_N(N+2t) \Big\}
\, \Longrightarrow \,
\Xi^{\rm sin}(t)
\nonumber\\
&&\Xi_{N}^{\rm edge}(t)
\equiv \Big\{ X_1(N^{1/3}+2t)-a_N(t), \cdots,
X_N(N^{1/3}+2t)-a_N(t) \Big\}
\, \Longrightarrow \,
\Xi^{\rm Ai}(t), \quad
N \to \infty
\nonumber
\end{eqnarray}
with $a_N(t)$ defined by (\ref{eqn:aNs}),
in the sense of finite-dimensional distributions,
where the multitime correlation functions
of $\Xi^{\rm sin}(t)$ and $\Xi^{\rm Ai}(t)$
are given by $\rho_{\rm sin}$ of (\ref{eqn:bulk1})
and $\rho_{\rm Ai}$ of (\ref{eqn:edge1}),
respectively.
When we apply this lemma to the 
$N$-particle system of noncolliding BM,
$\Xi^{\bf X}_{N}(t)$, which is a determinantal process
associated with the Hermite kernel
$\widetilde{\mbK}_N$ given by (\ref{eqn:tKN1}),
$C$ depends on $N$. 
Careful estimation gives an upper bound
on $C$, which is uniform in $N$,
and tightness is established both in
$\{ \Xi^{\rm bulk}_N(t) \}_{N \in \N}$
and $\{ \Xi^{\rm edge}_N(t) \}_{N \in \N}$.
\vskip 0.5cm

\noindent{\it Proof of Lemma \ref{lem:continuity}.} \quad
Here we use the notation
$$
\cR_{m,n}(\x_{m+n})=
\cR(0, \{x_1, \cdots, x_m \}; t, \{x_{m+1},\dots, x_{m+n}\}).
$$
The left-hand-side of (\ref{eqn:KolCri}) is given by
$\partial^4 
\Phi_t(-f,f;\alpha, \alpha)/ \partial^4 \alpha|_{\alpha=0}$
from the two-time generating function given by
(\ref{eqn:ttchara0}) and (\ref{eqn:ttchara1}),
which equals
\begin{eqnarray}
&&\int_{{\Lambda}^4}d\x_4 \ \prod_{i=1}^4 f(x_i)
\Big\{ \cR_{4,0}(\x_4) - 4 \cR_{3,1}(\x_4) + 6 \cR_{2,2}(\x_4) 
-4\cR_{1,3}(\x_4) + \cR_{0,4}(\x_4) \Big\}
\nonumber\\
&+&\int_{{\Lambda}^3}d\x_3 \ \prod_{i=1}^3 f(x_i)
\Big\{ 2\sum_{j=1}^3 f(x_j) \cR_{3,0}(\x_3) 
+6 (-f(x_1)-f(x_2)+f(x_3)) \cR_{2,1}(\x_3) 
\nonumber\\
&&\qquad +6 (f(x_1)-f(x_2)-f(x_3)) \cR_{1,2}(\x_3) 
+2\sum_{j=1}^3 f(x_j) \cR_{0,3}(\x_3) \Big\}
\nonumber\\
&+&\int_{{\Lambda}^2}d\x_2 \ \prod_{i=1}^2 f(x_i)
\Big\{  (2f(x_1)^2 + 3f(x_1)f(x_2)+2f(x_2)^2)(\cR_{2,0}(\x_2) 
+ \cR_{0,2}(\x_2))
\nonumber\\
&&\qquad +(-4f(x_1)^2 + 6f(x_1)f(x_2)-4f(x_2)^2)\cR_{1,1}(\x_2) \Big\}
\nonumber\\
&+&\int_{\Lambda}d x_1 \ f(x_1)^4 \Big\{\cR_{1,0}(x_1) + \cR_{0,1}(x_1)
\Big\}.
\nonumber
\end{eqnarray}
Since 
$
\cR_{3,0}(x_1,x_2,x_3) = \cR_{0,3}(x_1,x_2,x_3),
\cR_{2,1}(x_1,x_2,x_3) = \cR_{1,2}(x_3,x_2,x_1),
$
{\it etc.},
the above quantity is twice of
\begin{eqnarray}
I &\equiv& \int_{{\Lambda}^4}d\x_4 \ \prod_{i=1}^4 f(x_i)
\Big\{ \cR_{4,0}(\x_4) -4 \cR_{3,1}(\x_4) + 3 \cR_{2,2}(\x_4) \Big\}
\nonumber\\
&+&\int_{{\Lambda}^3}d\x_3 \ \prod_{i=1}^3 f(x_i)
\Big\{ 2\sum_{j=1}^3 f(x_j) \cR_{3,0}(\x_3) 
+6 (-f(x_1)-f(x_2)+f(x_3)) \cR_{2,1}(\x_3) \Big\}
\nonumber\\
&+&\int_{{\Lambda}^2}d\x_2 \ \prod_{i=1}^2 f(x_i)
\Big\{  (2f(x_1)^2 + 3f(x_1)f(x_2)+2f(x_2)^2)\cR_{2,0}(\x_2)
\nonumber\\
&&\qquad +(-2f(x_1)^2 + 3f(x_1)f(x_2)-2f(x_2)^2)\cR_{1,1}(\x_2)  
\Big\}
\nonumber\\
&+&\int_{\Lambda}d x_1 \ f(x_1)^4 \cR_{1,0}(x_1).
\nonumber
\end{eqnarray}
We put 
${\bf M}_{m,n}(\x_{m+n})
={\bf M}(t, x_{m+1},\dots, x_{m+n}| x_1,\dots,x_m)$,
${\bf D}(\x_n)={\bf D}(t, \x_n|\x_n)$
and
\begin{eqnarray}
&&\cD_{1,1}(\x_2) 
= -\delta_t(x_2,x_1)\sG_t(x_1,x_2),
\nonumber\\
&&{ \cD}_{2,1}(\x_3) 
= \sum_{i=1}^2 (-1)^{i+3} \delta_t(x_3,x_i) 
\det {\bf M}_{2,1}(\x_3)^{\{3\} \{i\}},
\nonumber\\
&&{ \cD}_{3,1}(\x_4) 
= \sum_{i=1}^3 (-1)^{i+4} \delta_t(x_4,x_i) 
\det {\bf M}_{3,1}(\x_4)^{\{4\} \{i\}},
\nonumber\\
&&{ \cD}_{2,2}(\x_4) = \sum_{i=1}^2 \sum_{j=3}^4 (-1)^{i+j} 
\delta_t(x_j,x_i) \det {\bf M}_{2,2}(\x_4)^{\{j\} \{i\}},
\nonumber\\
&&\widehat{\cD}_{2,2}(\x_4) 
=\det {\bf D}(\x_4)_{\{3,4\} \{1,2\}} 
\det {\bf M}_{2,2}(\x_4)^{\{3,4\} \{1,2\}}.
\nonumber
\end{eqnarray}
From (\ref{eqn:expansion1}) we have 
\begin{eqnarray}
&&\cR_{1,0}(x_1)=\cG_{1,0}(x_1),
\nonumber\\
&&\cR_{1,1}(\x_2)
=\cG_{1,1}(\x_2) - \cD_{1,1}(\x_2),
\nonumber\\
&&\cR_{2,1}(\x_3)
= \cG_{2,1}(\x_3) - \cD_{2,1}(\x_3),
\nonumber\\
&&\cR_{3,1}(\x_4)
= \cG_{3,1}(\x_4)- \cD_{3,1}(\x_4),
\nonumber\\
&& \cR_{2,2}(\x_4)
= \cG_{2,2}(\x_4) - \cD_{2,2}(\x_4) + \widehat{\cD}_{2,2}(\x_4).
\nonumber
\end{eqnarray}
We divide $I$ into four terms $I=\sum_{j=1}^{4} I_j$ with
\begin{eqnarray}
I_1 &\equiv& \int_{\Lambda}d x_1 \ f(x_1)^4 \cG_{1,0}(x_1)
\nonumber\\
&&\quad - \int_{{\Lambda}^2}d\x_2 \ \prod_{i=1}^2 f(x_i)
(-2f(x_1)^2 + 3f(x_1)f(x_2)-2f(x_2)^2)\cD_{1,1}(\x_2),
\nonumber\\
I_2 &\equiv&\int_{{\Lambda}^2}d\x_2 \ \prod_{i=1}^2 f(x_i)
\Big\{(2f(x_1)^2 + 3f(x_1)f(x_2)+2f(x_2)^2)\cG_{2,0}(\x_2)
\nonumber\\
&&\qquad +(-2f(x_1)^2 + 3f(x_1)f(x_2)-2f(x_2)^2)\cG_{1,1}(\x_2) 
\Big\}
\nonumber\\
&&\quad - \int_{{\Lambda}^3}d\x_3 \ \prod_{i=1}^3 f(x_i)
6 (-f(x_1)-f(x_2)+f(x_3)) \cD_{2,1}(\x_3)
\nonumber\\
&&\quad +\int_{{\Lambda}^4}d\x_4 \ \prod_{i=1}^4 f(x_i)3 
\widehat{\cD}_{2,2}(\x_4),
\nonumber\\
I_3 &\equiv&\int_{{\Lambda}^3}d\x_3 \ \prod_{i=1}^3 f(x_i)
\Big\{ 2\sum_{j=1}^3 f(x_j) \cG_{3,0}(\x_3) 
+6 (-f(x_1)-f(x_2)+f(x_3)) \cG_{2,1}(\x_3) \Big\}
\nonumber\\
&&\quad - \int_{{\Lambda}^4}d\x_4 \ \prod_{i=1}^4 f(x_i)
\Big\{ -4 \cD_{3,1}(\x_4) + 3 \cD_{2,2}(\x_4)  \Big\},
\nonumber\\
I_4 &\equiv& \int_{{\Lambda}^4}d\x_4 \ \prod_{i=1}^4 f(x_i)
\Big\{ \cG_{4,0}(\x_4) -4 \cG_{3,1}(\x_4) + 3 \cG_{2,2}(\x_4)  
\Big\}.
\nonumber
\end{eqnarray}

Using (\ref{eqn:DG=rho}) we have
\begin{eqnarray}
I_1 &=& -2 \int_{{\Lambda}^2} 
f(x_1)f(x_2)^3 \sG_t(x_1,x_2)\delta_t(x_1,x_2)dx_1dx_2
\nonumber\\
&& \quad 
-2 \int_{{\Lambda}^2} 
f(x_1)^3 f(x_2) \sG_t(x_1,x_2)\delta_t(x_1,x_2)dx_1dx_2
\nonumber
\\
&& \quad +3 \int_{{\Lambda}^2} 
f(x_1)^2f(x_2)^2 \sG_t(x_1,x_2)\delta_t(x_1,x_2)dx_1dx_2
+  \int_{\Lambda} f(x_1)^4 \rho(x_1)dx_1
\nonumber
\\
&=&\frac{1}{2} \int_{{\Lambda}^2} 
\Big\{ f(x_1)-f(x_2) \Big\}^4 
\sG_t(x_1,x_2)\delta_t(x_1,x_2)dx_1dx_2.
\nonumber
\end{eqnarray}
By the assumption (\ref{eqn:Bound1}) we have
\begin{equation}
I_1 \leq C_1 t^2, \quad t\in (0,1],
\label{eqn:boundI1}
\end{equation}
where $C_1$ does not depend on $t$.

Since ${\sG}_{t}(x,y)=K(x,y)
+\partial {\sG}_t(x,y)/\partial t |_{t=0} \, t
+{\cal O}(t^2)$, we have
for any $k, \ell \in \N$
$$
\det {\bf M}_{k, \ell}(\x_{k+\ell})
= \det {\bf M}_{k-1, \ell+1}(\x_{k+\ell})
+{\cal O}(t^2).
$$
Then
$$
I_{j}=\tilde{I}_{j}+{\cal O}(t^2),
\quad j=2,3,4,
$$
where
\begin{eqnarray}
\tilde{I}_2 &=&\int_{{\Lambda}^2}d\x_2 \ \prod_{i=1}^2 f(x_i)
6f(x_1)f(x_2)\cG_{2,0}(\x_2)
\nonumber\\
&&- \int_{{\Lambda}^3}d\x_3 \ \prod_{i=1}^3 f(x_i)
6 (-f(x_1)-f(x_2)+f(x_3)) \cD_{2,1}(\x_3)
+ \int_{{\Lambda}^4}d\x_4 \ \prod_{i=1}^4 f(x_i)3 
\widehat{\cD}_{2,2}(\x_4),
\nonumber\\
\tilde{I}_3 &=& - \int_{{\Lambda}^4}d\x_4 \ \prod_{i=1}^4 f(x_i)
\{ -4 \cD_{3,1}(\x_4) + 3 \cD_{2,2}(\x_4)  \},
\nonumber
\end{eqnarray}
and $\tilde{I}_{4}=0$.
Since the estimate (\ref{eqn:boundI1}) was obtained,
it is enough to show that 
\begin{equation}
\tilde{I}_j \leq C_{j} t^2,
\quad t \in (0, 1], \quad j=2, 3
\label{eqn:boundIj}
\end{equation}
for the proof, 
where $C_j$'s do not depend on $t$.
Using eqs. (\ref{eqn:DG=rho}) and (\ref{eqn:DG=K}),
we obtain
\begin{eqnarray}
\tilde{I}_{2} &=&\int_{{\Lambda}^4}d\x_4 \ 
\left| \matrix{
\sG_{t}(x_{4}, x_1) & \sG_{t}(x_{4}, x_2) \cr
\sG_{t}(x_3, x_{1}) & \sG_{t}(x_{3}, x_2)
} \right|
\Bigg\{ 6f(x_1)^2 f(x_2)^2 \delta_t(x_1,x_4)\delta_t(x_2,x_3)
\nonumber\\
&&\quad + 3 f(x_1)f(x_2)f(x_3)f(x_4) 
\left| \matrix{
\delta_{t}(x_{4}, x_1) & \delta_{t}(x_{4}, x_2) \cr
\delta_{t}(x_3, x_{1}) & \delta_{t}(x_{3}, x_2)
} \right|
\nonumber\\
&& \left. + 6 f(x_1)f(x_2)f(x_3)(-f(x_1)-f(x_2)+f(x_3)) 
\left| \matrix{
\delta_{t}(x_{4}, x_1) & \delta_{t}(x_{4}, x_2) \cr
\delta_{t}(x_3, x_{1}) & \delta_{t}(x_{3}, x_2)
} \right| \right\}
\nonumber\\
&=&\int_{{\Lambda}^4}d\x_4 \ 
\left| \matrix{
\sG_{t}(x_{4}, x_1) & \sG_{t}(x_{4}, x_2) \cr
\sG_{t}(x_3, x_{1}) & \sG_{t}(x_{3}, x_2)
} \right|
\left| \matrix{
\delta_{t}(x_{4}, x_1) & \delta_{t}(x_{4}, x_2) \cr
\delta_{t}(x_3, x_{1}) & \delta_{t}(x_{3}, x_2)
} \right|
3F(\x_4),
\nonumber
\end{eqnarray}
where
$$
F(\x_4)= f(x_1)f(x_2) \Big\{
f(x_1)f(x_2)+f(x_3)f(x_4)-2f(x_1)f(x_3)-2f(x_2)f(x_3)+2f(x_3)^2
\Big\}.
$$
By simple calculation we see that
\begin{eqnarray}
&&F(x_1,x_2,x_3,x_4) + F(x_1,x_2,x_4,x_3)
+F(x_3,x_4,x_1,x_2) + F(x_3,x_4,x_2,x_1)
\nonumber\\
&&=-2(f(x_1)-f(x_3))(f(x_3)-f(x_2))(f(x_2)-f(x_4))(f(x_4)-f(x_1)).
\nonumber
\end{eqnarray}
Then 
\begin{eqnarray}
&&\int_{{\Lambda}^4}d\x_4 \ 
\left| \matrix{
\sG_{t}(x_{4}, x_1) & \sG_{t}(x_{4}, x_2) \cr
\sG_{t}(x_3, x_{1}) & \sG_{t}(x_{3}, x_2)
} \right|
\left| \matrix{
\delta_{t}(x_{4}, x_1) & \delta_{t}(x_{4}, x_2) \cr
\delta_{t}(x_3, x_{1}) & \delta_{t}(x_{3}, x_2)
} \right|
F(\x_4)
\nonumber\\
&=&\int_{{\Lambda}^4}d\x_4 \ \delta_{t}(x_{4}, x_1)\delta_{t}(x_{3}, x_2)
(f(x_1)-f(x_4))(f(x_2)-f(x_3))
\nonumber\\
&&
\times\left| \matrix{
\sG_{t}(x_{4}, x_1) & \sG_{t}(x_{4}, x_2) \cr
\sG_{t}(x_3, x_{1}) & \sG_{t}(x_{3}, x_2)
} \right|
(f(x_1)-f(x_3))(f(x_2)-f(x_4))
\nonumber\\
&=&\int_{{\Lambda}^4}d\x_4 \ \delta_{t}(x_{4}, x_1)\delta_{t}(x_{3}, x_2)
(f(x_1)-f(x_4))(f(x_2)-f(x_3))
\nonumber\\
&&\times\Big\{ 
(\sG_{t}(x_{4}, x_1)\sG_{t}(x_{3}, x_2)-\sG_{t}(x_1, x_2)\sG_{t}(x_{3}, x_4))
(f(x_1)-f(x_4))(f(x_2)-f(x_3))
\nonumber\\
&&\qquad \qquad -\sG_{t}(x_{4}, x_1)\sG_{t}(x_{3}, x_2)
(f(x_1)f(x_4)+f(x_2)f(x_3))
\nonumber\\
&&\qquad \qquad +\sG_{t}(x_{2}, x_4)\sG_{t}(x_{1}, x_3)
(f(x_1)f(x_4)+f(x_2)f(x_3))
\Big\}
\nonumber\\
&=&\int_{{\Lambda}^4}d\x_4 \ \delta_{t}(x_{4}, x_1)\delta_{t}(x_{3}, x_2)
(f(x_1)-f(x_4))(f(x_2)-f(x_3))
\nonumber\\
&&\times\Big\{ 
(\sG_{t}(x_{4}, x_1)\sG_{t}(x_{3}, x_2)-\sG_{t}(x_1, x_2)\sG_{t}(x_{3}, x_4))
(f(x_1)-f(x_4))(f(x_2)-f(x_3))
\nonumber\\
&&+(\sG_{t}(x_{2}, x_4)\sG_{t}(x_{1}, x_3)
-\sG_{t}(x_{2}, x_1)\sG_{t}(x_{4}, x_3))
(f(x_1)f(x_4)+f(x_2)f(x_3))
 \Big\}
\nonumber\\
&=&\int_{{\Lambda}^4}d\x_4 \ \delta_{t}(x_{4}, x_1)\delta_{t}(x_{3}, x_2)
(f(x_1)-f(x_4))(f(x_2)-f(x_3))
\nonumber\\
&&\times\Big[ 
(\sG_{t}(x_{4}, x_1)\sG_{t}(x_{3}, x_2)-\sG_{t}(x_1, x_2)\sG_{t}(x_{3}, x_4))
(f(x_1)-f(x_4))(f(x_2)-f(x_3))
\nonumber\\
&&+\left( \frac{\partial}{\partial x_1} \sG_{t}(x_{1}, x_2)
\frac{\partial}{\partial x_2}\sG_{t}(x_{1}, x_2)
-\sG_{t}(x_{2}, x_1)
\frac{\partial^2}{\partial x_1 \partial x_2}\sG_{t}(x_{1}, x_2)\right)
\nonumber\\
&&\qquad\times (f(x_1)f(x_4)+f(x_2)f(x_3))(x_4-x_1)(x_3-x_2)
 \Big] +{\cal O}(t^2).
\nonumber
\end{eqnarray}
By the assumption (\ref{eqn:Bound1}) we obtain
(\ref{eqn:boundIj}) for $j=2$.

Since
\begin{eqnarray}
&&(-1)^{1+4} \delta_t(x_1,x_4) 
\det {\bf M}_{3,1}(x_1,x_2,x_3,x_4)^{\{1\} \{4\}}
= (-1)^{2+4} \delta_t(x_2,x_4) 
\det {\bf M}_{3,1}(x_2,x_1,x_3,x_4)^{\{2\} \{4\}}
\nonumber\\
&&\qquad 
=(-1)^{3+4} \delta_t(x_3,x_4) 
\det {\bf M}_{3,1}(x_3,x_2,x_1,x_4)^{\{3\} \{4\}},
\nonumber\\
&&(-1)^{1+4} \delta_t(x_1,x_4) 
\det {\bf M}_{2,2}(x_1,x_2,x_3,x_4)^{\{1\} \{4\}}
=(-1)^{1+3} \delta_t(x_1,x_3) 
\det {\bf M}_{2,2}(x_1,x_2,x_4,x_3)^{\{1\} \{3\}}
\nonumber\\
&&\qquad =(-1)^{2+3} \delta_t(x_2,x_3) 
\det {\bf M}_{2,2}(x_2,x_1,x_4,x_3)^{\{2\} \{3\}}
\nonumber\\
&& \qquad
=(-1)^{2+4} \delta_t(x_2,x_4) 
\det {\bf M}_{2,2}(x_2,x_1,x_3,x_4)^{\{2\} \{4\}},
\nonumber
\end{eqnarray}
we have
\begin{eqnarray}
\tilde{I_3} &=& 12 \int_{{\Lambda}^4}d\x_4 \ \prod_{i=1}^4 f(x_i) 
\delta_t(x_1,x_4)
\nonumber\\
&&\times\left[ -
\left| \matrix{
K(x_{2}, x_1) & \rho (x_2) & K(x_2, x_3) \cr
K(x_{3}, x_1) & K(x_3, x_2) & \rho (x_3) \cr
\sG_{t}(x_{4}, x_1) & \sG_{t}(x_{4}, x_2) &  \sG_{t}(x_{4}, x_3)
} \right|
+
\left| \matrix{
K(x_{2}, x_1) & \rho (x_2) & \sG_{-t}(x_2, x_3) \cr
\sG_t (x_{3}, x_1) & \sG_t (x_3, x_2) & \rho (x_3) \cr
\sG_{t}(x_{4}, x_1) & \sG_{t}(x_{4}, x_2) &  K(x_{4}, x_3)
} \right| \right].
\nonumber
\end{eqnarray}
Then we have (\ref{eqn:boundIj}) for $j=3$.
This completes the proof. \qed

\subsection{Bilinear Forms}

Since the Fredholm determinant of the dual operator
coincides with that of the original operator,
the reversibility (and also the stationarity)
of the present determinantal processes
is guaranteed.

Osada constructed $\mX$-valued reversible processes,
which have determinantal point fields as
their reversible stationary measures
by Dirichlet form approach.
With the Palm measure $\mu^{z}$
the Dirichlet form of Osada can be written as
\begin{equation}
{\cal E}(F,G)= \int_{{\Lambda}}dz\rho(z)
\int_{\mX} \mu^z(d\eta)
\frac{\partial}{\partial z}F(\{ z\}\cup\eta)
\frac{\partial}{\partial z}G(\{z\} \cup\eta)
\label{eqn:Dirichlet}
\end{equation}
for local smooth functions $F, G$
\cite{Osa04}(see also \cite{Spo87}).
By the general theory of Dirichlet forms
\cite{FOT94}, his processes are diffusion processes
({\it i.e.}, continuous strong Markov processes).
For our class of determinantal processes,
we found the following fact 
(Proposition \ref{thm:bilinear}).
It suggests that our determinantal 
processes are identified with Osada's processes.

A function $F$ on the configuration space $\mX$ is said to be
polynomial, if it is written of the form
$F(\xi)=\widetilde{F}(\langle f_1, \xi \rangle,
\langle f_2, \xi \rangle, \cdots,
\langle f_k, \xi \rangle)$
with a polynomial function $\widetilde{F}$
on ${\R}^{k}, k \in \N_0$,
where $f_j \in C_{0}^{\infty}(\R),
1 \leq j \leq k$.
Let $\wp$ be the set of all polynomial functions on
$\mX$, which is a dense subset of 
$L^{2}(\mX, \mu)$ ;
the space of square integrable
functions on $\mX$ with the determinantal point field $\mu$.

\begin{prop}
\label{thm:bilinear}
Let $\Xi(t)$ be the determinantal process,
whose multitime correlation functions are
given by (\ref{eqn:rhogen1}). 
Then for $F, G \in \wp$ we have
\begin{equation}
-\frac{\partial}{\partial t}
\E [F(\Xi(0)) G(\Xi(t))]
\Big|_{t=0}
= {\cal E}(F,G),
\label{eqn:bilinear}
\end{equation}
where ${\cal E}$ is given by (\ref{eqn:Dirichlet}).
\end{prop}
\vskip 0.5cm
Markov property of determinantal processes
has been studied by Borodin and Olshanski \cite{BO06,BO07b}.
But we can not prove that the
infinite determinantal processes in our class are Markovian.
The equivalence between our processes and
Osada's is not yet established.

To show this proposition, it is enough to consider the case
that $F$ and $G$ are of the form
$$
F(\xi)=\exp \left(
\sum_{m=1}^{k} \alpha_m \sum_{x \in \xi} f_m(x) \right),
\quad
G(\xi)=\exp \left(
\sum_{m=1}^{k} \beta_m \sum_{x \in \xi} g_m(x) \right)
$$
with $k \in \N_0$,
$\alpha_m, \beta_m \in \R$,
$f_m, g_m \in C_0^{\infty}(\R),
1 \leq m \leq k$.
Then if we set
$\chi_{1}(x)= \exp \left(
\sum_{m=1}^{k} \alpha_{m} f_m(x) \right)-1$
and
$\chi_{2}(x)= \exp \left(
\sum_{m=1}^{k} \beta_{m} g_m(x) \right)-1$,
$$
F(\xi)=\prod_{x\in\xi}(1+\chi_{1} (x)),
\quad
G(\xi)=\prod_{x\in\xi}(1+\chi_{2} (x)).
$$
The left-hand-side of (\ref{eqn:bilinear})
equals 
\begin{eqnarray}
&&
-\frac{\partial}{\partial t}\Phi_{t}
(\{f_m \}, \{g_m \}; \{\alpha_m \}, \{\beta_m \})
\Big|_{t=0}
\nonumber\\
&=&-\sum_{m=0}^{\infty} \sum_{n=0}^{\infty} 
\frac{1}{m! n!} 
\int_{{\Lambda}^{m}} d \x \int_{{\Lambda}^{n}} d {\bf y} \
\prod_{j=1}^{m} \chi_{1}(x_{j})
\prod_{k=1}^{n} \chi_{2}(y_{k})
\frac{\partial}{\partial t}
\cR(0,\{\x_{m}\};t, \{\y_{n}\})\bigg|_{t=0}. 
\nonumber
\end{eqnarray}
Since the Palm measure $\mu^z$ is a determinantal point 
field associated with
$\cR^z(\{\x_m \})$ (see (\ref{eqn:palm}) )
the right-hand-side of (\ref{eqn:bilinear}) equals
\begin{eqnarray}
&&\int_{{\Lambda}}dz\rho(z)
\frac{\partial \chi_{1} (z)}{\partial z}
\frac{\partial \chi_{2}(z)}{\partial z}
\int_{\mX} \mu^z(d\eta)F(\eta)G(\eta)
\nonumber\\
&&=\int_{{\Lambda}}dz\rho(z)
\frac{\partial \chi_{1} (z)}{\partial z}
\frac{\partial \chi_{2}(z)}{\partial z}
\sum_{m=0}^{\infty} \frac{1}{m!} 
\int_{{\Lambda}^m} d \x_m  \
\prod_{j=1}^{m} 
\Big\{ \chi_{1}(x_{j})\chi_{2}(x_{j})+\chi_{1}(x_{j})
+\chi_{2}(x_{j}) \Big\}
\cR^z(\{\x_{m}\}).
\nonumber
\end{eqnarray}
Hence Proposition \ref{thm:bilinear} 
can be derived from the following lemma.

\begin{lem}
\label{lem:key1}
For any $m, n \in \N$ we have
\begin{eqnarray}
&& \lim_{t\to 0}\int_{{\Lambda}^m} d \x_m \int_{{\Lambda}^n} d \y_n \
\prod_{j=1}^{m} \chi_{1}(x_{j})
\prod_{k=1}^{n} \chi_{2}(y_{k})
\frac{\partial}{\partial t}\cR(0,\{\x_m\};t,\{\y_n\})
\nonumber\\
&&= - \sum_{\ell=1}^{m \wedge n}
\frac{m ! n!}{(m-\ell)! (n-\ell)! (\ell-1)!}
\int_{{\Lambda}}dz\rho(z)
\frac{\partial \chi_{1} (z)}{\partial z}
\frac{\partial \chi_{2}(z)}{\partial z}
\nonumber\\
&&\qquad \times 
\int_{{\Lambda}^{m-\ell}} d \x_{m-\ell} 
\int_{{\Lambda}^{n-\ell}} d \y_{n-\ell} 
\int_{{\Lambda}^{\ell-1}} d \w_{\ell-1} \
\prod_{j=1}^{m-\ell} \chi_{1}(x_{j})
\prod_{k=1}^{n-\ell} \chi_{2}(y_{k})
\prod_{i=1}^{\ell-1} \chi_{1}(w_i)\chi_{2}(w_i)
\nonumber\\
&& \hskip 5cm \times
\cR^z(\{\x_{m-\ell}\} \cup \{\y_{n-\ell}\} \cup \{\w_{\ell-1}\}).
\nonumber
\end{eqnarray}
\end{lem}
\vskip 5mm

\noindent{\it Proof of Lemma \ref{lem:key1}.} \quad
We use the expansion formula (\ref{eqn:expansion1})
of the two-time correlation function
to calculate the integral
$$
I=\int_{{\Lambda}^m}d\x_m \int_{{\Lambda}^n}d\y_n
\prod_{i=1}^m \chi_{1} (x_i) \prod_{j=1}^n \chi_{2} (y_j)
\cR(0, \{\x_{m}\};t,\{\y_{n}\}).
$$
By permutation invariance of the integrand 
we put $a_i =m-\ell+i,b_i=i$
with $i=1,2,\dots, \ell$, and then
\begin{eqnarray}
I&=&\int_{{\Lambda}^m}d\x_m \int_{{\Lambda}^n}d\y_n
\prod_{i=1}^m \chi_{1} (x_i) \prod_{j=1}^n \chi_{2} (y_j)
\det {\bf M}(t, \y_n| \x_{m})
\nonumber\\
&+&\sum_{\ell=1}^{m\wedge n} 
{m \choose \ell}{n\choose \ell}
\int_{{\Lambda}^m}d\x_m \int_{{\Lambda}^n}d\y_n
\prod_{i=1}^m \chi_{1} (x_i) \prod_{j=1}^n \chi_{2} (y_j)
\nonumber\\
&& \hskip 6cm \times
\det {\bf D}(t, \y_n|\x_m)_{\{\a_\ell\} \{\b_\ell\}}
\det {\bf M}(t, \y_n|\x_m)^{\{m+\b_\ell\}\{\a_\ell\}}
\nonumber\\
&=&\int_{{\Lambda}^m}d\x_m \int_{{\Lambda}^n}d\y_n
\prod_{i=1}^m \chi_{1} (x_i) \prod_{j=1}^n \chi_{2} (y_j)
\det {\bf M}(t, \y_{n}| \x_{m})
\nonumber\\
&+&\sum_{\ell=1}^{m\wedge n} 
\frac{m! n!}{\ell !(m-\ell)!(n-\ell)!}
\int_{{\Lambda}^m}d\x_m \int_{{\Lambda}^n}d\y_n
\prod_{i=1}^m \chi_{1} (x_i) \prod_{j=1}^n \chi_{2} (y_j)
\prod_{k=1}^\ell \delta_t(x_{m-\ell+k},y_k)
\det {\bf M}^{\ell}_t(\x_{m};\y_{n}),
\nonumber
\end{eqnarray}
where
${\bf M}^{\ell}_t(\x_{m};\y_{n})
={\bf M}(t, \y_n|\x_m)^{\{m+1, \cdots, m+\ell\} 
\{m-\ell+1, \cdots, m\}}$.
By using (\ref{eqn:Gd-I}), (\ref{eqn:G-K}),
(\ref{eqn:Gt2}) and (\ref{eqn:dG-dG}),
we see that 
\begin{equation}
\lim_{t\to 0}\frac{\partial}{\partial t}
\int_{{\Lambda}^m}d\x_m \int_{{\Lambda}^n}d\y_n
\prod_{i=1}^m \chi_{1} (x_i) \prod_{j=1}^n \chi_{2} (y_j)
\det {\bf M}(t, \y_n| \x_{m})=0.
\label{eqn:l=0}
\end{equation}
From (\ref{eqn:GtA1}) and (\ref{eqn:Gt2})
\begin{small}
\begin{eqnarray}
&&\frac{\partial}{\partial t}
\int_{{\Lambda}^m}d\x_m \int_{{\Lambda}^n}d\y_n
\prod_{i=1}^m \chi_{1} (x_i) \prod_{j=1}^n \chi_{2} (y_j)
\prod_{k=1}^\ell \delta_t(x_{m-\ell+k},y_k)
\det {\bf M}^{\ell}_t(\x_{m};\y_{n})
\nonumber\\
&&=\sum_{p=1}^\ell 
\int_{{\Lambda}^m}d\x_m \int_{{\Lambda}^n}d\y_n
\prod_{i=1}^m \chi_{1} (x_i) \prod_{j=1}^n \chi_{2} (y_j)
\frac{\partial^2}{\partial y_p^2} \delta_t(x_{m-\ell+p},y_p)
\prod_{k=1,k\not= p}^\ell \delta_t(x_{m-\ell+k},y_k)
\det {\bf M}^{\ell}_t(\x_{m};\y_{n})
\nonumber\\
&&+\sum_{p=1}^{m-\ell} 
\int_{{\Lambda}^m}d\x_m \int_{{\Lambda}^n}d\y_n
\prod_{i=1}^m \chi_{1} (x_i) \prod_{j=1}^n \chi_{2} (y_j)
\prod_{k=1}^\ell \delta_t(x_{m-\ell+k},y_k)
\nonumber\\
&&\times\det 
\left( \matrix{
K(x_{1}, x_{1}) & \cdots &0& \cdots & K(x_{1}, x_{m-\ell}) 
&{\sG}_{t}(x_{1}, y_{1}) & \cdots&{\sG}_{t}(x_{1}, y_{n}) \cr
K(x_{2}, x_{1}) & \cdots &0& \cdots & K(x_{2}, x_{m-\ell}) 
&{\sG}_{t}(x_{2},y_{1}) & \cdots & {\sG}_{t}(x_{2}, y_{n}) \cr
   &   \cdots &  & & \cdots &  \cr
K(x_{m}, x_{1}) & \cdots &0& \cdots & K(x_{m}, x_{m-\ell}) 
& {\sG}_{t}(x_{m}, y_{1}) 
& \cdots& {\sG}_{t}(x_{m}, y_{n}) \cr
\sG_{-t}(y_{\ell+1}, x_{1}) & \cdots 
& \partial_{y_{\ell+1}}^2 {\sG}_{-t}(y_{\ell+1}, x_p)
& \cdots & \sG_{-t}(y_{\ell+1}, x_{m-\ell}) 
& K(y_{\ell+1},y_{1}) & \cdots & K(y_{\ell+1}, y_{n}) \cr
\sG_{-t}(y_{\ell+2}, x_{1}) & \cdots &
\partial_{y_{\ell+2}}^2 {\sG}_{-t}(y_{\ell+2}, x_p)
& \cdots & \sG_{-t}(y_{\ell+2}, x_{m-\ell}) 
& K(y_{\ell+2}, y_{1}) & \cdots & K(y_{\ell+2}, y_{n}) \cr
   &   \cdots &  & & \cdots & \cr
\sG_{-t}(y_{n}, x_{1}) & \cdots &
\partial_{y_{n}}^2 {\sG}_{-t}(y_{n}, x_p)
& \cdots
& \sG_{-t}(y_{n}, x_{m-\ell}) 
& K(y_{n}, y_{1}) & \cdots & K (y_{n},y_{n}) 
} \right)
\nonumber\\
&&-\sum_{p=1}^{n} 
\int_{{\Lambda}^m}d\x_m \int_{{\Lambda}^n}d\y_n
\prod_{i=1}^m \chi_{1} (x_i) \prod_{j=1}^n \chi_{2} (y_j)
\prod_{k=1}^\ell \delta_t(x_{m-\ell+k},y_k)
\nonumber\\
&&\times\det
\left( \matrix{
K(x_{1}, x_{1}) & \cdots & K(x_{1}, x_{m-\ell}) 
&{\sG}_{t}(x_{1}, y_{1}) & \cdots&
\partial_{y_p}^{2} {\sG}_{t}(x_{1}, y_{p}) & \cdots &
{\sG}_{t}(x_{1}, y_{n}) \cr
K(x_{2}, x_{1}) & \cdots & K(x_{2}, x_{m-\ell}) 
&{\sG}_{t}(x_{2},y_{1}) & \cdots &
\partial_{y_p}^{2} {\sG}_{t}(x_{2}, y_{p}) & \cdots &
 {\sG}_{t}(x_{2}, y_{n}) \cr
   &   \cdots &  & & \cdots &  \cr
K(x_{m}, x_{1}) & \cdots & K(x_{m}, x_{m-\ell}) 
& {\sG}_{t}(x_{m}, y_{1}) 
& \cdots& 
\partial_{y_p}^{2} {\sG}_{t}(x_{m}, y_{p}) & \cdots &
{\sG}_{t}(x_{m}, y_{n}) \cr
\sG_{-t}(y_{\ell+1}, x_{1}) & \cdots & \sG_{-t}(y_{\ell+1}, x_{m-\ell}) 
& K(y_{\ell+1},y_{1}) & \cdots & 
0& \cdots &
K(y_{\ell+1}, y_{n}) \cr
\sG_{-t}(y_{\ell+2}, x_{1}) & \cdots &\sG_{-t}(y_{\ell+2}, x_{m-\ell}) 
& K(y_{\ell+2}, y_{1}) & \cdots &
0& \cdots &
 K(y_{\ell+2}, y_{n}) \cr
   &   \cdots &  & & \cdots & \cr
\sG_{-t}(y_{n}, x_{1}) & \cdots & \sG_{-t}(y_{n}, x_{m-\ell}) 
& K(y_{n}, y_{1}) & \cdots &
0& \cdots & K (y_{n},y_{n}) 
} \right)
\nonumber\\
&&\equiv I_1(t) + I_2(t) -I_3(t),
\label{I(1)+I(2)-I(3)}
\end{eqnarray}
\end{small}
where we have used the abbreviation
$\partial^2_y = \partial^2/\partial y^2$.
By partial integration
\begin{small}
\begin{eqnarray}
I_1(t) &=& \sum_{p=1}^\ell 
\int_{{\Lambda}^m}d\x_m \int_{{\Lambda}^n}d\y_n
\prod_{i=1}^m \chi_{1} (x_i) \prod_{j=1, j\not=p}^n \chi_{2} (y_j)
\prod_{k=1}^\ell \delta_t(x_{m-\ell+k},y_k)
\frac{\partial^2}{\partial y_p^2}
\bigg\{ \chi_{2} (y_p) \det {\bf M}^{\ell}_t(\x_{m};\y_{n})\bigg\}
\nonumber\\
&=&\sum_{p=1}^\ell 
\int_{{\Lambda}^m}d\x_m \int_{{\Lambda}^n}d\y_n 
\prod_{i=1}^m \chi_{1} (x_i) \prod_{j=1}^n \chi_{2} (y_j)
\prod_{k=1}^\ell \delta_t(x_{m-\ell+k},y_k)
\frac{\partial^2}{\partial y_p^2}
\det {\bf M}^{\ell}_t(\x_{m};\y_{n})
\nonumber\\
&&+ \sum_{p=1}^\ell 
\int_{{\Lambda}^m}d\x_m \int_{{\Lambda}^n}d\y_n
\prod_{i=1}^m \chi_{1} (x_i) \prod_{j=1, j\not=p}^n \chi_{2} (y_j)
\prod_{k=1}^\ell \delta_t(x_{m-\ell+k},y_k)
\nonumber\\
&&\qquad\times
\bigg\{\frac{\partial^2 \chi_{2} (y_p)}{\partial y_p^2}
\det {\bf M}^{\ell}_t(\x_{m};\y_{n})
+2 \frac{\partial \chi_{2} (y_p)}{\partial y_p} 
\frac{\partial \det {\bf M}^{\ell}_t(\x_{m};\y_{n})}{\partial y_p}\bigg\}
\nonumber\\
&\equiv& I_{11}(t)+I_{12}(t).
\nonumber
\end{eqnarray}
\end{small}
By simple calculation 
\begin{small}
\begin{eqnarray}
&&I_{11}(t) + I_2(t) -I_3(t)
\nonumber\\
&&=\sum_{p=1}^\ell 
\int_{{\Lambda}^m}d\x_m \int_{{\Lambda}^n}d\y_n
\prod_{i=1}^m \chi_{1} (x_i) \prod_{j=1}^n \chi_{2} (y_j)
\prod_{k=1}^\ell \delta_t(x_{m-\ell+k},y_k)
\nonumber\\
&&\times \det
\left( \matrix{
K(x_{1}, x_{1}) & \cdots & K(x_{1}, x_{m-\ell}) 
&{\sG}_{t}(x_{1}, y_{1}) & \cdots&
0 & \cdots &
{\sG}_{t}(x_{1}, y_{n}) \cr
K(x_{2}, x_{1}) & \cdots & K(x_{2}, x_{m-\ell}) 
&{\sG}_{t}(x_{2},y_{1}) & \cdots &
0 & \cdots &
 {\sG}_{t}(x_{2}, y_{n}) \cr
   &   \cdots &  & & \cdots &  \cr
K(x_{m}, x_{1}) & \cdots & K(x_{m}, x_{m-\ell}) 
& {\sG}_{t}(x_{m}, y_{1}) & \cdots&
0 & \cdots &
 {\sG}_{t}(x_{m}, y_{n}) \cr
\sG_{-t}(y_{\ell+1}, x_{1}) & \cdots & \sG_{-t}(y_{\ell+1}, x_{m-\ell}) 
& K(y_{\ell+1},y_{1}) & \cdots &
\partial_{y_{p}}^{2} K(y_{\ell+1}, y_p) & \cdots &
 K(y_{\ell+1}, y_{n}) \cr
\sG_{-t}(y_{\ell+2}, x_{1}) & \cdots &\sG_{-t}(y_{\ell+2}, x_{m-\ell}) 
& K(y_{\ell+2}, y_{1}) & \cdots &
\partial_{y_{p}}^{2} K(y_{\ell+2}, y_p) & \cdots &
 K(y_{\ell+2}, y_{n}) \cr
   &   \cdots &  & & \cdots & \cr
\sG_{-t}(y_{n}, x_{1}) & \cdots & \sG_{-t}(y_{n}, x_{m-\ell}) 
& K(y_{n}, y_{1}) & \cdots &
\partial_{y_{p}}^{2} K(y_{n}, y_p) & \cdots &
 K (y_{n},y_{n}) 
} \right)
\nonumber\\
&&+\sum_{p=1}^{m-\ell} 
\int_{{\Lambda}^m}d\x_m \int_{{\Lambda}^n}d\y_n
\prod_{i=1}^m \chi_{1} (x_i) \prod_{j=1}^n \chi_{2} (y_j)
\prod_{k=1}^\ell \delta_t(x_{m-\ell+k},y_k)
\nonumber\\
&&\times\det 
\left( \matrix{
K(x_{1}, x_{1}) & \cdots &0& \cdots & K(x_{1}, x_{m-\ell}) 
&{\sG}_{t}(x_{1}, y_{1}) & \cdots&{\sG}_{t}(x_{1}, y_{n}) \cr
K(x_{2}, x_{1}) & \cdots &0& \cdots & K(x_{2}, x_{m-\ell}) 
&{\sG}_{t}(x_{2},y_{1}) & \cdots & {\sG}_{t}(x_{2}, y_{n}) \cr
   &   \cdots &  & & \cdots &  \cr
K(x_{m}, x_{1}) & \cdots &0& \cdots & K(x_{m}, x_{m-\ell}) 
& {\sG}_{t}(x_{m}, y_{1}) 
& \cdots& {\sG}_{t}(x_{m}, y_{n}) \cr
\sG_{-t}(y_{\ell+1}, x_{1}) & \cdots 
& \partial_{y_{\ell+1}}^2 {\sG}_{-t}(y_{\ell+1}, x_p)
& \cdots & \sG_{-t}(y_{\ell+1}, x_{m-\ell}) 
& K(y_{\ell+1},y_{1}) & \cdots & K(y_{\ell+1}, y_{n}) \cr
\sG_{-t}(y_{\ell+2}, x_{1}) & \cdots &
\partial_{y_{\ell+2}}^2 {\sG}_{-t}(y_{\ell+2}, x_p)
& \cdots & \sG_{-t}(y_{\ell+2}, x_{m-\ell}) 
& K(y_{\ell+2}, y_{1}) & \cdots & K(y_{\ell+2}, y_{n}) \cr
   &   \cdots &  & & \cdots & \cr
\sG_{-t}(y_{n}, x_{1}) & \cdots &
\partial_{y_{n}}^2 {\sG}_{-t}(y_{n}, x_p)
& \cdots
& \sG_{-t}(y_{n}, x_{m-\ell}) 
& K(y_{n}, y_{1}) & \cdots & K (y_{n},y_{n}) 
} \right)
\nonumber\\
&&-\sum_{p=\ell +1}^{n} 
\int_{{\Lambda}^m}d\x_m \int_{{\Lambda}^n}d\y_n
\prod_{i=1}^m \chi_{1} (x_i) \prod_{j=1}^n \chi_{2} (y_j)
\prod_{k=1}^\ell \delta_t(x_{m-\ell+k},y_k)
\nonumber\\
&&\times\det
\left( \matrix{
K(x_{1}, x_{1}) & \cdots & K(x_{1}, x_{m-\ell}) 
&{\sG}_{t}(x_{1}, y_{1}) & \cdots&
\partial_{y_p}^{2} {\sG}_{t}(x_{1}, y_{p}) & \cdots &
{\sG}_{t}(x_{1}, y_{n}) \cr
K(x_{2}, x_{1}) & \cdots & K(x_{2}, x_{m-\ell}) 
&{\sG}_{t}(x_{2},y_{1}) & \cdots &
\partial_{y_p}^{2} {\sG}_{t}(x_{2}, y_{p}) & \cdots &
 {\sG}_{t}(x_{2}, y_{n}) \cr
   &   \cdots &  & & \cdots &  \cr
K(x_{m}, x_{1}) & \cdots & K(x_{m}, x_{m-\ell}) 
& {\sG}_{t}(x_{m}, y_{1}) 
& \cdots& 
\partial_{y_p}^{2} {\sG}_{t}(x_{m}, y_{p}) & \cdots &
{\sG}_{t}(x_{m}, y_{n}) \cr
\sG_{-t}(y_{\ell+1}, x_{1}) & \cdots & \sG_{-t}(y_{\ell+1}, x_{m-\ell}) 
& K(y_{\ell+1},y_{1}) & \cdots & 
0& \cdots &
K(y_{\ell+1}, y_{n}) \cr
\sG_{-t}(y_{\ell+2}, x_{1}) & \cdots &\sG_{-t}(y_{\ell+2}, x_{m-\ell}) 
& K(y_{\ell+2}, y_{1}) & \cdots &
0& \cdots &
 K(y_{\ell+2}, y_{n}) \cr
   &   \cdots &  & & \cdots & \cr
\sG_{-t}(y_{n}, x_{1}) & \cdots & \sG_{-t}(y_{n}, x_{m-\ell}) 
& K(y_{n}, y_{1}) & \cdots &
0& \cdots & K (y_{n},y_{n}) 
} \right).
\nonumber
\end{eqnarray}
\end{small}
By using (\ref{eqn:Gd-I}) and (\ref{eqn:dG-dG}) we have
\begin{equation}
\lim_{t\to 0}\bigg\{I_{11}(t) + I_2(t) -I_{3}(t) \bigg\}=0.
\label{Is->0}
\end{equation}
Suppose that $x_{m-\ell+i} \to y_i$, $i=1,2,\dots \ell$,
\begin{eqnarray}
&&\det {\bf M}^{\ell}_t(\x_{m};\y_{n})\to 
\cR (\{\x_{m-\ell}\} \cup \{\y_n\}),
\nonumber
\\
&&\frac{\partial \det {\bf M}^{\ell}_t(\x_{m};\y_{n})}{\partial y_p}
\to \frac{1}{2}
\frac{\partial \cR (\{\x_{m-\ell}\} \cup \{\y_{n}\})}{\partial y_p}
\quad \mbox{in} \quad t \to 0
\nonumber
\end{eqnarray}
for any $p=1,2,\dots, \ell$.
Hence we have
\begin{eqnarray}
\lim_{t\to 0} I_{12}(t)
&=&\sum_{p=1}^\ell 
\int_{{\Lambda}^{m-\ell}}d\x_{m-\ell} \int_{{\Lambda}^n}d\y_n
\prod_{i=1}^{m-\ell} \chi_{1} (x_i) 
\prod_{k=1, k\not=p}^{\ell}\chi_{1}(y_k)\chi_{2} (y_k)
\prod_{j=\ell+1}^n \chi_{2} (y_j)
\nonumber\\
&&\qquad\times \chi_{1}(y_p)
\bigg\{ \frac{\partial^2 \chi_{2} (y_p)}{\partial y_p^2}
\cR (\{\x_{m-\ell}\} \cup \{\y_n\})
+\frac{\partial \chi_{2} (y_p)}{\partial y_p} 
\frac{\partial \cR (\{\x_{m-\ell}\} \cup \{\y_{n}\})}
{\partial y_p}
\bigg\}
\nonumber\\
&=& \sum_{p=1}^\ell 
\int_{{\Lambda}^{m-\ell}}d\x_{m-\ell} \int_{{\Lambda}^n}d\y_n
\prod_{i=1}^{m-\ell} \chi_{1} (x_i) 
\prod_{k=1, k\not=p}^{\ell}\chi_{1}(y_k)\chi_{2} (y_k)
\prod_{j=\ell+1}^n \chi_{2} (y_j)
\nonumber\\
&&\qquad\times \chi_{1}(y_p)
\frac{\partial}{\partial y_p} 
\bigg\{ \frac{\partial \chi_{2} (y_p) }{\partial y_p}
\cR (\{\x_{m-\ell}\} \cup \{\y_n\})
\bigg\}
\nonumber\\
&=&-\sum_{p=1}^\ell 
\int_{{\Lambda}^{m-\ell}}d\x_{m-\ell} \int_{{\Lambda}^n}d\y_n
\prod_{i=1}^{m-\ell} \chi_{1} (x_i) 
\prod_{k=1, k\not=p}^{\ell}\chi_{1}(y_k)\chi_{2} (y_k)
\prod_{j=\ell+1}^n \chi_{2} (y_j)
\nonumber\\
&&\qquad\times
\frac{\partial \chi_{1}(y_p)}{\partial y_p}
\frac{\partial \chi_{2} (y_p) }{\partial y_p}
\cR (\{\x_{m-\ell}\} \cup \{\y_n\})
\nonumber\\
&=&-\ell 
\int_{{\Lambda}^{m-\ell}}d\x_{m-\ell}\int_{{\Lambda}}dz
\int_{{\Lambda}^{\ell-1}}d\w_{\ell-1} \int_{{\Lambda}^{n-\ell}}d\y_{n-\ell}
\prod_{i=1}^{m-\ell} \chi_{1} (x_i) 
\prod_{k=1}^{\ell-1}\chi_{1}(y_k)\chi_{2} (y_k)
\prod_{j=1}^n \chi_{2} (w_j)
\nonumber\\
&&\qquad\times
\frac{\partial \chi_{1}(z)}{\partial z}
\frac{\partial \chi_{2} (z) }{\partial z}
\cR (\{\x_{m-\ell}\} \cup \{z\} \cup \{\w_{\ell-1}\}
\cup \{\y_{n-\ell}\})
\nonumber\\
&=&-\ell \int_{{\Lambda}}dz \cR(z)
\frac{\partial \chi_{1}(z)}{\partial z}
\frac{\partial \chi_{2} (z) }{\partial z}
\int_{{\Lambda}^{m-\ell}}d\x_{m-\ell}
\int_{{\Lambda}^{\ell-1}}d\w_{\ell-1} \int_{{\Lambda}^{n-\ell}}d\y_{n-\ell}
\nonumber\\
&&\qquad\times
\prod_{i=1}^{m-\ell} \chi_{1} (x_i) 
\prod_{k=1}^{\ell-1}\chi_{1}(y_k)\chi_{2} (y_k)
\prod_{j=1}^n \chi_{2} (w_j)
\cR^z (\{\x_{m-\ell}\} \cup \{\w_{\ell-1}\}
\cup \{\y_{n-\ell}\}).
\label{limI(12)}
\end{eqnarray}
Combining (\ref{eqn:l=0}), (\ref{I(1)+I(2)-I(3)}),
(\ref{Is->0}) and (\ref{limI(12)}),
we obtain Lemma \ref{lem:key1}.
\qed

\SSC{CONCLUDING REMARKS}

The study on noncolliding BM 
and determinantal processes reported in this paper
will be extended in several directions.
We would like to give some of future problems below.
\begin{description}
\item{(1)} \quad
In Section 4  we let the initial configuration
$\nu_0$ be the GUE-eigenvalue distribution,
and shown that the system is determinantal.
If we let $\nu_0$ be the eigenvalue distribution
in the Gaussian orthogonal ensemble (GOE),
\begin{equation}
\nu_0^{\rm GOE}(\x^{(0)})
=\frac{1}{C'_N} t_0^{-N(N+1)/4} 
e^{-|\x^{(0)}|^2/2t_0} h_N(\x^{(0)}),
\quad \x^{(0)} \in \W_N,
\label{eqn:GOE1}
\end{equation}
(\ref{eqn:nut1}) becomes
\begin{eqnarray}
\nu_t(\x) &\propto& h_{N}(\x) \int_{\W_N}
f_N(t-t_0, \x| \x^{(0)}) d \x^{(0)}
\nonumber\\
&=& \frac{1}{N!} h_{N}(\x) 
\int_{\R^{N}} f_{N}(t-t_0, \x|\x^{(0)}) \,
{\rm sgn}(h_N(\x^{(0)})) d \x^{(0)}.
\label{eqn:GOE2}
\end{eqnarray}
Instead of the Heine identity (\ref{eqn:Heine}),
we should use 
the de Bruin identity \cite{deB55}
$$
\int_{\W_N} d \x \,
\det_{1 \leq j, k \leq N} \Big[ g_{j}(x_k) \Big]
={\rm Pf}_{1 \leq j, k \leq N}
\left[ \int_{\R} dx \int_{\R} d \widetilde{x} \,
{\rm sgn}(\widetilde{x}-x) g_{j}(x) g_k(\widetilde{x}) 
\right]
$$
for integrable continuous functions
$g_j, 1 \leq j \leq N$.
As shown in Appendix A of \cite{KT07}, for example,
the generating function of multitime correlation
functions is then expressed by the Fredholm Pfaffian
\cite{Rai00} and the system becomes a Pfaffian process,
in the sense that any multitime correlation function
is given by a Pfaffian.
Such Pfaffian processes have been studied by
many authors 
\cite{BS03,Sos03b,Sos04,BR04,SI04,Ferrari04,IS05}.
The systems studied in \cite{FNH99,NKT03,N03,KNT04} are
also Pfaffian processes, since the `quaternion
determinantal expressions' of correlation functions, 
introduced and developed by Dyson, Mehta, Forrester, and Nagao 
\cite{D70,Mehta89,Meh04,NF99,Nag01},
are readily transformed to Pfaffian expressions.
As implied by Proposition \ref{thm:GUE}, the system
exhibits a transition from GOE distribution to GUE
distribution \cite{MP83,PM83,KT02}.
Continuity of sample paths and general characterization
of infinite Pfaffian processes will be interesting
problems.
The case with other initial distribution
(in particular, when it has continuous parameters)
will be interesting \cite{KT04,KT07}.

\item{(2)} \quad
As explained in Sections 1 and 3, the present
noncolliding BM is the $h$-transform of the
absorbing BM in the Weyl chamber 
(\ref{eqn:Weyl1}) of type $A_{N-1}$.
We can find appropriate $h$-transforms
of the absorbing BMs in the Weyl chambers
of types $C_N$ and $D_{N}$.
The obtained noncolliding diffusion processes are
stochastic versions of non-standard random matrix
ensembles, which were called the class C and class D, respectively,
by Altland and Zirnbauer \cite{AZ96,AZ97,Zir96}.
The stochastic version of the chiral GUE,
realized by the noncolliding squared Bessel process,
was studied by K\"onig and O'Connell \cite{KO01},
which is also obtained as an $h$-transform
of the absorbing BM in the Weyl chamber of type $C_N$.
See \cite{Gra99,KT04} for more details.
Systematic classifications of determinantal
and Pfaffian processes will be important.

\item{(3)} \quad
There are many other examples of finite and infinite
determinantal processes, which are not considered
in the present paper.
Markov processes on partitions (Young diagrams)
have been studied and 
determinantal processes associated with other types
of projections than ours have been reported
\cite{OR03,BO06,BO07a}.
The determinantal processes, whose kernels
are expressed using multiple orthogonal functions 
({\it e.g.} the Pearcey kernel),
are discussed in 
\cite{BH98a,BH98b,BK04a,BK04b,ABK05,IS05,
TW06,BK07,OR07}.
The consideration given in Sections 6 and 7
in the present paper should be generalized.
\end{description}

\vskip 1cm
\noindent{\bf ACKNOWLEDGMENTS}

The present authors would like to thank
H. Osada, Y. Takahashi, A. Borodin and G. Olshanski
for useful discussions on determinantal processes.
M.K. is supported in part by
the Grant-in-Aid for Scientific Research 
(KIBAN-C, No.17540363) of Japan Society for
the Promotion of Science.
H.T. is supported in part by
the Grant-in-Aid for Scientific Research 
(KIBAN-C, No.19540114) of Japan Society for
the Promotion of Science.

\begin{small}

\end{small}

\begin{thebibliography}{99}
\bibitem{AM99}
M. Adler and P. van Moerbeke,
The spectrum of coupled random matrices,
{\it Ann. Math.} {\bf 149}:921-976 (1999).

\bibitem{AM05}
M. Adler and P. van Moerbeke,
PDF's for the joint distributions of the
Dyson, Airy and Sine processes,
{\it Ann. Probab.} {\bf 33}:1326-1361 (2005).

\bibitem{AZ96}
A. Altland and M. R. Zirnbauer, 
Random matrix theory of a chaotic Andreev quantum dot,
{\it Phys. Rev. Lett.} {\bf 76}:3420-3424 (1996).

\bibitem{AZ97}
A. Altland and M. R. Zirnbauer, 
Nonstandard symmetry classes
in mesoscopic normal-superconducting hybrid structures,
{\it Phys. Rev.} B {\bf 55}:1142-1161 (1997).

\bibitem{ABK05}
A. I. Aptekarev, P. M. Bleher and A. B. J. Kuijlaars,
Large $n$ limit of Gaussian random matrices with
external source, part II,
{\it Commun. Math. Phys.} {\bf 259}:367-389 (2005).

\bibitem{AHG06}
G. Arcioni, S. de Haro and P. Gao,
Diffusion model for SU(N) QCD screening,
{\it Phys. Rev.} D {\bf 73}:074508 (2006).

\bibitem{BBDS06}
J. Baik, A. Borodin, P. Deift and T. Suidan,
A model for the bus system in Cuernavaca (Mexico),
{\it J. Phys. A: Math. Gen.} {\bf 39}:8965-8975 (2006).

\bibitem{Bal00}
A. B. Balantekin,
Character expansions, Itzykson-Zuber integrals, and
the QCD partition function,
{\it Phys. Rev.} D {\bf 62}:085017/1-8 (2000).

\bibitem{Bal02}
A. B. Balantekin,
Character expansions for the orthogonal and symplectic
groups,
{\it J. Math. Phys.} {\bf 43}:604-620 (2002).

\bibitem{Bat53}
H. Bateman,
{\it Higher Transcendental Functions}, 
(A. Erd\'elyi Ed.), Vol. 2, (McGraw Hill, 1953).

\bibitem{Baym69}
G. Baym,
{\it Lectures on Quantum Mechanics},
(W. A. Benjamin, Inc., 1969).

\bibitem{Bil99}
P. Billingsley,
{\it Convergence of Probability Measures},
2nd ed.,
(John Willey \& Sons, Inc., 1999).

\bibitem{Bog06}
N. M. Bogoliubov,
XX0 Heisenberg chain and random walks,
{\it J. Math. Sci.} {\bf 138}:5636-5643 (2006).

\bibitem{BO06}
A. Borodin and G. Olshanski,
Markov processes on partitions,
{\it Probab. Th. Rel. Fields}, {\bf 135}:84-152 (2006).

\bibitem{BO07a}
A. Borodin and G. Olshanski,
Asymptotics of Plancherel-type random partitions,
{\it Journal of Algebra}, {\bf 313}:40-60 (2007).

\bibitem{BO07b}
A. Borodin and G. Olshanski,
Infinite-dimensional diffusions as limits
of random walks on partitions,
{\sf arXiv:math.PR/07061034}.

\bibitem{BR04}
A. Borodin and E. M. Rains,
Eynard-Mehta theorem, Schur process, and their pfaffian analogs,
{\it J. Stat. Phys.} {\bf 121}:291-317 (2005).

\bibitem{BS02}
A. N. Borodin and P. Salminen,
{\it Handbook of Brownian Motion - Facts and Formulae}, 2nd ed.,
(Birkh\"auser, Basel, 2002).

\bibitem{BS03}
A. Borodin and A. Soshnikov,
Janossy densities. I. Determinantal ensembles,
{\it J. Stat. Phys.}{\bf 113}:595-610 (2003).

\bibitem{BH98a}
E. Br\'ezin and S. Hikami,
Universal singularity at the closure of a gap 
in a random matrix theory,
{\it Phys. Rev.} E {\bf 57}:4140-4149 (1998).

\bibitem{BH98b}
E. Br\'ezin and S. Hikami,
Level spacing of random matrices in an external sources,
{\it Phys. Rev.} E {\bf 58}:7176-7185 (1998).

\bibitem{BK04a}
P. M. Bleher and A. B. J. Kuijlaars,
Random matrices with external sources and multiple 
orthogonal polynomials,
{\it Int. Math. Research Notices} {\bf 2004}:109-129 (2004).

\bibitem{BK04b}
P. M. Bleher and A. B. J. Kuijlaars,
Large $n$ limit of Gaussian random matrices with 
external sources, part I,
{\it Commun. Math. Phys.} {\bf 252}:43-76 (2004).

\bibitem{BK07}
P. M. Bleher and A. B. J. Kuijlaars,
Large $n$ limit of Gaussian random matrices with 
external sources, part III: double scaling limit,
{\it Commun. Math. Phys.} {\bf 270}:481-517 (2007).

\bibitem{deB55}
N. G. de Bruijn, 
On some multiple integrals involving determinants,
{\it J. Indian Math. Soc.} {\bf 19}:133-151 (1955).

\bibitem{deG68}
P.-G. de Gennes,
Soluble model for fibrous structures with steric constraints,
{\it J. Chem. Phys.} {\bf 48}:2257-2259 (1968).

\bibitem{deH05}
S. de Haro,
Chern-Simons theory, 2d Yang-Mills, and Lie algebra
wanders,
{\it Nucl. Phys.} {\bf B 730}:312-351 (2005).

\bibitem{Dirac58}
P. A. M. Dirac,
{\it The Principles of Quantum Mechanics},
4th ed.,
(Oxford University Press, 1958).

\bibitem{Doo84}
J. L. Doob,
{\it Classical Potential Theory and its Probabilistic Counterpart},
(Springer, New York, 1984).

\bibitem{Dys62}
F. J. Dyson, 
A Brownian-motion model for the eigenvalues of a random matrix,
{\it J. Math. Phys.} {\bf 3}:1191-1198 (1962).

\bibitem{D70}
F. J. Dyson, 
Correlation between the eigenvalues of a random matrix,
{\it Commun. Math. Phys.} {\bf 19}:235-250 (1970).

\bibitem{EK06}
P. Eichelsbacher and W. K\"onig,
Ordered random walks,
{\sf arXiv:math.PR/0610850}.

\bibitem{EG95} 
J. W. Essam and A. J. Guttmann, 
Vicious walkers and directed polymer networks in
general dimensions,
{\it Phys. Rev.} E {\bf 52}:5849-5862 (1995).

\bibitem{Ferrari04}
P. L. Ferrari, 
Polynuclear growth on a flat substrate and
edge scaling of GOE eigenvalues, 
{\it Commun. Math. Phys.} {\bf 252}:77-109 (2004).

\bibitem{FS03}
P. L. Ferrari and H. Spohn,
Step fluctuations for a faceted crystal,
{\it J. Stat. Phys.} {\bf 113}:1-46 (2003).

\bibitem{Fis84}
M. E. Fisher, 
Walks, walls, wetting, and melting,
{\it J. Stat. Phys.} {\bf 34}:667-729 (1984).

\bibitem{FNH99}
P. J. Forrester, T. Nagao and G. Honner,
Correlations for the orthogonal-unitary and
symplectic-unitary transitions at the
hard and soft edges,
{\it Nucl. Phys.} {\bf B553[PM]}:601-643 (1999).

\bibitem{FOT94}
M. Fukushima, Y. Oshima and M. Takeda,
{\it Dirichlet Forms and Symmetric Markov Processes},
(Walter de Gruyter, Berlin, 1994).

\bibitem{Ful97}
W. Fulton,
{\it Young Tableaux with Applications
to Representation Theory and Geometry},
(Cambridge Univ. Press, Cambridge, 1997).

\bibitem{GV85}
I. Gessel and G. Viennot, 
Binomial determinants, paths, and hook length formulae,
{\it Adv. in Math.} {\bf 58}:300-321 (1985).

\bibitem{Gra99}
D. J. Grabiner, 
Brownian motion in a Weyl chamber,
non-colliding particles, and random matrices,
{\it Ann. Inst. Henri Poincar\'e,
Probab. Stat.} {\bf 35}:177-204 (1999).

\bibitem{GS92}
G. R. Grimmett and D. R. Stirzaker,
{\it Probability and Random Processes},
2nd ed.,
(Clarendon Press, Oxford, 1992).

\bibitem{HO07}
J. Harnad and A. Yu. Orlov,
Fermionic construction of tau functions and
random processes,
{\it Physica} D (2007, in press); 
{\sf arXiv:math-ph/07041157}.

\bibitem{IS05}
T. Imamura and T. Sasamoto, 
Polynuclear growth model with external source
and random matrix model with deterministic source,
{\it Phys. Rev.} E {\bf 71}:041606 (2005).

\bibitem{Joh02}
K. Johansson,
Non-intersecting paths, random tilings and random matrices,
{\it Probab. Th. Rel. Fields}, {\bf 123}:225-280 (2002).

\bibitem{Joh03}
K. Johansson, 
Discrete polynuclear growth and determinantal processes,
{\it Commun. Math. Phys.} {\bf 242}:277-329 (2003).

\bibitem{Kal97}
O. Kallenberg,
{\it Foundations of Modern Probability},
(Springer, 1997).

\bibitem{KS91}
I. Karatzas and S. E. Shreve,
{\it Brownian Motion and Stochastic Calculus},
2nd ed.,
(Springer, 1991).

\bibitem{KM59}
S. Karlin and J. McGregor,
Coincidence probabilities,
{\it Pacific J. Math.}:{\bf 9} 1141-1164 (1959).

\bibitem{K00}
M. Katori,
Percolation transitions and wetting transitions in
stochastic models,
{\it Brazilian J. Phys.} {\bf 30}:83-96 (2000).

\bibitem{KNT04}
M. Katori, T. Nagao and H. Tanemura,
Infinite systems of non-colliding Brownian particles,
{\it Adv. Stud. in Pure Math.} {\bf 39}
``{\it Stochastic Analysis on Large Scale Interacting Systems}", 
pp.283-306, 
(Mathematical Society of Japan, Tokyo, 2004);
{\sf arXiv:math.PR/0301143}.

\bibitem{KT02}
M. Katori and H. Tanemura,
Scaling limit of vicious walks and two-matrix model,
{\it Phys. Rev.} E {\bf 66}:011105/1-12 (2002).

\bibitem{KT03a}
M. Katori and H. Tanemura,
Functional central limit theorems for vicious walkers,
{\it Stoch. Stoch. Rep.} {\bf 75}:369-390 (2003);
{\sf arXiv:math.PR/0203286}.

\bibitem{KT03b}
M. Katori and H. Tanemura,
Noncolliding Brownian motions and Harish-Chandra
formula,
{\it Elect. Comm. in Probab.} {\bf 8}:112-121 (2003).

\bibitem{KT04}
M. Katori and H. Tanemura,
Symmetry of matrix-valued stochastic processes and
noncolliding diffusion particle systems,
{\it J. Math. Phys.} {\bf 45}:3058-3085 (2004).

\bibitem{KT05}
M. Katori and H. Tanemura,
Nonintersecting paths, Noncolliding diffusion processes
and representation theory,
{\it RIMS Kokyuroku} {\bf 1438}:83-102 (2005);
{\sf arXiv:math.PR/0501218}.

\bibitem{KT07}
M. Katori and H. Tanemura, 
Infinite systems of non-colliding generalized meanders
and Riemann-Liouville differintegrals,
{\it Probab. Th. Rel. Fields}, {\bf 138}:113-156 (2007).

\bibitem{KTNK03}
M. Katori, H. Tanemura, T. Nagao and N. Komatsuda,
Vicious walk with a wall, noncolliding meanders,
chiral and Bogoliubov-de Gennes random matrices,
{\it Phys. Rev.} E  {\bf 68}:021112/1-16 (2003).

\bibitem{KO01}
W. K\"onig and N. O'Connell,
Eigenvalues of the Laguerre process as non-colliding 
squared Bessel process,
{\it Elec. Comm. in Probab.} {\bf 6}:107-114 (2001).

\bibitem{Kra06}
C. Krattenthaler,
Watermelon configurations with wall interaction:
exact and asymptotic results,
{\it J. Phys. Conf. Series}
{\bf 42}:179-212 (2006).

\bibitem{KGV00}
C. Krattenthaler, A. J. Guttmann and X. G. Viennot,
Vicious walkers, friendly walkers and Young tableaux: II.
With a wall,
{\it J. Phys. A: Math. Phys.} {\bf 33}:8835-8866 (2000).

\bibitem{Lin73}
B. Lindstr\"om, 
On the vector representations of induced matroids,
{\it Bull. London Math. Soc.} {\bf 5}:85-90 (1973).

\bibitem{Mac82}
I. G. Macdonald,
Some conjectures for root systems,
{\it SIAM J. Math. Anal.} {\bf 13}:988-1007 (1982).

\bibitem{Mac95}
I. G. Macdonald,
{\it Symmetric Functions and Hall Polynomials}, 2nd ed., 
(Oxford Univ. Press, Oxford, 1995).

\bibitem{Mehta89}
M. L. Mehta, 
{\it Matrix Theory},
(Editions de Physique, Orsay, 1989).

\bibitem{Meh04}
M. L. Mehta, 
{\it Random Matrices}, 3rd ed.,
(Elsevier, Amsterdam, 2004).

\bibitem{MP83}
M. L. Mehta and A. Pandey,
On some Gaussian ensemble of Hermitian matrices,
{\it J. Phys. A: Math. Gen.} {\bf 16}:2655-2684 
(1983).

\bibitem{Nag01}
T. Nagao, 
Correlation functions for multi-matrix models
and quaternion determinants,
{\it Nucl. Phys.} {\bf B602}:622-637 (2001).

\bibitem{N03}
T. Nagao, 
Dynamical correlations for vicious random walk with a wall, 
{\it Nucl. Phys.}{\bf B658[FS]}:373-396 (2003).

\bibitem{NF98b}
T. Nagao and P. J. Forrester, 
Multilevel dynamical correlation function for Dyson's
Brownian motion model of random matrices,
{\it Phys. Lett.} {\bf A247}:42-46 (1998).

\bibitem{NF99}
T. Nagao and P. J. Forrester, 
Quaternion determinant expressions for multilevel
dynamical correlation functions of parametric
random matrices,
{\it Nucl. Phys.} {\bf B563[PM]}:547-572 (1999).

\bibitem{NKT03}
T. Nagao, M. Katori and H. Tanemura,
Dynamical correlations among vicious random walkers,
{\it Phys. Lett.} {\bf A 307}:29-35 (2003).

\bibitem{OR03}
A. Okounkov and N. Reshetikhin,
Correlation function of Schur process
with application to local geometry of 
a random 3-dimensional Young diagram,
{\it J. Amer. Math. Soc.} {\bf 16}:581-603 (2003).

\bibitem{OR07}
A. Okounkov and N. Reshetikhin,
Random skew plane partitions and the Pearcey process,
{\it Commun. Math. Phys.} {\bf 269}:571-609 (2007).

\bibitem{Osa96}
H. Osada,
Dirichlet form approach to infinite-dimensional 
Wiener processes with singular interactions,
{\it Commun. Math. Phys.} {\bf 176}:117-131 (1996).

\bibitem{Osa04}
H. Osada, 
Non-collision and collision properties of 
Dyson's model in infinite dimension and other stochastic dynamics
whose equilibrium states are determinantal random point fields, 
{\it Adv. Stud. in Pure Math.} {\bf 39}
``{\it Stochastic Analysis on Large Scale Interacting Systems}", 
pp.325-343, 
(Mathematical Society of Japan, Tokyo, 2004).

\bibitem{PM83}
A. Pandey and M. L. Mehta,
Gaussian ensembles of random Hermitian intermediate between
orthogonal and unitary ones,
{\it Commun. Math. Phys.} {\bf 87}:449-468 (1983).

\bibitem{PS02}
M. Pr\"ahofer and H. Spohn, 
Scale invariance of the PNG droplet and the Airy process,
{\it J. Stat. Phys.} {\bf 108}:1071-1106 (2002).

\bibitem{Rai00}
E. M. Rains, 
Correlation functions for symmetrized increasing
subsequences,
{\sf arXiv:math.CO/0006097}.

\bibitem{RY98}
D. Revuz and M. Yor,
Continuous Martingales and Brownian Motion. 3rd ed.,
(Springer, Now York, 1998).

\bibitem{RS93}
L. C. G. Rogers and Z. Shi,
Interacting Brownian particles
and the Wigner law,
{\it Probab. Th. Rel. Fields} {\bf 95}:555-570 (1993).

\bibitem{SI04}
T. Sasamoto and T. Imamura,
Fluctuations of the one-dimensional polynuclear growth model
in half space,
{\it J. Stat. Phys.} {\bf 115}:749-803 (2004).

\bibitem{Sel44}
A. Selberg,
Bemerkninger om et multiplet integral,
{\it Norsk Matematisk Tidsskrift} 
{\bf 26}:71-78 (1944).

\bibitem{ST03}
T. Shirai and Y. Takahashi,
Random point fields associated with certain
Fredholm determinants I:
fermion, Poisson and boson point process,
{\it J. Funct. Anal.} 
{\bf 205}:414-436 (2003).

\bibitem{Sos00}
A. Soshnikov, 
Determinantal random point fields,
{\it Russian Math. Surveys} {\bf 55}:923-975 (2000).

\bibitem{Sos03b}
A. Soshnikov, 
Janossy densities. II. Pfaffian ensemble,
{\it J. Stat. Phys.} {\bf 113}:611-622 (2003).

\bibitem{Sos04}
A. Soshnikov, 
Janossy densities of coupled random matrices,
{\it Commun. Math. Phys.} {\bf 251}:447-471 (2004).

\bibitem{Spo87}
H. Spohn H, 
Interacting Brownian particles:
a study of Dyson's model, in
{\it Hydrodynamic Behavior and Interacting Particle Systems}
(G. Papanicolaou Ed.), 
IMA Volumes in Mathematics and its Applications {\bf 9},
pp. 151-179,
(Springer, 1987).

\bibitem{Sta99}
R. P. Stanley,
{\it Enumerative Combinatorics},
vol.2,
(Cambridge University Press, Cambridge, 1999).

\bibitem{Ste90}
J. R. Stembridge,
Nonintersecting paths, pfaffians, and the plane partitions,
{\it Adv. in Math.} {\bf 83}:96-131 (1990).

\bibitem{Sze75}
G. Szeg\"o,
{\it Orthogonal Polynomials}, 4th ed.,
(American Mathematical Society, 1975).

\bibitem{TY03}
H. Tanemura and N. Yoshida,
Localization transition of $d$-friendly walkers,
{\it Probab. Th. Rel. Fields}, {\bf 125}:593-608 (2003).

\bibitem{Tit62}
E. C. Titchmarsh,
{\it Eigenfunction Expansions Associated with 
Second-order Differential Equations},
Part I, 2nd ed.,
(Clarendon Press, Oxford, 1962).

\bibitem{TW94}
C. A. Tracy and H. Widom,
Level-spacing distributions and the Airy kernel,
{\it Commun. Math. Phys.} {\bf 159}:151-174 (1994).

\bibitem{TW03}
C. A. Tracy and H. Widom,
A system of differential equations for the Airy process,
{\it Elect. Commun. Probab.} {\bf 8}:93-98 (2003).

\bibitem{TW04}
C. A. Tracy and H. Widom,
Differential equations for Dyson processes,
{\it Commun. Math. Phys.} {\bf 252}:7-41 (2004).

\bibitem{TW06}
C. A. Tracy and H. Widom,
The Pearcey process,
{\it Commun. Math. Phys.} {\bf 263}:381-400 (2006).

\bibitem{TW07}
C. A. Tracy and H. Widom,
Nonintersecting Brownian excursions,
{\it Ann. Appl. Probab.} {\bf 17}:953-979 (2007).

\bibitem{Zir96}
M. R. Zirnbauer,
Riemannian symmetric superspaces and their origin
in random-matrix theory,
{\it J. Math. Phys.} {\bf 37}:4986-5018 (1996).

\end{thebibliography}
\end{document}